\documentclass[11pt,openany]{book}
\usepackage{epsfig}

\newenvironment{proof}{{\it Proof:\/}}{$\Box$\vskip 0.08in}
\usepackage{amssymb}

 \newtheorem{theorem}{Theorem}[section]
 \newtheorem{lemma}[theorem]{Lemma}
\newtheorem{conjecture}[theorem]{Conjecture}
 \newtheorem{corollary}[theorem]{Corollary}
 \newtheorem{remark}[theorem]{Remark}
\newtheorem{proposition}[theorem]{Proposition}
\newtheorem{exercise}[theorem]{Exercise}

\newtheorem{definition}[theorem]{Definition}

\newtheorem{example}[theorem]{Example}

\newtheorem{formulla}[theorem]{}

\newcommand{\lk}{{\mbox{ lk }}}

\newcommand{\row}[2]{{\mbox{$#1_1,#1_2,\ldots,#1_{#2}$}}} 
\newcommand{\kwad}{\#}

\newcommand{\spn}{{\mbox{ span }}}

\newcommand{\pct}[1]{}%{{\Psfig{figure=#1}}}
\newcommand{\Psfig}[1]{{\mbox{$\ \ $}}}

\begin{document}
\renewcommand{\thechapter}{\Roman{chapter}}
\thispagestyle{empty}

%\maketitle

\
\vspace{0.5in}
 \begin{center}
 {\LARGE\bf KNOTS}\\
{\bf From combinatorics of knot diagrams to combinatorial topology
based on knots}
% Cambridge University Press, to appear, 2006, 2004
\end{center}

\vspace*{0.3in}

\centerline{Warszawa, November 30, 1984 -- Bethesda, October 31, 2004}
\vspace*{0.3in}

 \begin{center}
                      {\LARGE \bf J\'ozef H.~Przytycki}
\end{center}

\vspace*{0.3in}

%\newpage
\ \\
{\LARGE  LIST OF CHAPTERS}:\ \\
\ \\
{\LARGE \bf Chapter I: \ Preliminaries }\\
\ \\
{\LARGE \bf Chapter II:\ History of Knot Theory}\\
\ \\
{\LARGE \bf Chapter III:\ Conway type invariants }\\
\ \\ 
{\LARGE \bf Chapter IV:\  Goeritz and Seifert matrices}\\ \ \\
{\LARGE \bf Chapter V:\ Graphs and links}\\
{\bf This e-print. Chapter V starts at page 5}\\
\ \\
{\LARGE \bf Chapter VI:\ Fox $n$-colorings, Rational moves, Lagrangian tangles
and Burnside groups}\ \\
\ \\
{\LARGE \bf Chapter VII:\ Symmetries of links}\ \\
\ \\
{\LARGE \bf Chapter VIII:\ Different links with the same 
Jones type polynomials}\ \\
\ \\
{\LARGE \bf Chapter IX:\ Skein modules} \\
\ \\
{\LARGE \bf Chapter X:\ Khovanov Homology: categorification of the Kauffman
bracket relation}\\ 
{\bf e-print: http://arxiv.org/pdf/math.GT/0512630 }\\ \ \\
{\LARGE \bf Appendix I.\ }\ \\ 
\ \\
{\LARGE \bf Appendix II.\ }\\ \ \\
{\LARGE \bf Appendix III.\ }\\
%\centerline{\LARGE \bf LIST OF CHAPTERS GO HERE}
\
\newline

\ \\
{\LARGE \bf Introduction}\\
\ \\
This book is
about classical Knot Theory, that is, about
the position of a circle (a knot) or of a number of disjoint circles
(a link) in the space $R^3$ or in the sphere $S^3$.
We also venture into Knot Theory in general 3-dimensional
manifolds.

The book has its predecessor in Lecture Notes on Knot Theory,
which was published in Polish\footnote{The
Polish edition was prepared for the ``Knot Theory" mini-semester
at the Stefan Banach Center, Warsaw, Poland, July-August, 1995.}
in 1995  \cite{P-18}.
A rough translation of the Notes (by J.Wi\'sniewski) was
ready by the summer of 1995. It differed from the Polish edition
with the addition of
the full proof of Reidemeister's theorem. While I couldn't find
time to refine the translation and prepare the final manuscript,
I was adding new material and rewriting existing
chapters. In this way I created a new book based on the Polish 
Lecture Notes 
but expanded 3-fold.
Only the first part of Chapter III (formerly Chapter II),
on Conway's algebras is essentially unchanged from the Polish book
and is based on preprints \cite{P-1}.

As to the origin of the Lecture Notes, I was teaching an advanced course
in theory of 3-manifolds and Knot Theory at Warsaw University and it
was only natural to write down my talks (see Introduction to  (Polish)
Lecture Notes).
I wrote the proposal for the Lecture Notes by the  December 1, 1984
deadline.
In fact I had to stop for a while our work on generalization of
the Alexander-Conway and Jones polynomials in order to submit
the proposal. From that time several excellent books on Knot Theory
have been published on various level and for various readership.
This is reflected in my choice of material for the book -- knot theory
is too broad to cover every aspect in one volume. I decided to
concentrate on topics on which I was/am doing an active research.
Even with this choice the full account of skein module theory is
relegated to a separate book (but broad outline is given in
Chapter IX).

In the first Chapter we offer historical perspective to
the mathematical theory of knots, starting from the first
 precise approach to Knot Theory by Max Dehn
and Poul Heegaard in the Mathematical Encyclopedia \cite{D-H} 1907.
We start the chapter by introducing
lattice knots and polygonal knots.
The main part of the chapter is devoted to the proof of
Reidemeister's theorem which allows combinatorial treatment of
Knot Theory.

In the second Chapter we offer the history of Knot 
Theory starting from the ancient Greek tract on surgeon's slings, 
through Heegaard's thesis
relating knots with the field of {\it analysis situs}
(modern algebraic topology) newly developed by Poincar\`e,
 and ending with the Jones polynomial and related knot invariants.

In the third Chapter we discuss invariants of Conway type; that is,
invariants which have the following property: the values of the invariant
for oriented links $L_0$ and $L_-$  determine
its value for the link $L_+$ (similarly, the values of the invariant
for $L_0$ and $L_+$ determine its value for $L_-$).
The diagrams of oriented links $L_0$, $L_-$ and $L_+$ are different
only at small disks as pictured in Fig.~0.1.\ \\
%\pagebreak
%\vspace*{0.8in}
\centerline{\psfig{figure=L+L-L0.eps,height=1.8cm}}
\begin{center}
                    Fig.~0.1
\end{center}

Some classical invariants of knots turn out to be invariants
of Conway type.  
These include the number of components,
the global linking number, the normalized Alexander polynomial
(Conway polynomial), the signature, the Jones polynomial, and
its 2-variable generalization known as the Jones-Conway or
Homflypt polynomial\footnote{Actually, it seems that there is
no fixed name for this invariant; the following names are also used:
Conway-Jones, Flypmoth, Homfly, skein, Thomflyp, twisted Alexander,
generalized Jones, two variables Jones.}.
Sikora's proof that Conway algebras do not give any invariants
of links, stronger than  Jones-Conway polynomial, is given
in  
%Section XX of 
Chapter 3.
In the second part we will also discuss generalizations
of Conway type invariants which are obtained by adding an extra diagram
$L_\infty$, see Fig.~0.3.
We will also discuss Kauffman's method of constructing
invariants of links.\ \\
\centerline{\psfig{figure=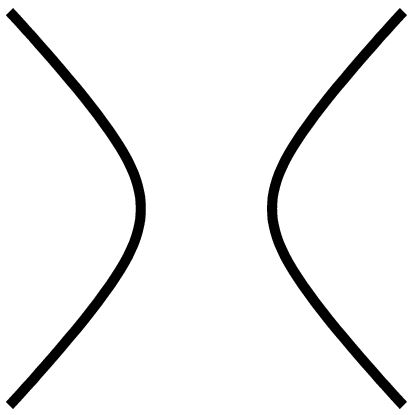,height=1.3cm}}
\begin{center}
                              Fig.~0.3
\end{center}
In the fourth Chapter of the notes we will describe
Goeritz and Seifert matrices and their relations to the
Jones type link invariants.
%Maybe here Traczyk's formula for signature of alternating links?
%At this point we will apply some elementary techniques
%of algebraic topology.
%A reader interested only in combinatorial aspects of the theory
%may want to consider merely browsing this part or skipping it entirely.

In the fifth Chapter we present applications
of graphs in Knot Theory and we prove two classical conjectures of Tait
by using the Jones and Kauffman polynomials.
%Is second Tait conjecture, and signature, following from Jones of
%2-cabling of alternating link?
Moreover, we discuss two important classes
of links: alternating links and
their generalizations, adequate links.

In the sixth Chapter we discuss several
open problems in classical Knot Theory and we develop techniques
that allow us to study them:  Lagrangian tangles
and Burnside groups.

In the seventh Chapter we examine symmetries of links.
New polynomial invariants provide us with very efficient criteria
for studying symmetric links.
As an application we give a partial
characterization of knots which are obtained from trajectories
of a point in a 3-dimensional billiard.

In the eighth Chapter we analyze various methods of constructing
different links with the same Jones type polynomials.
We demonstrate how ideas from the Graph Theory and statistical mechanics are
fruitfully applied in the Knot Theory.

In the ninth Chapter, we propose a generalization of Jones-type invariants
to any 3-dimensional manifold via a construction of skein modules.
Our method leads to {\it algebraic topology based
on knots}\footnote{During the first half of the XX century the branch
of topology which is now named {\it algebraic topology} was called
{\it combinatorial topology}. This name motivates the subtitle of
this book.}. We sketch the theory in this chapter and the full account
will be described in the sequel book that is under preparation
\cite{P-30}. %\cite{Pr-??}

In the tenth Chapter, we describe Khovanov homology of links in $S^3$.
We study the size (thickness) of them and their torsion part.
Subsequently we describe generalization of Khovanov homology to
a 3-manifold being an $I$-bundle over a surface. In this case we
relate Khovanov homology with the Kauffman bracket skein module discussed
in the ninth Chapter.

The book is supplemented with three appendices...
\ \  \ \ SEE Introduction before CHAPTER I.

\setcounter{chapter}{4}

\chapter{Graphs and links}\label{V}
%Jesli section lub chapter ma za dlugi tytul to zaraz po
%\chapter{...} albo \section{...} piszesz
%\markboth{\hfil{\sc ...}\hfil}
%{\hfil{\sc ...}\hfil} I wypelniasz kropkowane miejsca odpowiednimi 
%skrotami nazwy section i chapter.
\centerline{Bethesda, October 31, 2004}
In Chapter V we present several results which demonstrate a close 
connection and useful exchange of ideas between graph theory and 
knot theory. These disciplines were shown to be related from the 
time of Tait (if not Listing) but the great flow of ideas started 
only after Jones discoveries. The first deep relation in this new 
trend was demonstrated by Morwen Thistlethwaite and we describe several 
results by him in this Chapter. We also present 
%We present, in this e-print, 
results from two preprints \cite{P-P-0,P-34}, 
in particular we sketch 
two generalizations of the Tutte polynomial of graphs,
 $\chi(G;x,y)$, or, more precisely, the deletion-contraction method 
which Tutte polynomial utilize.
The first generalization considers, instead of graphs, general objects 
called setoids or group systems. The  second one deals with completion 
of the expansion of a graph with respect to subgraphs.
 We are motivated here by finite type 
invariants of links developed by Vassiliev and Gusarov along the line 
presented in \cite{P-9} (compare Chapter IX).
The dichromatic Hopf algebra, described in Section 2, have its 
origin in Vassiliev-Gusarov theory mixed with work of G. Carlo-Rota and 
his former student (now professor at GWU) W. Schmitt.
%This e-print is a preliminary version of the first part of Chapter 5 
%of the book ({\it KNOTS, 
% from combinatorics of knot diagrams to combinatorial topology
%based on knots}) which I am preparing for Cambridge University Press.

\section{Knots, graphs and their polynomials}\label{V.1}
%\centerline{Preliminary version of the translation for Cambridge Univ.Press.}
%\markboth{\hfil{\sc Graphs and links}\hfil} {\hfil{\sc Knots}\hfil}

In this section we discuss relations between graph and knot theories.
We describe several applications of graphs to knots.
In particular we consider various  
interpretations of the Tutte polynomial of graphs in knot theory.
This serves as an introduction to the subsequent sections
where we prove two of the classical conjectures of Tait \cite{Ta}.
In the present section we rely mostly on 
\cite{This-1,This-5} and \cite{P-P-1}.

By a graph $G$ we understand a finite set $V(G)$ of vertices
together with a finite set of edges $E(G)$. To any edge we associate 
a pair of (not necessarily distinct) vertices which we call 
endpoints of the edge.
We allow that the graph $G$ has multiple edges and loops
(Fig.1.1)\footnote{In terms of algebraic topology a graph is 
a 1-dimensional CW-complex. Often it is called a pseudograph and the 
word ``graph" is reserved for a 1-dimensional simplicial complex, 
that is,  loops and multiple edges are not allowed. 
We will use in such a case the 
term a {\it simple} (or {\it classical}) graph. 
If multiple edges are allowed but loops 
are not we use often the term a {\it multigraph}, \cite{Bo-1}}.
 A loop is an edge with one endpoint.

\vspace*{1.5in}\centerline{\Psfig{figure=Rys.1.1}}
\centerline{{\psfig{figure=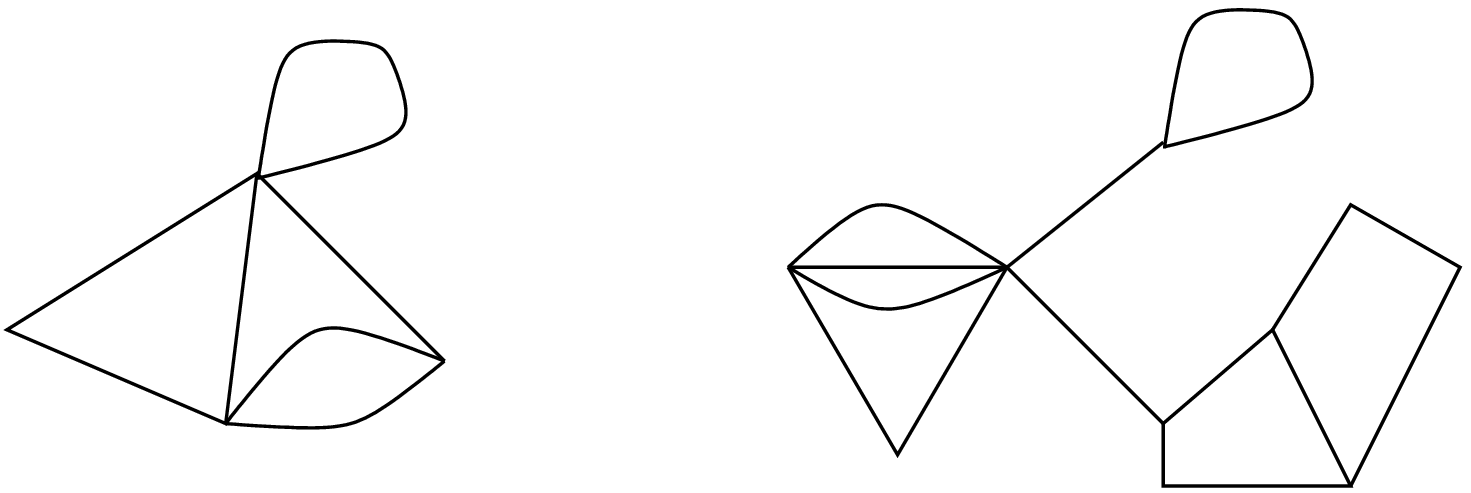,height=2.5cm}}}  
\begin{center}
Fig.~1.1
\end{center}

By $p_0(G)$ we denote the number of components of the graph $G$
and by $p_1(G)$ we denote its cyclomatic number, i.e.~the minimal number
of edges which have to be removed
from the graph in order to get a graph without 
cycles.\footnote{In terms of algebraic 
topology $p_0(G)$ and $p_1(G)$ are equal to
dimensions of homology groups $H_0(G)$ and $H_1(G)$, respectively. In
this context the notation $b_0$ and $b_1$ is used and numbers are called
the Betti numbers.}
A connected graph without cycles (i.e.~$p_0=1$, $p_1=0$) is called a tree.
If $G$ has no cycles , i.e.~$p_1=0$, then the graph $G$ is called a forest.
By a spanning tree (resp. forest) of the graph $G$ we understand 
a tree (resp. forest) in $G$ which contains all vertices of $G$. 
By an isthmus of $G$ we understand an edge of $G$, 
removal of which increases the number of components of the graph. 

To a given graph we can associate a polynomial in various ways.
The first such a polynomial, called the chromatic polynomial of a graph,
was introduced by Birkhoff in 1912 \cite{Birk}\footnote
{J.B.Listing, in 1847\cite{Lis}, introduced polynomial of knot diagrams.
For a graph $G$, the Listing polynomial, denoted by $JBL(G)$,
 can be interpreted as follows:
$JBL(G)=\Sigma a_i(G)x^i$ where $a_i(G)$ is the number of vertices in 
$G$ of valency $i$.}.
For a natural number $\lambda$, the chromatic polynomial, denoted 
by $C(G,\lambda)$, counts the number 
of possible ways of coloring the vertices of $G$ in $\lambda$ colors
in such a way that each edge has endpoints colored in different colors
(compare Exercise 1.14).
The chromatic polynomial was generalized 
by Whitney and Tutte \cite{Tut-1}.

%\begin{definition}\label{V.1.1}
\begin{definition}\label{4:1.1}

The following conditions define the Tutte polynomial\footnote{Tutte called
this polynomial the {\it dichromat}.} 
$\chi (G;x,y)\in Z[x,y]$ of a graph $G$: 

\begin{enumerate}
\setcounter{enumi}{-1}

\item[(1)]
 $\chi (\bullet) = \chi (\emptyset ) = 1$

\item[(2)]
 $\chi(\bullet\!\!\!\!\longrightarrow\!\!\!\!\bullet) = x$

\item[(3)]
 $\chi (\bullet\!\!\bigcirc) = y$

\item[(4)]
%$\chi(G_1 \!\!\!\!\bullet G_2)$
$\chi (G_1 * G_2) = \chi (G_1) \chi (G_2)$,
where the product $G_1 * G_2$ is obtained from $G_1$ and $G_2$
by identifying two vertices, each one chosen on each of the two 
graphs\footnote{The product $*$ depends on the choice of base point 
vertices which are identified. The precise notation should be 
$(G_1,v_1)*(G_2,v_2)$.}.

\item[(5)]
 $\chi (G_1\sqcup G_2) = \chi (G_1)\chi (G_2)$, where $\sqcup$
denotes the disjoint sum of graphs.

\item[(6)]
 $\chi (G) = \chi (G-e) + \chi (G/e)$, where $e$ 
is an edge which is neither a loop nor an isthmus
and $G/e$ denotes contracting of the edge $e$, i.e.~a graph
which is obtained from $G$ by removing $e$
and identifying its endpoints.
\end{enumerate}
\end{definition}

Before we show that Tutte polynomial is well defined we 
suggest the following exercise.

\begin{exercise}\label{V.1.2}
Prove that $\chi (T_{i,j})=x^iy^j$, where $T_{i,j}$ is a connected graph
obtained from a tree of $i$ edges by adding $j$ loops to it. 
\end{exercise}
In this exercise we use the Euler's lemma that every tree has a vertex 
of degree 1, where degree (or valency) of a vertex is the number of 
%adjacent 
incident edges (counting a loop twice).

In order to prove the existence of the Tutte 
polynomial\footnote{Impatient reader can prove existence quickly by 
first ordering edges of $G$ and then using formula (6) for edges, 
in chosen ordering, till one reaches trees with loops for which 
the formula from Exercise 1.2 is applied. 
Then one checks that changing ordering of
edges preserve the polynomial; compare Fig. 1.4.}  
we will consider 
a slightly more general polynomial invariant of graphs, which is closely
related to link polynomials. 
Namely, we will define the Kauffman bracket polynomial 
and we will compare it with the Tutte polynomial.

\begin{definition}\label{V.1.3}
The Kauffman bracket polynomial $\langle G\rangle$ of the graph 
$G$ ($\langle G\rangle\in Z[\mu,A,B]$) 
is defined inductively by the following formulas:

\begin{enumerate}
\item[(1)]
 $\langle\bullet\rangle = 1$
\item[(2)]
 $\langle G_1\sqcup G_2\rangle = \mu\langle G_1\rangle\langle
G_2\rangle$

\item [(3)]
$\langle G\rangle = B\langle G-e\rangle + A\langle G//e\rangle$
where $G//e = G/e$ if $e$ is not a loop, and if $e$ is a loop, then 
by $G//e$ we understand a graph with the edge $e$ removed and one 
``free'' vertex added.
\end{enumerate}
\end{definition}

The Kauffman bracket polynomial of a graph is uniquely defined as our
rules allows computation of a polynomial for every graph. It is well 
defined because it can be given by a single formula which satisfies
our rules.

\begin{lemma}\label{4:1.4}
$$\langle G\rangle = \sum_{S \in 2^{E(G)}}
\mu^{p_0(G-S)+p_1(G-S)-1} A^{|E(G) - S|} B^{|S|}$$
where $S$ is an arbitrary set of edges of $G$, including the empty
set, and $G-S$
denotes a graph obtained from $G$ by removing
all these edges.
\end{lemma}

\begin{exercise}\label{4:1.5}
Prove that the formula for $\langle G\rangle$ introduced
in Lemma 1.4 satisfies all conditions which are set up
in the Definition 1.3. In particular, show that
if $G$ is a tree with loops 
then $\langle G\rangle = (A+B\mu)^a (B+\mu A)^b$, 
where $a$ is the number of edges in the tree and 
$b$ is the number of loops.

\end{exercise}

\begin{theorem}\label{4:1.6}
The following identity holds
$$\langle G\rangle = 
\mu^{p_0 (G)-1} B^{p_1 (G)} A^{E(G)-p_1(G)} \chi (G;x,y)$$ 
where $x = \frac{A+\mu B}{A}$ and $ y = \frac{B+\mu A}{B}$.
\end{theorem}

Proof. Using properties of the Kauffman bracket
polynomial, one can verify 
easily that $\chi(G,x,y)$ computed from the theorem satisfies
the conditions of Definition 1.1 (c.f.~\cite{P-P-1}).
In particular if $e$ is an edge of $G$ which is neither
an isthmus nor a loop and assuming that the theorem holds
for $G-e$ and $G/e$, we obtain: 
$\mu^{p_0 (G)-1} B^{p_1 (G)} A^{E(G)-p_1(G)} \chi (G;x,y)=
\langle G\rangle = B\langle G-e \rangle + A\langle G/e \rangle =$ \\
$B(\mu^{p_0 (G-e)-1} B^{p_1 (G-e)} A^{E(G-e)-p_1(G-e)}) \chi (G-e);,x,y)+$\\
$A(\mu^{p_0 (G/e)-1} B^{p_1 (G/e)} A^{E(G/e)-p_1(G/e)}) \chi (G/e;x,y)=$\\
$\mu^{p_0 (G)-1} B^{p_1 (G)} A^{E(G)-p_1(G)} (\chi (G-e;x,y) 
+ \chi (G/e;x,y)$.

%$B(\mu^{p_0 (G-e)-1} B^{p_1 (G-e)} A^{E(G-e)-p_1(G-e)}) \chi (G-e);,x,y)+\\
%A(\mu^{p_0 (G/e)-1} B^{p_1 (G/e)} A^{E(G/e)-p_1(G/e)}) \chi (G/e;x,y)=\\
%B\langle G-e \rangle + A\langle G/e \rangle = \langle G\rangle$.

There are some simple but very useful properties of Tutte 
polynomial which follow quickly from our definition and basic 
properties of 2-connected graphs.
\begin{definition}\label{V.1.7}
\begin{enumerate}
\item[(i)]
We say that a graph $G$ is
2-connected if it is connected and has no cut vertex, i.e. $G$ cannot
be expressed as $G_1*G_2$ with $G_i$ having more than one vertex or 
being a loop.
\item[(ii)] More generally we say that a graph $G$ is $n$-connected 
if it is $(n-1)$-connected and cannot be obtained from two graphs 
$G_1$ and $G_2$, each of at least $n$ vertices 
by gluing them together along $n-1$ vertices.
\end{enumerate}
\end{definition}
\begin{lemma}\label{V.1.8}
Let $e$ be any edge of a 2-connected graph $G$ then
\begin{enumerate}
\item[(i)] if $G$ has more than one edge then $G-e$ 
and $G/e$  are connected. 
\item[(ii)] Either $G-e$ or $G/e$ is $2$-connected\footnote{I have 
been informed by Robin Thomas that analogous
theorem holds for 3-connected graphs:\ Every 3-connected
graph $G$ on at least five vertices has
an edge $e$ such that the graph $G-e$ or $G/e$ is $3$-connected. There is
similar theorem for $4$-connected graphs but nothing is known for
$n>4$.}.
\item[(iii)] Let $H$ be any 2-connected subgraph of $G$ then 
one can obtain $H$ from $G$ by a sequence of deletions and contractions 
in such a way that every graph on the way between $G$ and $H$ is 
2-connected. 
\item[(iv)] If $H$ is a minor of $G$ that is $H$ can be obtained from 
$G$ by a sequence of deletions and contractions and $H$ is 2-connected 
then we can find such a sequence so that every graph on the way is 
2-connected.
\end{enumerate}
\end{lemma}
\begin{proof}
\begin{enumerate}
\item[(i)] If $G$ has more than one edge and $G-e$ was a disjoint sum 
of $G_1$ (which is not one vertex graph) and $G_2'$, then we take 
$G_2$ obtained from $G_2'$ by adding $e$ to it. 
Then $G= G_1*G_2$, the contradiction.
\item[(ii)] It holds for 1-edge graph so let assume that $G$ has at least 
two edges. Let us assume now that $G-e$ is not 2-connected and that $v$
is a vertex the removal of
which makes $G-e$ disconnected; see Fig.~1.2.
Note that  $v$ cannot be an endpoint of $e$.
Let $v_e$ be a vertex of $G/e$ obtained from
endpoints of $e$. Clearly $v_e$ cannot be a cut vertex of $G/e$, Fig.~1.2.
On the other hand $v_e$ is the only possible vertex which can be
a cut vertex of $G/e$ (any other cut vertex of $G/e$ would be also 
a cut vertex of $G$). Thus $G/e$ is 2-connected.
\ \\
\centerline{{\psfig{figure=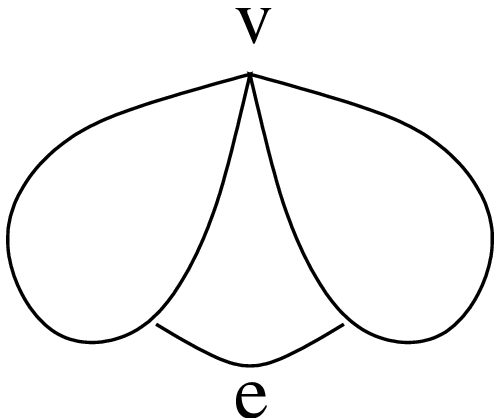,height=2.1cm}}\ \
{\psfig{figure=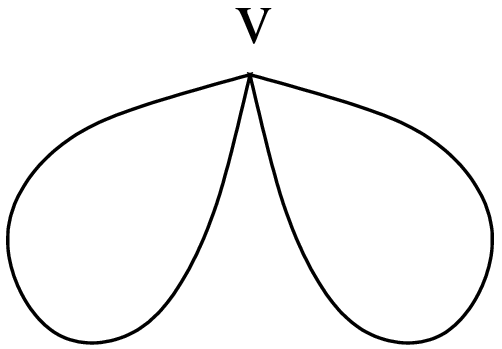,height=1.9cm}}\ \
{\psfig{figure=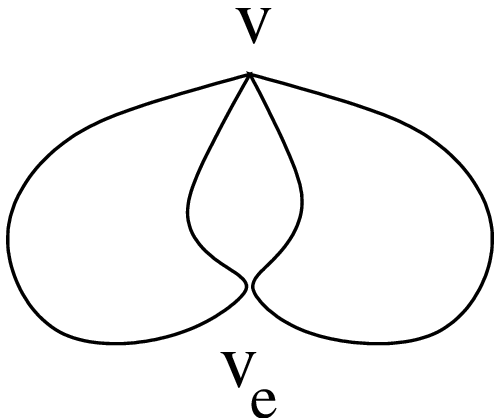,height=2.1cm}}
}
\begin{center}
Fig.~1.2.\ \  Graphs $G$, $G-e$ and $G/e$
\end{center}

\item[(iii)] We proceed by induction on the number of edges $E(G)-E(H)$.
Of course (ii) holds for $H=G$ so assume that $e$ is an edge in $G$ but 
not in $H$. If at least one vertex of $e$ is not in in $H$ then 
$H$ is a subgraph of $G-e$ and $G/e$  and we use an inductive assumption 
for that one which is 2-connected. If every edge in $E(G)-E(H)$ has 
 both endpoints on $H$ than deleting any edge of $E(G)-E(H)$ gives 
2-connected graph. One can  to visualize it by observing that 
adding an edge, which is not a loop, to a 2-connected graph ($H$ in our case)
 leads to a 2-connected graph. 
\item[(iv)] We modify inductive proof given in (iii) to this more general 
situation. As before assume that $e$ is an edge in $G$ but
not an edge of $H$. If at least one endpoint of $e$ is not in $H$ (one 
vertex in $H$ can correspond to several vertices in $G$) then $H$ is a 
minor of $G-e$ and $G/e$  and we use an inductive assumption
for that one which is 2-connected. If every edge in $E(G)-E(H)$ has
 both endpoints in $H$ and this endpoints are identified in $H$ then 
we contract this edge if $G/e$ is 2-connected. If $G/e$ is not 2-connected 
then $G$ can be decomposed into $G_1 \cup e \cup G_2$ as shown in Fig. 1.3 
with $G_1 \cup G_2$, $G_1 \cup e$ and $e \cup G_2$ 2-connected. \\

\ \\
\centerline{\psfig{figure=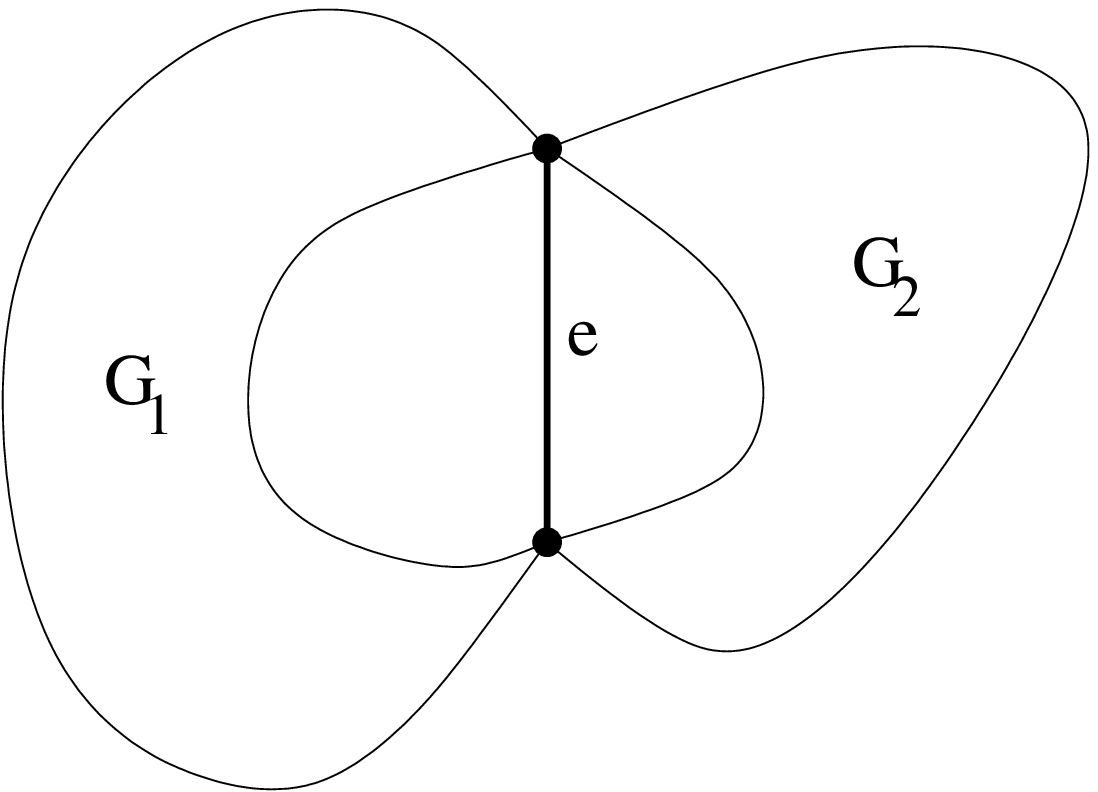,height=3.1cm}}
\begin{center}
Fig. 1.3.\ \  $G$ is 2-connected but $G/e$ is not
\end{center}
\ \\
Because 
endpoints of $e$ are identified in $H$ and $H$ is 2-connected therefore 
whole graph $G_i\cup e$ is minored (by deleting and contracting) to 
a point in $H$ for $i=1$ or $2$. Assume that it holds for $i=2$. 
Then we can use part (iii) of the lemma (or just inductive assumption)
to the 2-connected subgraph $G_1\cup e$ of $G$. Thus we can reach 
$G_1\cup e$ from $G$ via 2-connected graphs and $H$ is a minor of $G_1\cup e$ 
so we can use inductive assumption once more. Finally assume that every 
edge in $E(G)-E(H)$ has both endpoints on $H$ and this endpoints are
different in $H$. Therefore $H$ is a subgraph of $G$ and we can just 
delete these edges one by one (as in (iii)).
\end{enumerate}
\end{proof}

\begin{corollary}\label{V.1.9}
\ \\ Let $\chi (G;x,y)= \Sigma v_{ij}x^iy^j$. Then 
\begin{enumerate}
\item[(i)] $v_{ij}\geq 0$ and $v_{0,0}=0$ iff $|E(G)|>0$.

\item [(ii)] Let us assume that  $G$ is a 2-connected graph
with at least two edges,
in particular $G$ has neither a loop nor an isthmus. Then 
$v_{0,1} = v_{1,0} > 0$.

\item [(iii)] (a)If $G$ is a 2-connected graph with at least three
vertices then $v_{2,0}>0$.\\
(b) If $G$ is a 2-connected graph with at least three edges then 
$v_{0,2} + v_{2,0} - v_{1,1} > 0$.
\item[(iv)] If $G$ is a 2-connected graph which is neither
an $n$-gon nor 
a generalized theta curve\footnote{The generalized theta curve
 is dual to a polygon; compare Theorem 1.13.} (two vertices, connected
by $n$ edges) then $v_{1,1}>0$.
\item [(v)]
 If the graph $G$ has $\alpha$ isthmuses and $\beta$ loops
then $$\chi (G;x,y) = x^\alpha y^\beta \chi (G_1;x,y)\chi
(G_2;x,y),\ldots,\chi (G_t;x,y),$$
where $G_i$ are 2-connected components of $G$ with more than one edge.
\end{enumerate}
\end{corollary}

We extend our proposition in Exercise 1.11.

\begin{proof}
\begin{enumerate}
\item[(i)]
It follows from the definition of the Tutte polynomial (Def. 1.1).
\item[(ii)]
 We apply induction with respect to the number of edges
in the graph. We start with a graph 
%$\bullet\!\!\!\!\bigcirc\!\!\!\!\bullet$
\parbox{0.7cm}{\psfig{figure=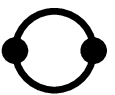,height=0.5cm}},
for which Corollary 1.9(ii) holds, that is
$\chi(
%\bullet\ \!\!\!\!\bigcirc\!\!\!\!\ \bullet\ 
\parbox{0.7cm}{\psfig{figure=2gon.eps,height=0.5cm}}
; x,y) = x+y$, 
and thus $v_{0,1} =v_{1,0} = 1$. 
Now let $G$ be an arbitrary 2-connected graph which has
$|E(G)|>2$ edges and we assume that for graphs with a smaller number
of edges the property 1.9 (ii) is true. Let $e$
be an arbitrary edge of $G$. Since $e$ is neither a loop nor an isthmus 
it follows that  $\chi (G) = \chi (G-e)+\chi (G/e)$. 
%%%%%%%%%%%%%%%%%%%%%%%%%%%%%%%%%%%%%%%%%%%%%%
%%%%%%%%%%%%%%%%%%%%%%%%%%%%%%%%%%%%%%%%%%%%%%%
Now, to prove that $v_{0,1},v_{1,0}>0$ we use 
Lemma 1.8 (either $G-e$ or $G/e$ is 2-connected) 
and the inductive assumption. To see that $v_{0,1}=v_{0,1}$ we 
use additionally the fact that if a graph $G-e$ or $G/e$ is not 
2-connected than it has an isthmus or a loop and then 
$v_{0,1}=v_{1,0}=0$.
%(it follows from (iii) that 
%for any connected but not 2-connected graph with at least 2 edges we have 
%$v_{0,1}=v_{1,0}=0$).

We can reformulate the idea of our proof in a more sophisticated manner 
by saying that we proved that every 2-connected graph $G$ with at least 
two edges has \ 
$\bullet\!\!\!\!\bigcirc\!\!\!\!\bullet$ \ 
as its 
minor (in the class of 2-connected graphs). 
%That is, $\bullet\!\!\!\!\bigcirc\!\!\!\!\bullet$ can be reached from
%$G$ by a sequence of deletions and contractions (in the class of 
%2-connected graphs). 
\item[(iii)] (a) If $G$ has $n$-gon as a subgraph $(n>2)$ then 
$\chi (G;x,y)$
contains as a summand the Tutte polynomial of the $n$-gon, that is
$x^{n-1}+...+x^2+x +y$ and $v_{2,0}>0$. Otherwise $G$ 
is a generalized theta curve  of $n>2$ edges (that is it has 2 vertices
connected by $n$ edges) and then 
$(G;x,y) = x+y+y^2+...+y^{n-1}$ so $v_{2,0}=0$.\\
(b) The formula $v_{0,2} + v_{2,0} - v_{1,1}=v_{1,0}$ 
holds for any graph with at least 3 edges. We check 
first that it holds for a tree with loops and for any graph with exactly 
3 edges. Then one induct on the number of edges using Tutte formula,  
Definition 1.1(6) (see \cite{Bo-2}, Exercise X.7.8 and its generalization 
by T.H.Brylawski\footnote{Brylawski's formula says that for a graph $G$ with 
more than $h$ edges we have the identity
$$\sum_{i=0}^h\sum_{j=0}^{h-i}(-1)^j{{h-i}\choose{j}}v_{i,j} = 0.$$}).
\item[(iv)] If $G$ is neither an $n$-gon nor the generalized theta 
curve then $G$ has the graph 
\parbox{0.7cm}{\psfig{figure=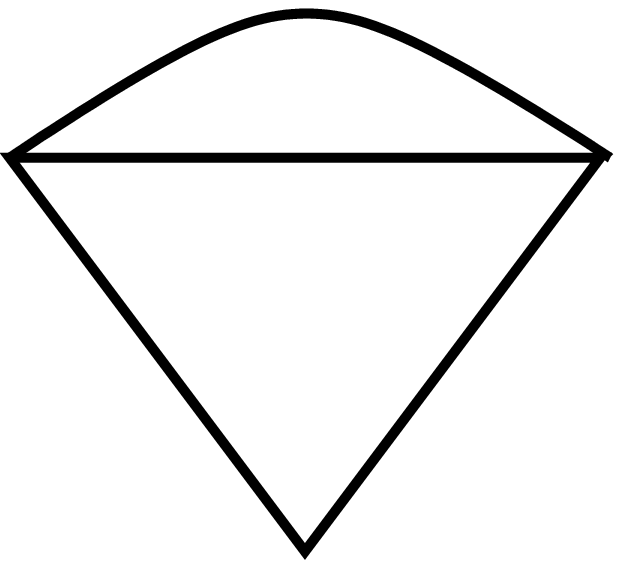,height=0.6cm}}
 as its minor (in fact $G$ contains an n-gon, ($n >2$, with two 
of its vertices connected by a path outside a polygon,
as a subgraph. Now $\chi 
(\parbox{0.7cm}{\psfig{figure=triangle2.eps,height=0.6cm}};x,y)= 
x^2+x+xy+y+y^2$ and we use 
Corollary (iii) and Tutte formula (6) of Definition 1.1 to complete 
the proof of (iv).
\item[(v)] It follows from Properties (2-4) of Definition 1.1.
\end{enumerate}
\end{proof} 
\begin{corollary}\label{V.1.10}
The numbers $v_{0,1}$ and $v_{1,0}$ are topological invariants
of the graph. That is, if a graph $G$ has at least two edges and $e$ is 
its edge then the subdivision of $e$ 
$(\stackrel{e}{\bullet\!\!\!\!\longrightarrow\!\!\!\!\bullet}\rightarrow
\stackrel{e_1\ \ \
e_2}{\bullet\!\!\!\!\longrightarrow\!\!\!\!\bullet\!\!\!\!\longrightarrow\!\!\!\!\bullet})$
changes neither $v_{0,1}$ nor $v_{1,0}$. 
\end{corollary}
\begin{proof}
 Let $G'$ be a graph obtained from $G$ by dividing the edge
$e$ into $e_1$ and $e_2$. For any graph with an edge, $v_{0,0} =0$. Now
if $e$ is an isthmus of $G$, then 
$\chi (G) = x \chi (G-e)$, thus $v_{1,0}(G) = v_{0,0}(G-e) = 0$ and
$v_{0,1}(G)=0$. Similarly $G'$ has an isthmus, so $v_{1,0}(G')=
v_{0,1}(G')= 0$. 
We can also give more general argument, based on the formula of
Corollary 1.8(iii), that is $G$ has at least two edges and is
connected but not 2-connected, then $v_{1,0}(G) =v_{0,1}(G)= 0$.

Assume now that $e$ is not an isthmus of $G$. Then $e_1$ is
neither an isthmus nor a loop of $G'$.  
Thus $\chi (G') = \chi (G' - e_1) +\chi (G'/e_1)$ but $G'/e_1
= G$ and $G' - e_1$ has an isthmus and therefore
$v_{0,1} (G' - e_1) = v_{1,0}(G' - e_1) = 0$. 
Finally $v_{1,0}(G')= v_{1,0}(G)$ and $v_{0,1} (G') =v_{0,1} (G)$.
This concludes the proof of Corollary 1.10. \ If $G$ has only one 
edge, $e$, then $e$ is either an isthmus  
and $\chi (G;x,y)=x$, $\chi (G';x,y)=x^2$,
or $e$ is a loop and $\chi (G;x,y)=y$, $\chi (G';x,y)=x+y$ and $v_{1,0}$  
is changed.
\end{proof}

\begin{exercise}\label{V.1.11}
Let $G$ be a connected graph with cyclomatic number  equal to $p_1(G)$,
and with $d(G)$ edges in every spanning tree 
($d(G)= |E(G)| - p_1(G) = |V(G)| -1$).
Show that:
\begin{enumerate}
\item[(1)] If $p$ is the number of loops in $G$ then
$v_{d,p} = 1$. Furthermore $v_{i,j}=0$ if $i> d(G)$.
\item[(2)] If $s$ is the number of isthmuses in $G$ then
$v_{s,p_1(G)} = 1$. Furthermore $v_{i,j}=0$ if $j> p_1(G)$.
\item[(3)] 
\begin{enumerate}
\item[(a)] 
If $x^iy^j$ is the maximal degree monomial dividing $\chi (G)$
then $G$ has $i$ isthmuses and $j$ loops.
\item[(b)]
The numbers $p_1(G),d(G),|E(G)|$ and $|V(G)|$ are determined by 
$\chi (G)$.
\end{enumerate} 
\item[(4)] Let  $G$ be a 2-connected graph with at least two edges. Then:\\  
 $v_{0,j} > 0$ if and only if $1\leq j\leq p_1(G)$ and
$v_{i,0} > 0$ if and only if $1\leq i\leq d(G)$. In particular 
$\chi (G;x,y)$
contains the summand $x^{d(G)}+...+x+y+...+y^{p_1(G)}$.
\item[(5)] If $G$ is a 3-connected graph of at least 4 vertices
then $\chi (G;x,y)$ contains as a summand the Tutte polynomial
of the complete graph on 4 vertices, $\chi (K_4;x,y) =
x^3+3x^2+2x + 4xy + 2y + 3y^2+y^3$. 
\item[(6)] Formulate analogue of part (5) for $4$- and $5$-connected 
planar graphs knowing that every $4$-connected planar graph with at least 
$5$ vertices has the octahedral graph as its minor (Fig. 1.9) and 
that every $5$-connected planar graph with at least 
$6$ vertices has the icosahedral graph as its minor \cite{Bo-1}. 
\end{enumerate}
\end{exercise}

Hint.\ The crucial fact we use in the inductive proof of Part (4) 
is Corollary 1.9(ii) ($G-e$ or $G/e$ is 2-connected).
In Part (5) we should show first that $K_4$ is a minor of every 3-connected 
graph with at least 4 vertices.

\begin{exercise}\label{V.1.12} Let $(G_1)^*_*(G_2)$ denote the 2 vertex 
product of graphs, that is we choose 2 vertices $v_i,w_i$ on $G_i$, 
$i=1,2$ and identify $v_1$ with $v_2$ and $w_1$ with $w_2$ (in full 
notation $(G_1,v_1,w_1)^*_*(G_2,v_2,w_2)$).
\begin{enumerate}
\item[(i)]  Find the formula for the Kauffman bracket $<(G_1)^*_*(G_2)>$ 
when \\ 
$<G_1>$, $<G^d_1>$ , $<G_2>$, $<G^d_2>$ are given. Here $G_i^d$ is 
the graph obtained from $G_i$ by identifying $v_i$ with $w_i$.
\item[(ii)] Show that the Kauffman bracket polynomial of $(G_1)^*_*(G_2)$ 
does not depend on the ordering of identified vertices, that is 
$<(G_1,v_1,w_1)^*_*(G_2,v_2,w_2)> $ \\ 
$= <(G_1,v_1,w_1)^*_*(G_2,w_2,v_2)>$. Borrowing terminology from Knot 
Theory we say that the second graph is obtained from the first by 
{\it mutation}\footnote{The term {\it Whitney twist} is occasionally used 
in graph theory but 
sometimes it means the operation which keeps the abstract graph and 
changes only its plane embedding.}
 and it is called the mutant of the first graph (Fig. XX present 
a pair of mutant graphs). 
\item[(iii)] Show that if $G_1$ and $G_2$ are 2-connected graphs and 
$v_1\neq w_1$, and $v_2\neq w_2$  then 
$(G_1,v_1,w_1)^*_*(G_2,v_2,w_2)$) is a 2-connected graph. 
\end{enumerate}
\end{exercise}

Below we outline the underlining ideas of the Tutte work on polynomial $\chi$ 
and relations to Knot Theory
(following \cite{P-P-0}).

Order edges of $G$:\ $e_1,e_2,...,e_E$. To find $\chi (G;x,y)$ we apply
deleting-contracting formula to edges of $G$ one by one (according to
our ordering) and never using an isthmus or a loop. Our computation
can be summarized by a binary computational tree, whose leaves are
trees with loops (as $G$ is connected). See figure below. 

%\vspace*{1.5in}
\centerline{{\psfig{figure=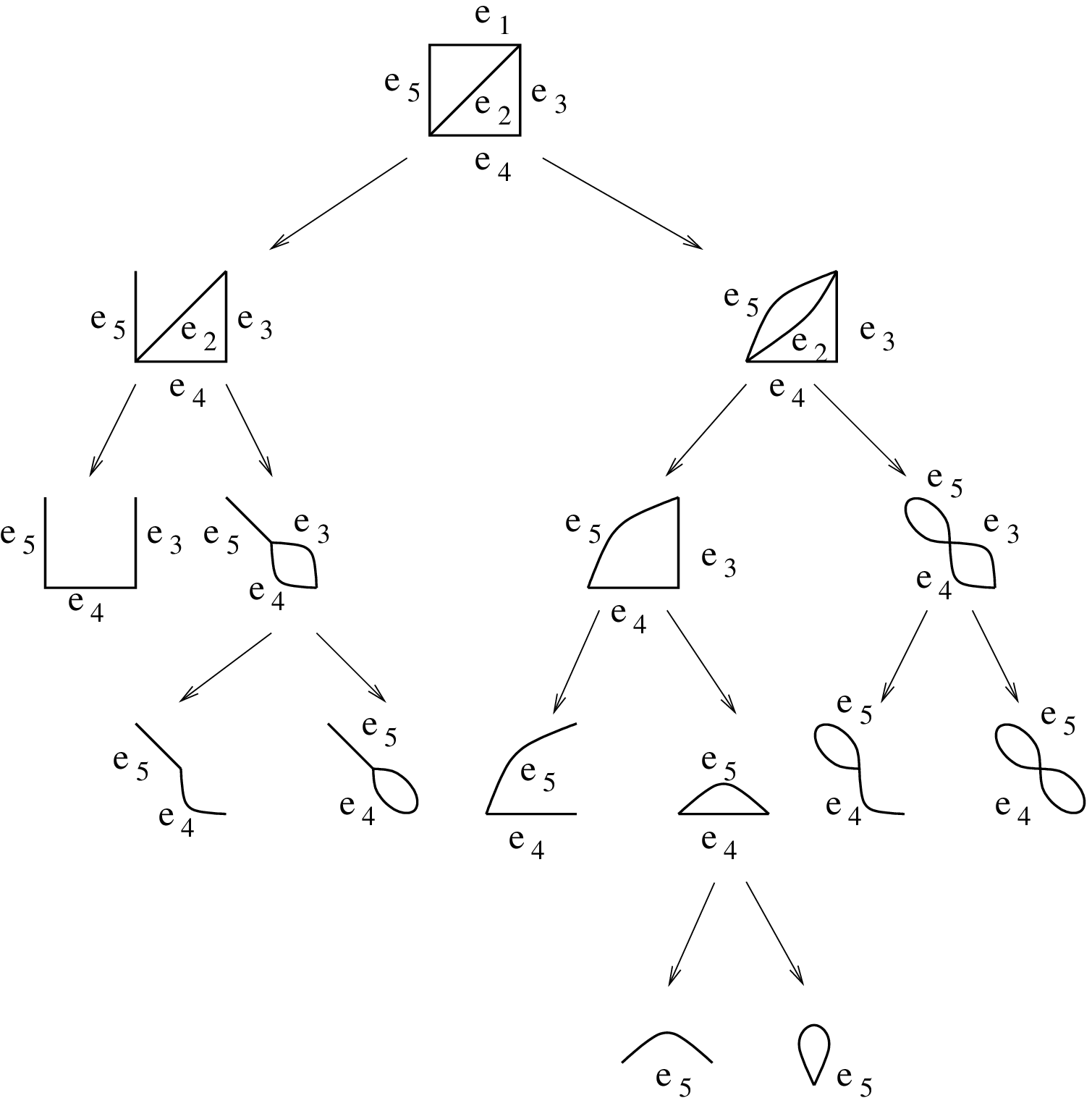,height=8.9cm}}}
\centerline{ Fig.1.4:\ Computational tree for the Tutte polynomial of
{\psfig{figure=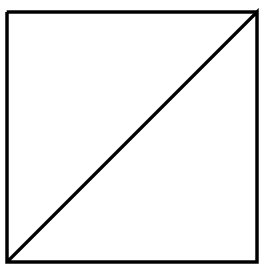,height=0.4cm}}.}

Leaves of the binary computational tree are in bijection with
spanning trees of $G$. For a leaf $F$, the associated spanning tree
is composed of isthmuses of $F$ (called internally active edges)
and edges of $G$ contracted on the way from $G$ to $F$ (called
internally inactive edges). This approach gives probably the
simplest description (and proof) of the celebrated Tutte formula
for the Tutte polynomial of a graph with ordered vertices:
$$\chi (G;x,y) = \sum_T x^{IA(T)}y^{EA(T)}$$
where the sum is taken over all spanning trees of $G$, and for
a spanning tree $T$ and associated leaf of the binary computational tree, $F$,
$IA(T)$ denote the number of internally active edges of $G$ that is
isthmuses of $F$, and $EA(T)$ denote the number of externally active 
edges of $G$ that is loops of $F$. 

With this setting the proof of 
the Tutte formula is an easy task. Also Exercise 1.11 follows easily, for
example no leaf can have more than $d(G)$ isthmuses and there
is exactly one leaf with $d(G)$ isthmuses and $p$ loops (we choose a path,
in the binary computational tree, composed only of $(p_1(G)-p)$ deletions.
This proves part (1) of 1.11.
\ \\

A graph is called planar if it can be embedded in a plane and
it is called plane if it is embedded in a plane.

For a plane graph $G$ we define its dual graph 
$G^\star$ in the following way:
If $G$ is connected then the vertices of $G^\star$ are connected 
components of $R^2-G$. To every edge $e$ of $G$
corresponds the dual edge $e^*$ of $G^\star$ joining vertices
(regions of $R^2-G$) separated by $e$, see Fig.~1.5.
In particular, $G$ and $G^\star$ have the same number of edges.
$G^\star$ can have different embeddings in a plane 
but $G^\star$ for a plane graph $G$ is uniquely 
defined in $S^2=R^2 \cup \infty$.
\ \\
\ \\
%\vspace*{1.5in}\centerline{\Psfig{figure=Rys.1.3}}
\centerline{{\psfig{figure=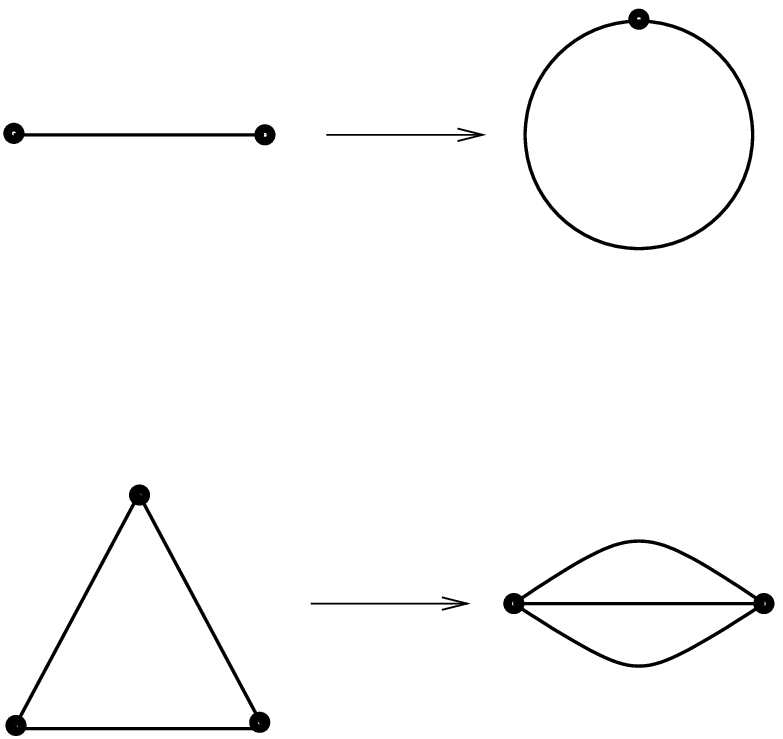,height=2.9cm}}} 
\begin{center}
Fig.~1.5
\end{center}

If the graph $G$ is not connected then
$G^\star$ is, by definition, a disjoint sum of graphs dual to
components of $G$.

For different embeddings of a connected planar graph $G$ we can get different 
duals
(even if $G$ is 2-connected), see Fig. 1.6. However if $G$  is a 
3-connected (Def. 1.7) 
planar graph then $G^\star $ is uniquely 
defined\footnote{It is known that a 3-connected planar graph 
has unique embedding in $R^2\cup \infty = S^2$, \cite{Tut-2}.}.

\ \\
\centerline{{\psfig{figure=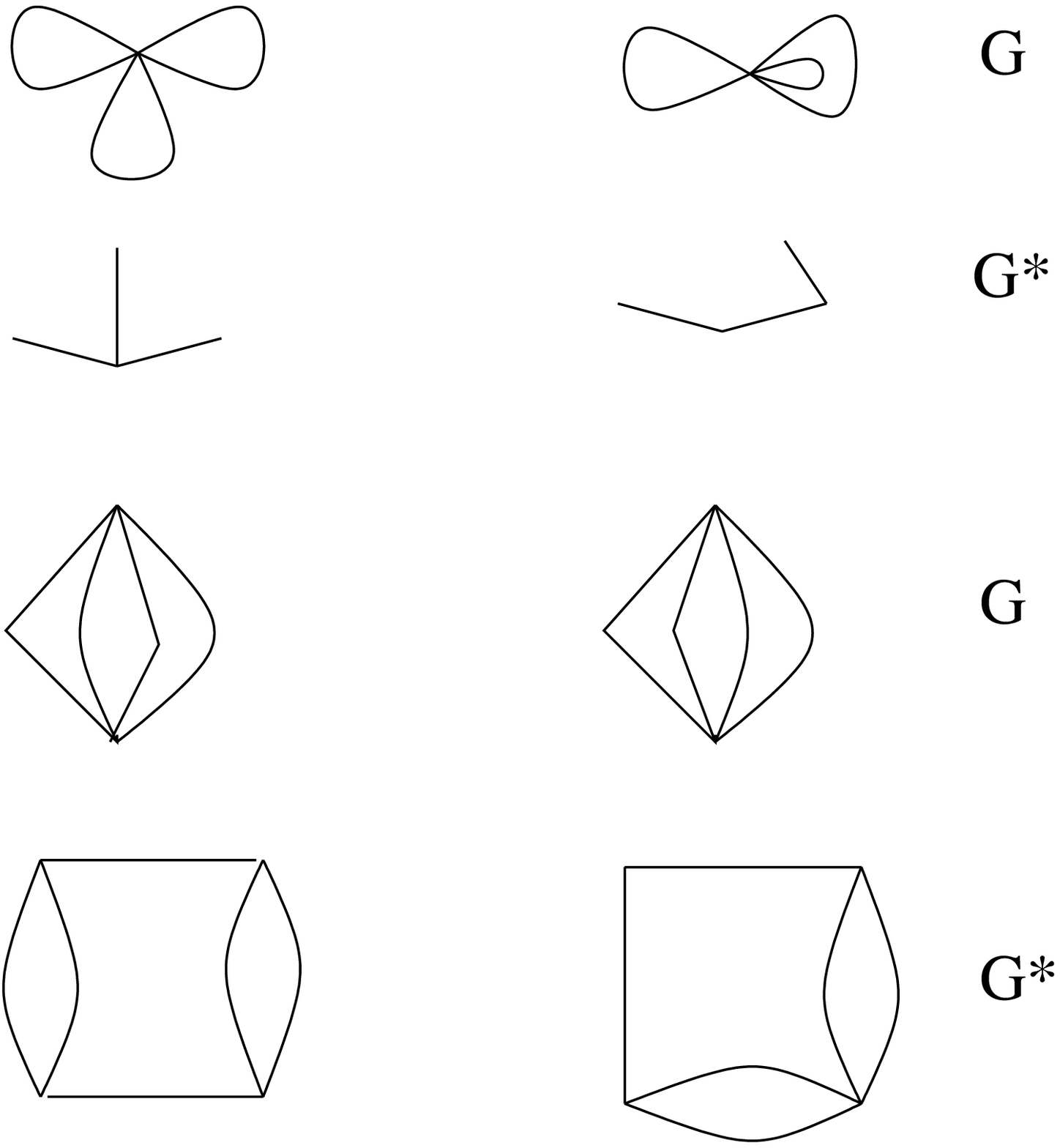,height=5.9cm}}}
\begin{center}
Fig.~1.6
\end{center}

\begin{exercise}\label{V.1.13} 
Prove that for the plane graph:
\begin{enumerate}
\item [(i)]
$$(G^\star)^\star = G.$$
\item [(ii)]
If $G$ is 2-connected then $G^\star$ is 2-connected.\\
Hint.\ Show that if $G=G_1*G_2$ then $G^\star = G_1^\star *G_2^\star$
\end{enumerate}
\end{exercise}

\begin{theorem}\label{V.1.14}
If $G$ is a planar graph then
\begin{enumerate}
\item[(1)] $\chi (G;x,y) = \chi (G^\star;y,x)$ 
\\
a similar identity holds for the Kauffman bracket polynomial
\item[(2)]
 $\langle G\rangle_{\mu,A,B} = \langle G^\star\rangle_{\mu,B,A}$
\end{enumerate}
\end{theorem}

Proof.
\begin{enumerate}
\item[(1)] First, let us note that a single edge is dual to a single 
loop (Fig.~1.5), so the result is true for a graph with one edge.
Next, we make an easy induction with respect to the number 
of edges in the graph.

\item[(2)]
 Similarly as in (1) we can apply induction or use directly 
Theorem 1.6.
%Theorem \ref{4:1.6}.
\end{enumerate}

\begin{exercise}\label{V.1.15}
Let us recall that for a given graph $G$ and a positive integer $\lambda$
we define $C(G,\lambda)$ to be equal to the number of possible ways of
coloring the vertices of $G$ in $\lambda$ colors
in such a way that the edges have endpoints colored
in different colors. Show that:
\begin{enumerate}

\item[(1)]
\begin{enumerate}
\item[(a)]
 If $G$ has $n$ vertices and no edges then
$C(G,\lambda) = \lambda^n$. If $G$ contains a loop then
$C(G,\lambda) = 0$.
\item [(b)] If an edge $e$ is not a loop then
$$C(G,\lambda) = C(G-e,\lambda) - C(G/e,\lambda)$$
\end{enumerate}
\item[(2)]
 Prove that $C(G,\lambda)$ and $\chi (G;x,y)$ 
are related by the formula 
$$C(G,\lambda)= (-1)^{|V(G)|-p_0(G)}\cdot\lambda^{p_0(G)}
\cdot\chi(G;1-\lambda,0)$$ 
(the chromatic polynomial is determined by the Tutte polynomial 
and the number of components of $G$. Recall that $|V(G)|$ is determined by 
the Tutte polynomial, Exercise 1.11(b)).
\item[(3)] 
\begin{enumerate}
\item[(a)] Show that if $G$ has at least one edge than $C(G,1)=0$.
\item[(b)] Show that the number\footnote{For a connected graph $G$, 
$v_{1,0}$ is named the {\it chromatic invariant} of $G$ \cite{Big}.} 
$v_{1,0}$ introduced in Corollary 1.9 is equal, 
up to the sign, to the derivative of the chromatic polynomial evaluated 
at $1$; we have $v_{1,0}= (-1)^{V(G)-(p_0(G)-1)}C'(G,\lambda)_{\lambda=1}$.
\end{enumerate}
\item[(4)] Show that $C(G,\lambda)$ and $<G>_{\mu,A,B}$
are related by the formula
$$C(G,\lambda)= 
(-1)^{|V(G)|-1}(\frac{B}{A})^{p_0(G)+1-p_1(G)}A^{-E(G)}<G>_{\mu,A,B}$$
%$$\mu^{p_0 (G)-1} B^{p_1 (G)} A^{E(G)-p_1(G)}C(G,\lambda)= (-1)^{|V(G)|-p_0(G)}\cdot\lambda^{p_0(G)}<G>_{\mu,A,B}$$ 
for $\mu = -\frac{B}{A}$,  
%$B+A\mu =0$, 
%$1-\lambda =x = 1+\frac{B}{A}\mu =  1- (\frac{B}{A})^2$,\\
$\lambda = (\frac{B}{A})^2= \mu^2$.
%\item[(5)] Relation between chromatic polynomial and Jones polynomial via 
%Kauffman bracket. They coincide if the longest cycle in the graph is bigger...
\end{enumerate}
\end{exercise}
The property (2) implies that $C(G,\lambda)$ is an invariant
polynomial (in variable $\lambda$) of the graph $G$ and it can be defined
by properties 1(a) and 1(b). The polynomial $C(G,\lambda)$ 
is called the chromatic polynomial of the graph $G$.

\begin{corollary}\label{V.1.16}
\begin{enumerate}
\item[(i)] $\chi(G;1,1) \geq |E(G)|$ for a connected graph without loops and 
isthmuses and the equality holds only  
for  a polygonal graph or the generalized theta curve, or the 
graph 
\parbox{2.2cm}{\psfig{figure=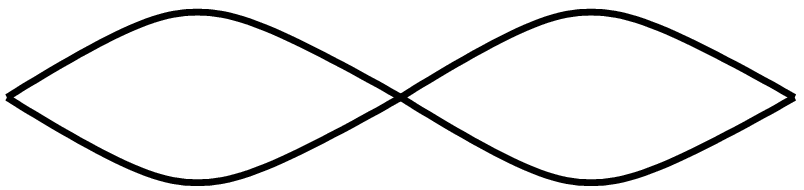,height=0.5cm}}. 
$\chi(G;1,1)$ is named the 
determinant or complexity of a connected graph and equal to the 
number of spanning trees of the graph, as explained before (compare 
also historical remarks in Subsection V.1.1). 
\item[(ii)] (Murasugi) The determinant of a non-split alternating link 
is no less than its crossing number (minimal number of crossings). 
Furthermore the equality holds only for a $(2,k)$ torus links and 
the connected sum of two Hopf links.
\item[(iii)] $\chi(G,1,1) \geq |E(G)| +10$  for a 3-connected graph with  
at least 4 vertices and the equality holds
for $K_4$. 
\item[(iv)] $\chi(G,1,1) \geq 2|E(G)| + 10$ if $G$ is 3-connected graph with 
at least 4 vertices and different from $K_4$.
%and the equality holds for...
% Or better $\chi(G,1,1) \geq 8|E(G)|-32$
\end{enumerate}
\end{corollary}

Part (ii) follows from (i) when diagrams are converted to links 
and a (monochromatic) graph translates into an alternating link.
This is explained below.

There are several ways of translating diagrams of links into 
graphs.
We begin with a classic one (introduced by Tait in 1876) 
which seems to be the most useful,
up to now.  We will consider 2-color graphs. That is: their 
edges will be colored in black and white and denoted by $b$ and $w$,
respectively. Frequently, in literature, black edges
are denoted positive ($+$) and white edges are called negative ($-$).

Now, given a connected diagram of a link $L$, we can color
connected components of the complement of the diagram in the plane
in black and white. We color them so that neighboring components
are colored in different colors --- exactly as in the construction
of Goeritz matrix (checkerboard coloring). 
Subsequently, we construct a planar graph
$G(L)$. Vertices of $G(L)$ represent black components of the divided
plane and edges represent crossings. Moreover, the edge associated
to a given crossing is either black or white depending on the situation
described in Fig.~1.7. \\
\ \\
%\vspace*{1.5in}\centerline{\Psfig{figure=Rys.1.4}}
\centerline{{\psfig{figure=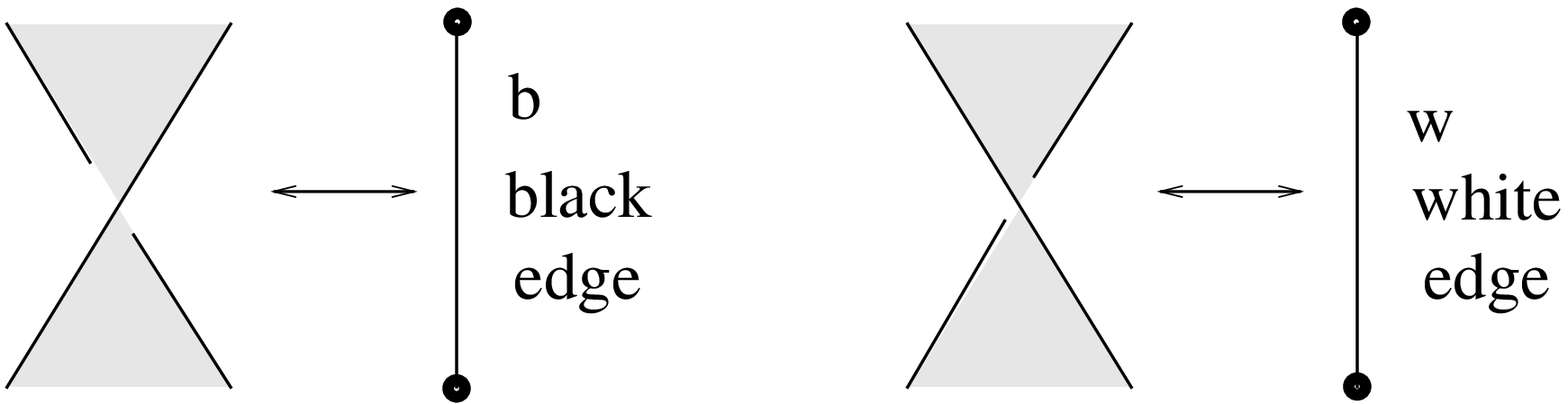,height=2.3cm}}} 
\begin{center}
Fig.~1.7
\end{center}

Examples of graphs associated to diagrams are pictured in Fig.~1.8.
\ \\
\ \\
%\vspace*{2.5in}\centerline{\Psfig{figure=Rys.1.5}}
\centerline{{\psfig{figure=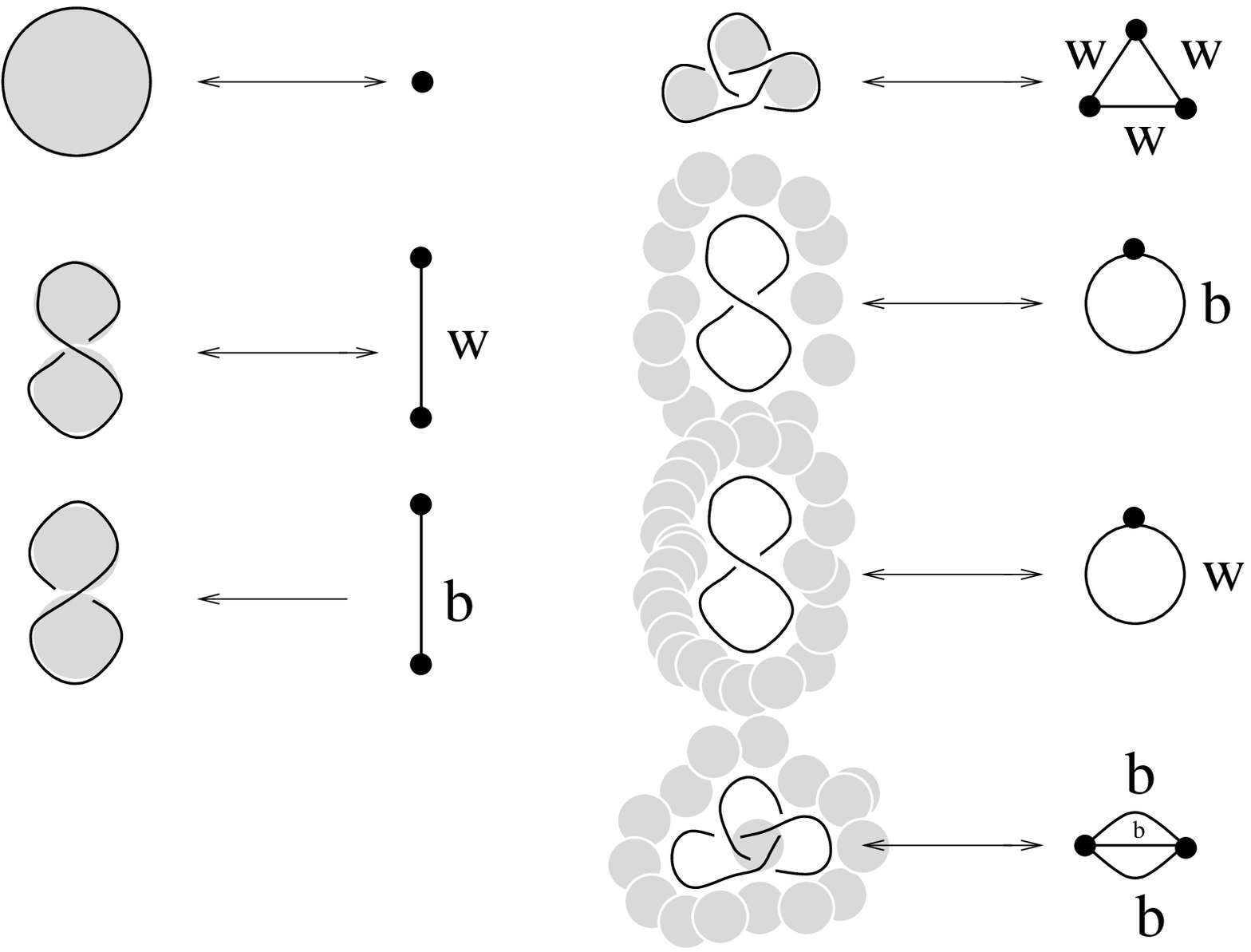,height=7.1cm}}}
\begin{center}
Fig.~1.8
\end{center}
\ \\
\centerline{{\psfig{figure=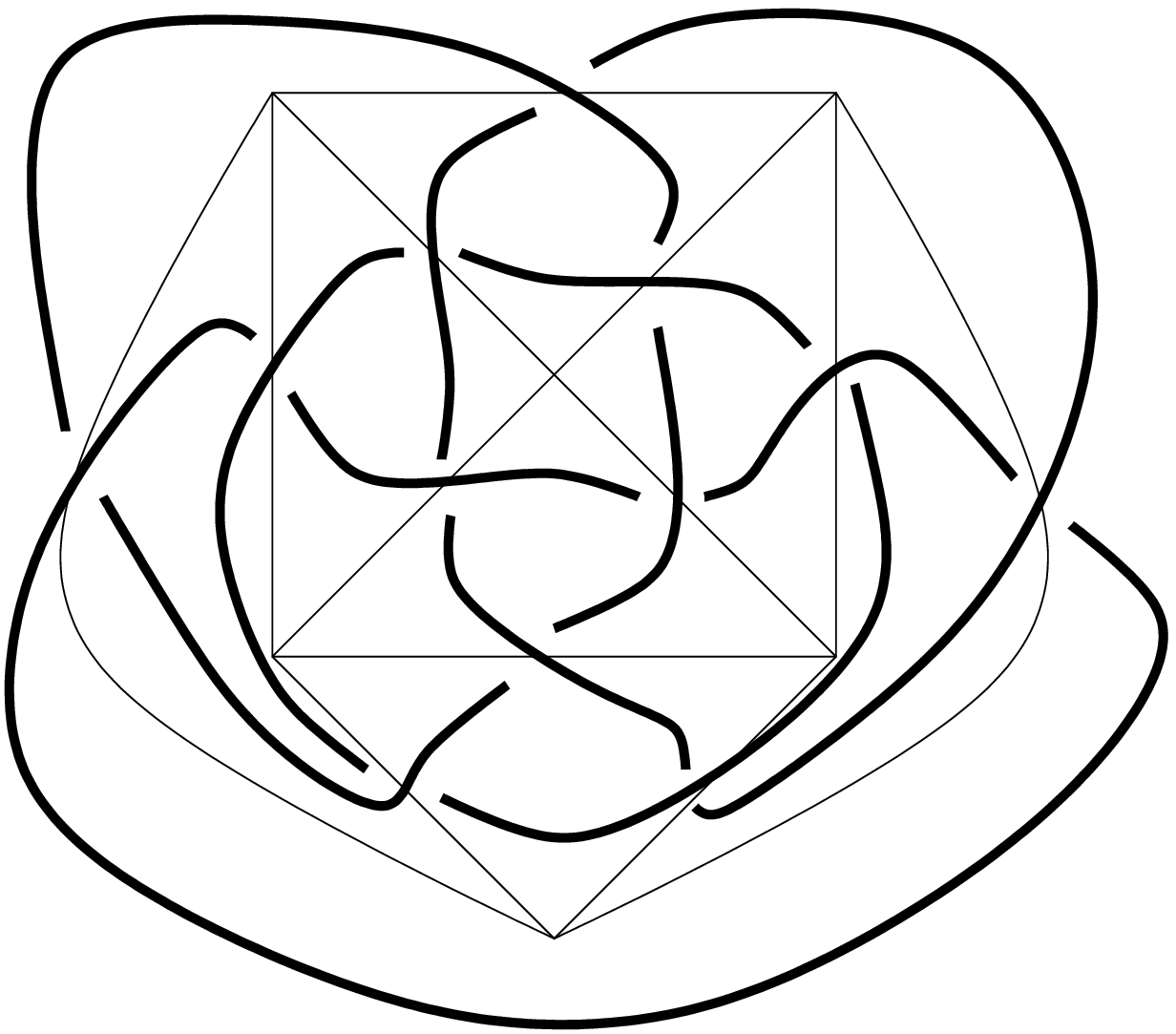,height=3.7cm}}}
\begin{center}
Fig. 1.9;\ Octahedral graph (with all $b$ edges) and 
the associated link diagram
\end{center}

We see that the graph associated to a diagram of a link does not depend
only on the diagram but also on the 
checkerboard colorings of the plane containing the diagram.

\begin{lemma}\label{V.1.16}
Let $L$ be a connected diagram of a link.
Then the related two checkerboard colorings of the plane 
yield two dual graphs and the duality interchanges
the colors of edges (black to white and vice versa).
\end{lemma}

The proof follows immediately from the construction
of $G(L)$.

If the diagram $L$ is oriented then the edges of $G(L)$ are not only 
colored in black or white but also signed ($+$ or $-$).
By definition, the sign of an edge is equal to the sign of the
crossing to which the edge is associated, (c.f.~Fig.1.10).\\
\ \\
%\vspace*{2in}\centerline{\Psfig{figure=Rys.1.6}}
\centerline{{\psfig{figure=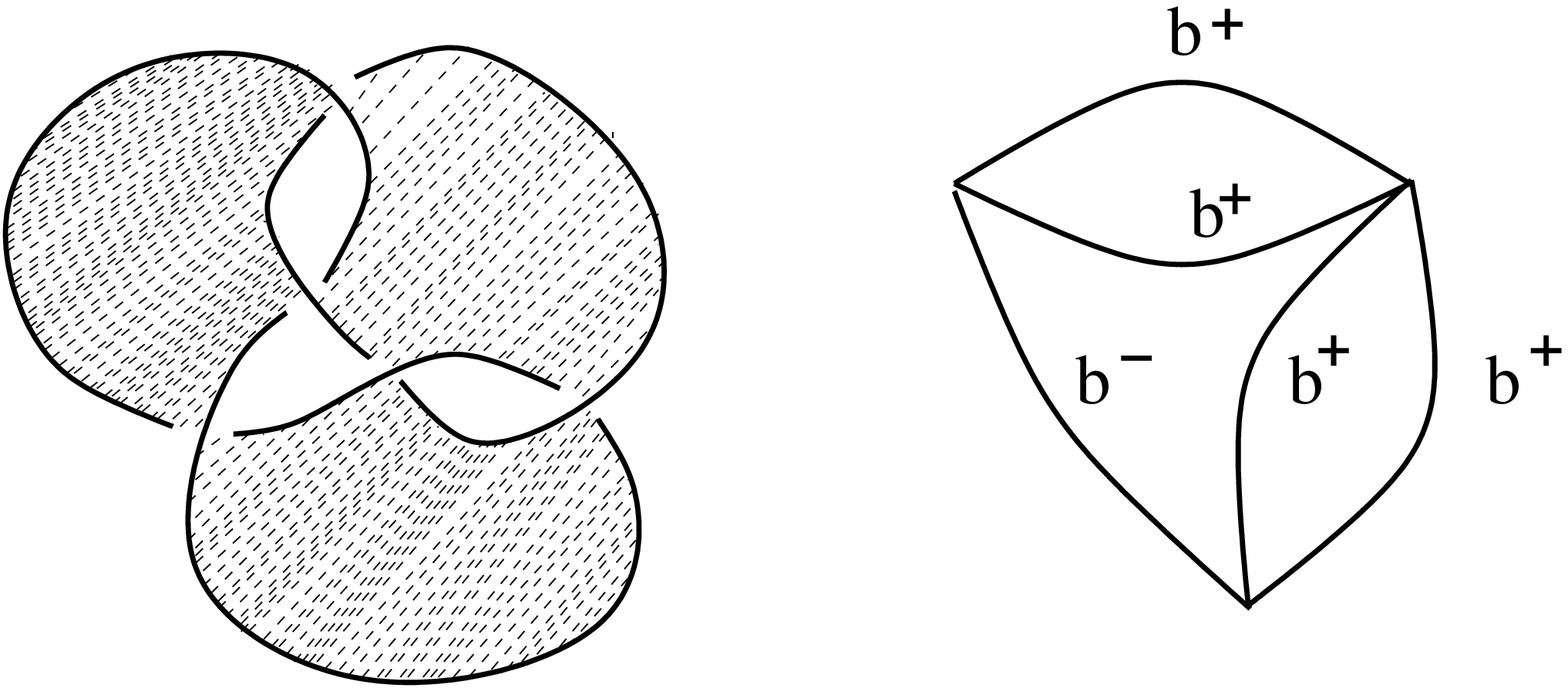,height=3.5cm}}}
\begin{center}
Fig.~1.10
\end{center}

Let us note that not all signed
2-color graphs are associated to
diagrams. The simplest example is the graph
$\stackrel{b^+}{\bullet\!\!\!\!\longrightarrow\!\!\!\!\bullet}$;
compare Section 5 and \cite{Ko}.

Another way of translating of signed 2-color planar graphs
to oriented links comes from an idea of Jaeger \cite{Ja-1} 
which was developed in \cite{P-P-1}, see also \cite{A-P-R}.

For a given 2-color ($b$ or $w$), signed ($+$ or $-$), planar graph
$G$ we associate an oriented diagram of a link $D(G)$ together
with a checkerboard coloring of the plane. We do it according
to the rules explained in Fig.~1.11:
\ \\
%\vspace*{2in}\centerline{\Psfig{figure=Rys.1.7}}
\ \\
\centerline{{\psfig{figure=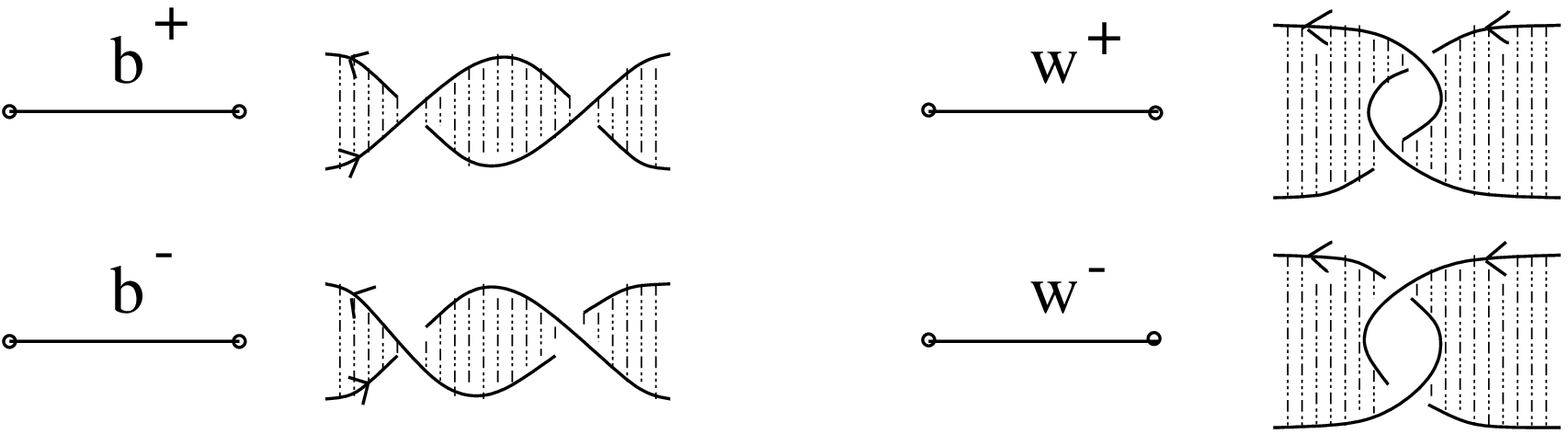,height=2.1cm}}}
\begin{center}
Fig.~1.11
\end{center}

Examples illustrating the construction of $D(G)$ 
are shown in Fig.~1.12.
\ \\
%\vspace*{3.5in}\centerline{\Psfig{figure=Rys.1.8}}
\ \\
\centerline{{\psfig{figure=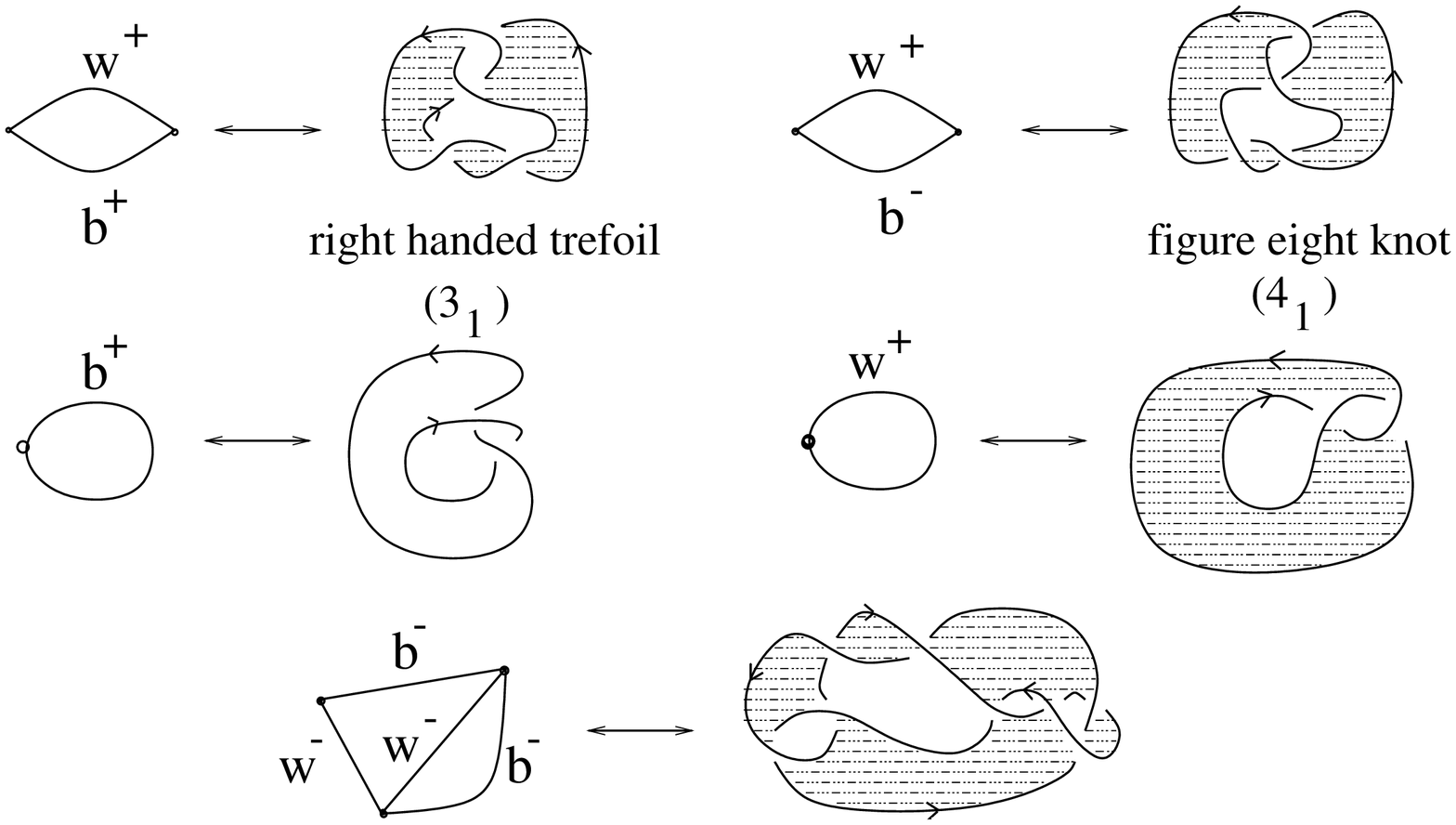,height=6.1cm}}}
\begin{center}
Fig.~1.12
\end{center}

Diagrams of knots which are of the form $D(G)$ for some
2-color signed graph are called matched diagrams. 
Every 2-bridge link has a matched diagram \cite{P-10} but 
probably it is not true that any link has a matched diagram 
however the existence of a counterexample is still an open problem
\footnote{A similar concept was considered before by J.H.Conway
who constructed knots which probably do not possess a matched diagram
\cite{Kir}.}.

\begin{conjecture}
\begin{enumerate}
\item[(i)] Not every link has a matched diagram.
\item[(ii)] Every oriented link is $t_3$-move 
(\parbox{5.2cm}{\psfig{figure=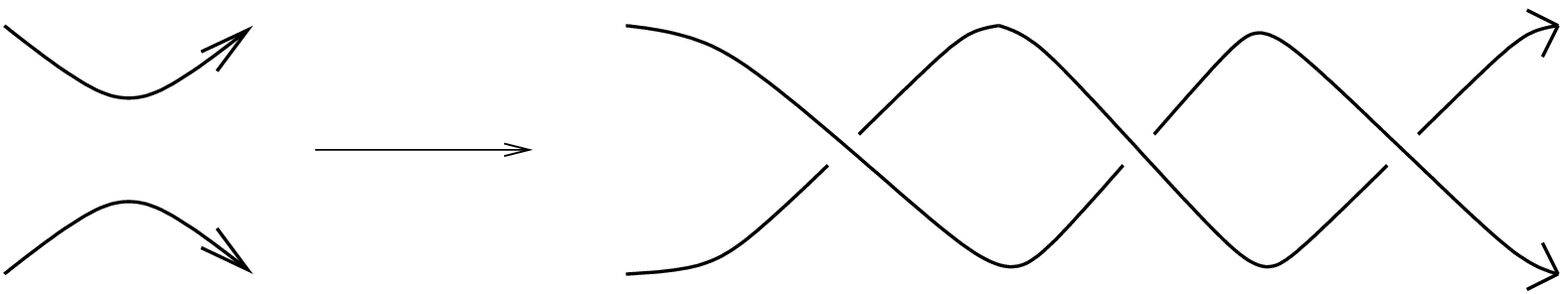,height=0.9cm}}) 
 equivalent to a link 
with a matched diagram.
\end{enumerate}
\end{conjecture}

We can further extend translation from plane graphs to links by 
considering weighted (by elements of $Q \cup \infty$) graphs and 
decorating the vertices of a medial graph\footnote{A medial graph 
$G^m$ of a plane graph $G$ is constructed by choosing vertices (of degree 
4) in the middle of edges of $G$ and connecting them along edges 
of $G$ as in the link diagram. In other words, $G^m$ is 
obtained from a diagram associated to the graph by identifying the
overcrossing with the underscrossing at every crossing of the diagram.}
by $\frac{p}{q}$-rational tangles  (compare Chapter VI). For integer tangles 
this translation was crucial in the proof that computing most of 
substitutions in Jones, Homflypt and Kauffman polynomials is NP-hard 
\cite{J-V-W}, compare Corollary 5.17.

\subsection{Polynomial invariants of chromatic graphs}
In this part we give, after \cite{P-P-1,P-P-2}, a historical introduction to
chromatic polynomials. We allow general weights on edges of a graph 
and develop formulas allowing a universal change of variables in 
dealing with various versions of polynomials.

$\chi(G;1,1)=\tau(G)$ denotes, as before, the complexity of 
the graph $G$, that is, the number of its spanning trees. 
Invariant $\tau$ was introduced and studied by 
Kirchhoff \cite{Kirch}. It has been noted in \cite{BSST} that if $e$ is 
an edge of $G$ that is not a loop  then 
$$ \tau(G) = \tau(G-e) + \tau(G/e)$$
As noted by Tutte (\cite{Tut-4}; p. 51), this equality had been long 
familiar to the authors of \cite{BSST}. The equality inspired Tutte 
to investigate all graphs invariants, $W(G)$, which satisfy 
the identity 
$$W(G) = W(G-e) + W(G/e)$$

This led to the discovery of the dichromatic polynomial and its variant, 
the Tutte polynomial \cite{Tut-1,Tut-3}.\footnote{H. Whitney 
\cite{Whit-1,Whit-2} was considering graph invariants $m_{i,j}$
which are essentially the coefficients of the dichromatic polynomial. 
He also analyzed closer the topological graph invariants $m_i$
which corresponds to the coefficient of the flow polynomial \cite{Whit-3}.
 R.M.Foster noticed \cite{Whit-1} that $ m_{i,j} $ invariants
satisfy $m_{i,j}(G) = m_{i,j}(G-e) + m_{i-1,j}(G/e)$.
 }  The ring of graphs from \cite{Tut-1}, obtained by taking the module 
of formal linear combinations of graphs and dividing this module 
by a submodule generated by 
deleting-contracting linear relations, can be thought as a precursor 
of skein modules of links discussed in Chapter IX.

C.M.Fortuin and P.W.Kastelyn generalized the dichromatic polynomial
to chromatic (weighted)
graphs \cite{F-K} (compare also O.J.Heilmann \cite{Hei}). The research of
\cite{F-K,Hei} was  motivated
by ``statistical mechanics" considerations. One should stress here 
that only slightly  earlier H.N.V.Temperley discovered that the
partition function for the Potts model is equivalent to 
the dichromatic polynomial
of the underlying graph \cite{Ess,T-L}.

The dichromatic polynomial for chromatic graphs gained new importance 
after the Jones discovery of new polynomial invariants of links 
and the observation of Thistlethwaite that the Jones polynomial 
of links is closely related to the Tutte polynomial of graphs.
Several researches rediscovered the dichromatic polynomial 
and analyzed its properties \cite{K-9,M-7,Tral-2,P-P-1,Yet,Zas}.

The following version of the dichromatic polynomial is motivated
by connections between graphs
and links.

A {\it chromatic graph} is a graph with a function $c$ on the edges, where
$c: (E(G) \to Z \times \{d,l\}$. The first element of the pair $c(e)$
is called the {\it color} and the second the {\it attribute} ($d$ - for
dark, $l$ for light) of the edge $e$. Note that chromatic graphs are
extensions of signed graphs were the attribute of an edge corresponds
to its sign (plus or minus) or $b$,$w$ colored graphs 
considered in the previous section. 
The {\it dual} to a connected chromatic plane graph $G$
is the graph $G^*=(V(G^*),E(G^*))$ where $V(G^*)$ and $E(G^*)$ are
defined as for non-chromatic graphs and the edge 
$e^*$ dual to $e$ has assigned the same color as $e$ and
the opposite attribute. Furthermore, $\bar G$ denotes the graph 
obtained from $G$ by reversing attributes of every edge (following 
knot theory analogy we say that 
$\overline G$ is a {\it mirror image} of $G$).

\begin{theorem}\label{R(G)}
%\marginpar{\tiny label: R(G)}
There exists an invariant of chromatic graphs $R(G) =$\\
$ R(G; \mu , r_1,r_2,A_i,B_i)$
which is uniquely defined by the following properties:
\begin{enumerate}
\item[(1)]
$R(T_n) = \mu ^{n-1}$; where $T_n$ is the $n$-vertex graph with no edges,
\item[(2)]
$$R(G) = (\frac{r_1}{\mu})^{\epsilon (d_i)}B_i R(G-d_i) +
r_2^{\delta(d_i)} A_i R(G/d_i)$$
$$R(G) = (\frac{r_1}{\mu})^{\epsilon (l_i)}A_i R(G-l_i) +
r_2^{\delta(l_i)} B_i R(G/l_i)$$
where
\[ \epsilon(e)  = \left \{ \begin{array}{ll}
                             0 & \mbox{if $e$ is not an isthmus}\\
                             1 & \mbox{if $e$ is  an isthmus}
\end{array}
\right. \]

\[ \delta(e)  = \left \{ \begin{array}{ll}
                             0 & \mbox{if $e$ is not a loop}\\
                             1 & \mbox{if $e$ is  a loop}
\end{array}
\right. \]
\end{enumerate}
\end{theorem}

Our variables have been chosen in such a way that the invariants 
for a plane graph $G$
and its dual $G^*$ are symmetric in the following sense:

\begin{lemma} \label{3.2}
%\marginpar{\tiny label : 3.2}
If $G$ is a plane graph then
$$R(G) = R(G; \mu , r_1,r_2,A_i,B_i) = R(G^*; \mu , r_2,r_1,A_i,B_i)$$
\end{lemma}

Note, that  $R(G)$ is a  2-isomorphism invariant of connected
chromatic graphs. Generally, when $G$ is not necessarily connected
and $\mu \neq 1$ then the polynomial measures also the number
of connected components of the graph. If we put $\mu = 1$ then 
the dichromatic polynomial, $R$, and its property
described in  Lemma
\ref {3.2} can be extended to matroids
(see \cite{Zas} for a full analysis of the Tutte polynomial of
colored matroids) or more generally to colored Tutte set
systems (see \cite{P-P-0}).

Let $S$ denote  a subset of edges of a graph $G$. 
By $(G:S)$ we denote the subgraph
of $G$ which includes all the vertices of $G$ but only edges in $S$.
The polynomial $R(G)$ has the following ``state model" expansion:

\begin{lemma} \label {3.3}
%\marginpar{\tiny label: 3.3 }
\begin{eqnarray*}
 & & R(G; \mu,r_1,r_2,A_i,B_i) = \\
  & & \mu ^{p_0(G)-1}\sum_{S \in 2^{E(G)}} 
r_1 ^{p_0(G:S)-p_0(G)}r_2^{p_1(G:S)}( \prod _{i=1}^n
A_i^{\alpha_i+\alpha_i'}
\cdot B_i^{\beta_i+\beta_i'})
\end{eqnarray*}
where the sum is taken over all subsets of $E(G)$, and $\alpha_i$
is the number of dark edges in $S$ of the $i^{th}$ color, $\alpha_i'$ is
the number of light edges in $E(G)-S$ of the
$ i^{th}$ color, $\beta_i$ is the number of dark edges in $E(G)-S$  of
the $i^{th}$ color, and $\beta_i'$ is the
number of light edges in $S$ of the $i^{th}$ color.
\end{lemma}

In the above lemma we consider a subset $S$ of the set of edges to be
the state of $G$ in the sense  that edges in $S$ are marked to be contracted
and the edges in $E(G) - S$ are marked to be deleted.

Below we list a few easy but useful properties of $R(G)$.

\begin{lemma}\label{prop of R(G)}
%\marginpar{\tiny label: prop of R(G)}
\
\begin{description}
\item
[(i)] $R(\overline  G; \mu ,r_1,r_2,A_i,B_i) = R(G;\mu ,r_1,r_2, B_i,A_i)$
\item
[(ii)] For any $i$, the number of $i^{th}$ colored edges 
of $G$ is equal to $ \alpha _i+ \alpha '_i+\beta _i + \beta '_i$ which 
is equal to  the highest power of $A_i$ in $R(G)$
\item[(iii)] If $G_1*G_2 $ is a one vertex product
 of $G_1$ and $G_2$ and $G_1 \sqcup G_2$ is
a disjoint sum of $G_1$ and $G_2$ then
$$R(G_1 \sqcup G_2) = \mu R(G_1*G_2) = \mu R(G_1)  R(G_2)$$
\item[(iv)] If $G$ is a loop or isthmus then we have
$$ R(\parbox{0.5cm}{\psfig{figure=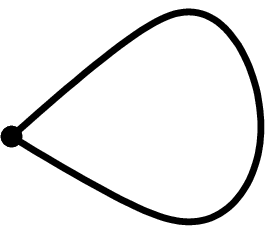,height=0.4cm}}
 _{d_i}) = B_i+r_2A_i$$
$$ R(\parbox{0.5cm}{\psfig{figure=loop.eps,height=0.4cm}}_{l_i}) = 
A_i+r_2B_i$$
$$ R(\parbox{0.5cm}{\psfig{figure=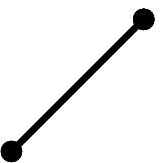,height=0.4cm}}_{d_i}) = 
A_i+r_1B_i$$
$$ R(\parbox{0.5cm}{\psfig{figure=interval.eps,height=0.4cm}}_{l_i}) = 
B_i+r_1A_i$$
where $d_i$ (resp.,$l_i$) denotes a dark (resp., a light) edge 
of the $i^{th}$ color.
\item
[(v)] If $Q(G;t,z)$ is the Traldi's version\footnote{Traldi's polynomial 
is characterized by the following properties:\ (i) $Q(T_n;t,z)=t^n$,\ 
(ii) $Q(G_1 \sqcup G_2;t,z)=Q(G_1;t,z)Q(G_2;t,z)$,\ (iii) If $e$ is 
not a loop then $Q(G;t,z)= Q(G-e;t,z) + w(e)Q(G/e;t,z)$,\ (iv) if 
$e$ is a loop then $Q(G;t,z)= (1+w(e)z)Q(G-e;t,z)$,\ (v) if $e$ is 
an isthmus then $Q(G;t,z)= (w(e)+t)Q(G/e;t,z)$,\ (vi) if the weight 
$w(e)$ of an edge $e$ is equal to zero then $Q(G;t,z)= Q(G-e;t,z)$.} 
of the dichromatic polynomial
\cite{Tral-2} then
$$Q(G;t,z) = \frac{tR(G;\mu ,r_1,r_2,A_i,B_i)}{\prod_{i} B_i^{E_i(G)}}$$
where $r_1 = \mu = t, r_2=z$,
$E_i(G)$ denotes
the number of $i^{th}$ colored edges in $G$,
and the weight, $w(e)$, of an edge $e$
of $G$ is defined by:
\[ w(e) = \left \{ \begin{array}{ll}
                    \frac{A_i}{B_i} & \mbox{if $e$ is a $d_i$ edge}\\
                    \frac{B_i}{A_i} & \mbox{if $e$ is an $l_i$ edge}
\end{array}
              \right.   \]
\end{description}
\end{lemma}

Note that both versions of the dichromatic polynomial are equivalent because,
by Lemma \ref{3.3} (ii),
$E_i(G)$ is determined by  $R(G)$. Furthermore $Q(G;t,z)$ determines $E_i(G)$
and  $p_0(G)$.

%  TO BE continued or deleted or contracted :-)

% A {\it chromatic graph} is a graph with a function $c$ on the edges, where 
%$c: (E(G) \to Z \times \{d,l\}$. The first element of the pair $c(e)$ 
%is called the {\it color} and the second the {\it attribute} ($d$ - for 
%dark, $l$ for light) of the edge $e$. Note that chromatic graphs are 
%extensions of signed graphs were the attribute of an edge corresponds 
%to its sign (plus or minus).
%%%%%%%%%%%%%%%%%%%%%%%%%%%%%%%%%%%%%%%%%%%%%%%%%%%%%%%%%%%%%%%%%%%%%%%%

\newpage
\section{Setoids and Dichromatic Hopf algebras}\label{V.2}
In this section we sketch two generalizations of the Tutte polynomial,
 $\chi(S)$, or, more precisely, deletion-contraction method.
The first generalization considers, instead of graphs, general object called
setoids or group system, in the second approach we work with graphs but
initial data are also graphs (finite type invariants of links have given
motivation here).

\begin{definition}\label{1}
A setoid $S=(E,T)$ is a pair composed of a set $E$ and a set of its
subsets $T\in 2^E$. By analogy with graph we call elements of $E$ -
edges and elements of $T$ -- trees (corresponding 
to spanning trees in a graph). An isthmus is an edge contained in all
trees of the setoid. A loop is an edge belonging to no tree.
The setoid $S-e$ is defined to be $(E-{e},T\cap 2^{E-{e}})$
that is trees of $S-e$ are elements of
$T$ which do not contain $e$. We say that $S-e$ is obtained from $S$
by a deleting operation. The setoid $S/e$ is defined to be $(E-{e},T'')$
where $t\in T''$ if $t\cup e$ is in $T$. We say that $S/e$ is obtained
from $S$ by a contracting operation.
For a setoid $S$ we associate complementary, or dual, setoid $S^*=(E,2^E-T)$.
If $E$ is finite, we say that
$S$ is a finite setoid. If $T$ is finite and every element of $T$
is finite we say that $S$ is finitely presented.

\end{definition}
To define the Tutte polynomial of a setoid, $\chi(S)$, we can
follow definition for graphs, except that one have to add
a special conditions to guarantee independence on the
orderings of edges. It is convenient to define special setoids
as setoids for which the result of computation of a polynomial does
not depend on the order of computation.
\begin{definition}
A special setoid is a setoid for which we can associate a polynomial
invariant of links, Tutte polynomial, $\chi(S) \in Z[x,y]$ satisfying:
\begin{enumerate}
\item[(i)] If $S$ has only one tree, then we put
$\chi(S)= x^iy^j$ where $i$ is the number of
elements in the tree and $j$ the number of elements
not in the tree. Furthermore for the empty $T$, we put $\chi(S)= 1$.
\item[(ii)] If $e$ is an edge of a setoid which is neither
an isthmus nor a loop then we have a deleting-contracting formula:
$$\chi(S)= \chi(S-e) + \chi(S/e) $$
\end{enumerate}
\end{definition}
Graphs are examples of special setoids.
More generally matroids are special setoids\footnote{In the case of matroids,
trees are called basis and isthmuses -- co-loops.}.
Below we describe a class of setoids slightly generalizing matroids.
\begin{definition}
We say that a setoid satisfies an exchange property (shortly E-setoid) if
for every tree $t$ and an edge $e\in t$ not an isthmus of $S$ there is an edge $
f$ outside $t$
such that $t- \{e\} \cup \{f\}$ is a tree. 
Furthermore for every edge $f$ not in $t$ not a loop of $S$
there is an edge $e$ in $t$ such that $t- \{e\} \cup \{f\}$ is a tree.
\end{definition}
Among properties of $E$-setoids we list a few of interest to us.
\begin{enumerate}
\item[(i)] A dual to an $E$-setoid is an $E$-setoid.
\item[(ii)] Exchange property is a hereditary property,
that is if $e$ is neither an isthmus nor a loop
of an $E$-setoid $S$ then $S-e$ and $S/e$ are $E$-setoids.
\item[(iii)] If $e$ and $f$ are neither  isthmuses nor a loops of
an $E$-setoid $S$ then
\begin{enumerate}
\item[(a)]
$f$ cannot be a loop of $S-e$ or isthmus of $S/e$.
\item[(b)]
If $f$ is an isthmus of of $S-e$ then $e$ is an isthmus of of $S-f$.
\item[(c)]
If $f$ is a loop of $S/e$ then $e$ is a loop of of $S/f$.
\end{enumerate}
\item[(iv)]
An $E$-setoid is a special setoids so has well defined Tutte polynomial.
\end{enumerate}
Among special setoids one should also mention symmetric setoids,
that is setoid whose set of trees is invariant under any permutation
of edges. In fact a symmetric setoid is an $E$-setoid.

If $S$ is a finite setoid with ordered edges
then the polynomial can be always computed, using the computational
binary tree build according to the ordering of edges with leaves
being setoids with one tree (exactly as we did in the case of graphs).
Fig.2.1 shows an examples of a computation.
%Fig.1.1 (Examples of computation).
Thus for a setoid with an ordering $\rho$ of edges
we have the well defined Tutte polynomial $\chi_{\rho}(S)$.
We can however associate an invariant to a setoid in many ways so it
does not depend on orderings.
\begin{proposition}
We have the following invariants of a finite setoid $S=(E,T)$.
\begin{enumerate}
\item[(1)]
The set of polynomials (with possible repetitions) $\{\chi_{\rho}(S)\}$
over all ordering of edges, $\rho$.
\item[(2)]
The greatest common divisor of polynomials from (1).
\item[(3)] The smallest common multiple of polynomials from (1).
\item[(4)] The ideal $I_S$ in $Z[x,y]$ generated by polynomials from (1).
\item[(5)] The algebraic set associated to $I_S$.
\item[(6)] The coordinate ring of the algebraic set of (5).
\item[(7)] The Tutte polynomial $\chi_{sym}(S)= \chi(S_{sym})$
that is a polynomial of the symmetrization $S_{sym}$ of the setoid $S$
defined by $S_{sym}=(E,T_{sym})$ where $T_{sym}$ is the smallest
set of trees containing $T$ and invariant under permutations of edges;
see an example of the computation
for $S=(\{e_1,e_2,e_3\},\{\{e_1\},\{e_2,e_3\}\})$ in Fig.2.1.
\end{enumerate}
\end{proposition}
\centerline{{\psfig{figure=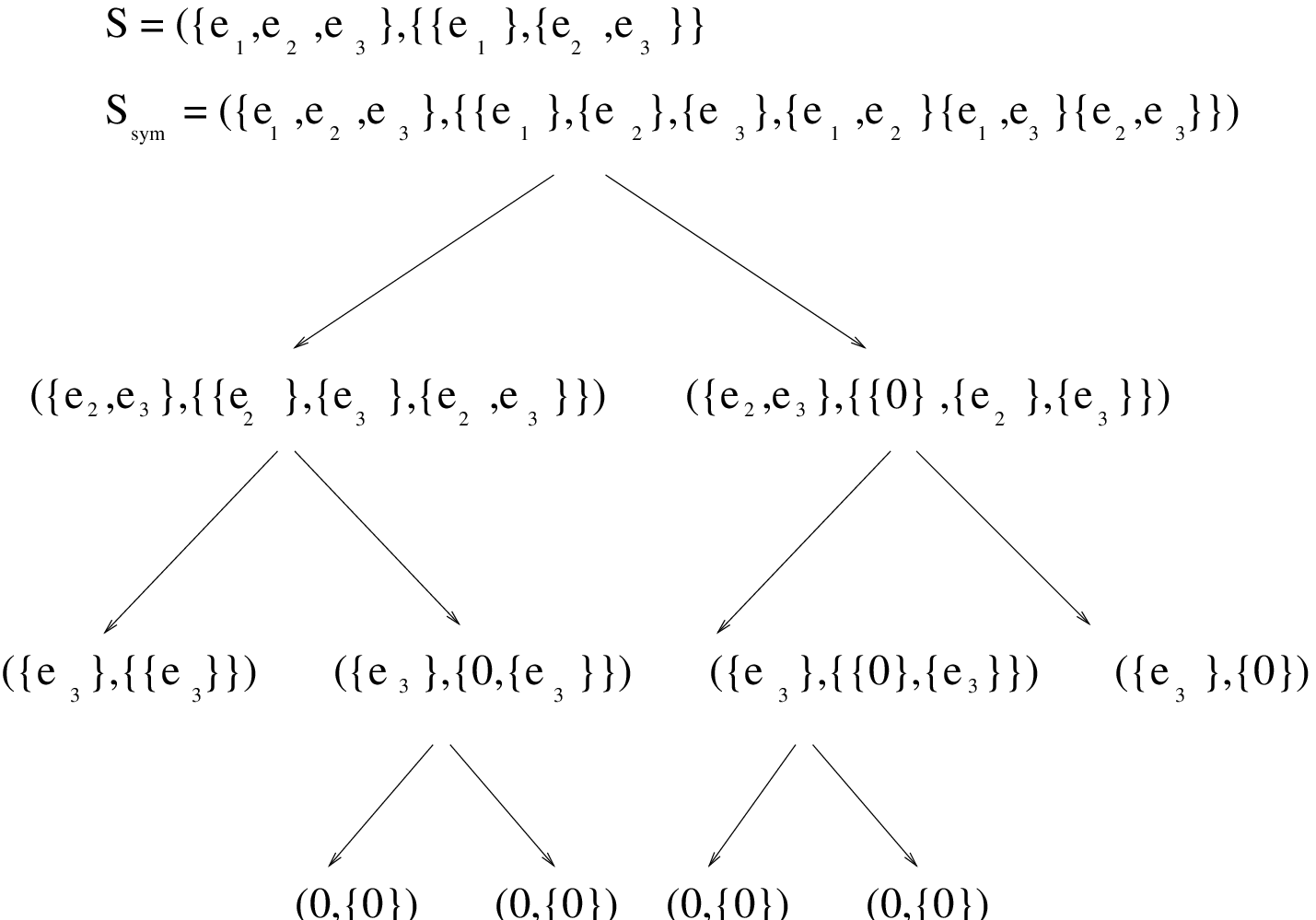,height=9.5cm}}}
\ \\
\centerline{Fig. 2.1\ $\chi_{sym}(S)= 4 + x + y$.}

We define a sum of setoids $S \cup S'$ as $(E\cup E', T\cup T')$ and
a product of setoids $S \star S'$ as $(E\cup E', T\times T')$.
We say that a setoid $S$ is 2-connected if it cannot be 
 obtained as a product of two setoids, $S=S' \star S''$
with at least one edge in each factor.

\begin{example} The setoid $S=(\{e_1,e_2\},\{\{e_1\},\{e_2\}\})$
is 2-connected.
\end{example}
\begin{exercise}
\begin{enumerate}
\item[(i)]
Show that if $S$ and $S'$ are $E$-setoids then  
$S \star S'$ is an $E$-setoid.
\item[(ii)]
Find conditions for a setoid which suffice to have:\
If $S$ is a 2-connected setoid of at least 2-edges
then $S-e$ or $S/e$ is 2-connected. What about $E$-setoids? 
Matroids?\footnote{Matroid is an $E$-setoid with all trees of
the same cardinality.}
\end{enumerate}
\end{exercise}

We will present now the second generalization of the deletion-contraction
method. We follow \cite{P-34} which in turn has been motivated by
Vassiliev-Gusarov invariants of knots \cite{P-9}. 
We present the idea for graphs,
however generalization to setoids is not difficult.

First we sketch the idea: we consider formal linear combinations of
finite graphs, $R{\cal G}$ with coefficients in a ring $R$.
Of course graphs (elements of ${\cal G}$) form a basis of $R{\cal G}$.
We introduce another basis of $R{\cal G}$  and then express
graphs as linear combinations of elements of the new basis. Coefficients
of this sum are graph invariants. Then we consider a filtration
of $R{\cal G}$ given by the new bases. The filtration allows us
to construct a Hopf algebra structure on the completion of $R{\cal G}$
with respect to the filtration.

Let $\cal G$ be the set of all finite graphs (up to isomorphism) and $R$
denote any commutative ring with unit.
Let $R\cal G$ denote the free $R$ module with
basis $\cal G$. We will enlarge the set of graphs $\cal G$
to
$\cal G'$ and add relations $\sim$ in $R\cal G'$ so
it reduces back to $R\cal G$.
Namely, $\cal G'$ is a set of graphs with two types of edges: classical and
special (or singular). $\cal G$ embeds in $\cal G'$ by interpreting elements
of $\cal G$ as having only classical edges.  Now, consider
in $R\cal G'$ relations, $\sim$, resolving special edges:

$G(e_s)= G(e)-(G-e)$ where $G(e)$ is a graph with a
classical edge $e$ and $G(e_s)$ is obtained from $G(e)$ by changing $e$ to
a special edge $e_s$. $G-e$ denotes, as usually, the graph obtained from $G(e)$
by deleting $e$.

Of course the embedding $R{\cal G} \subset R\cal G'$ induces the
$R$-isomorphism between $R\cal G$ and $R{\cal G'}/{\sim}$ and we will
usually identify these two modules. We just enlarged our $R{\cal G}$
by allowing graphs with special edges and then we express graphs
with special edges as linear combinations of classical graphs.
Of course special graphs ${\cal G}^s$ form also a basis of
$R{\cal G} = R{\cal G'}/{\sim} = R{\cal G}^s$, so we can use the base
change from ${\cal G}$ to ${\cal G}^s$. An algebraic structure (e.g.
bialgebra) simply expressed in ${\cal G}^s$ basis can look complicated
in ${\cal G}$ basis. In fact it allows us very simple interpretation
of an important Hopf algebra of Rota and Schmitt \cite{Schm-1}.

\begin{lemma}\label{Lemma IV.2.7}\
\begin{description}
\item [(a)]
(Change of basis). Express a graph $G \in \cal{G}$
as a linear combination of special graphs:
$$G = \sum_{H\in {\cal G}^s} a_H H$$
Then $a_H$ is equal to the number of embeddings of
$H$ in $G$ (embeddings which are bijections on vertices; a type of edges
is ignored). If we think about graphs in ${\cal G}^s$
as variables than the above formula can be called a pattern polynomial
of the graph $G$.
\item [(b)] Let $e_1,e_2,...,e_m$ be edges of a graph $G\in \cal{G}$. Then
$$G = \sum_{S\subset E(G)} G^{e_1,e_2,...,e_m}_{\epsilon_1,\epsilon_2,...,
\epsilon_m}$$
 where $\epsilon_i = 0$ or $-1$ 
($\epsilon_i = 0$ if $e_i \in S$ and  $-1$ otherwise)
and $G^{e}_{\epsilon}$, for
$\epsilon = 1$, $0$ or $-1$, denote three graphs in ${\cal G}'$
which differ only at the edge $e$ which is classical
for $\epsilon = 1$, special for $\epsilon = 0$ and deleted
for $\epsilon = -1$. 

\end{description}
\end{lemma}
\begin{proof}
Formula (b) follows by applying the formula
$G(e)= G(e_s)+(G-e)$ to every  edge of $G$; compare Example 2.8. \\
(a) is the interpretation of (b).
We can also prove (a) by an induction on the number
of the classical edges in a graph:\ consider an 
$R$-homomorphism $f:R{\cal G}' \mapsto R{\cal G}^s$ given for $G$ 
in ${\cal G}'$ by the formula $f(G)=\sum_{H\in {\cal G}^s}f_H(G)H$ 
where $f_H(G)$ is the number
of embeddings of $H$ in $G$ (bijective on vertices) with $H$ containing
all special edges of $G$. $f$ restricted to $R{\cal G}^s$ is therefore the
identity and it is immediate to check that $f_H(G(e))=f_H(G(e_s))+f_H(G-e),$
hence $G(e_s)-G(e)+(G-e)$ is in the kernel of $f$. Thus $f$ is an epimorphism
which descends to $f': R{\cal G}'/{\sim} \mapsto R{\cal G}^s.$
Because ${\cal G}^s$ generates $R{\cal G'}/{\sim}$, hence $f'$ is an 
$R$-isomorphism.
\end{proof}

\begin{example}\label{Example 1.3}
Consider two different connected graphs with 3-edges, $K_1$ and $K_2$.
We can use the relation $G(e)=G(e_s)+(G-e)$ to express $K_1$ and $K_2$ in
terms of graphs with only special edges. The binary computational resolving
tree and the result of the computation are shown in Figure 2.2.
\centerline{{\psfig{figure=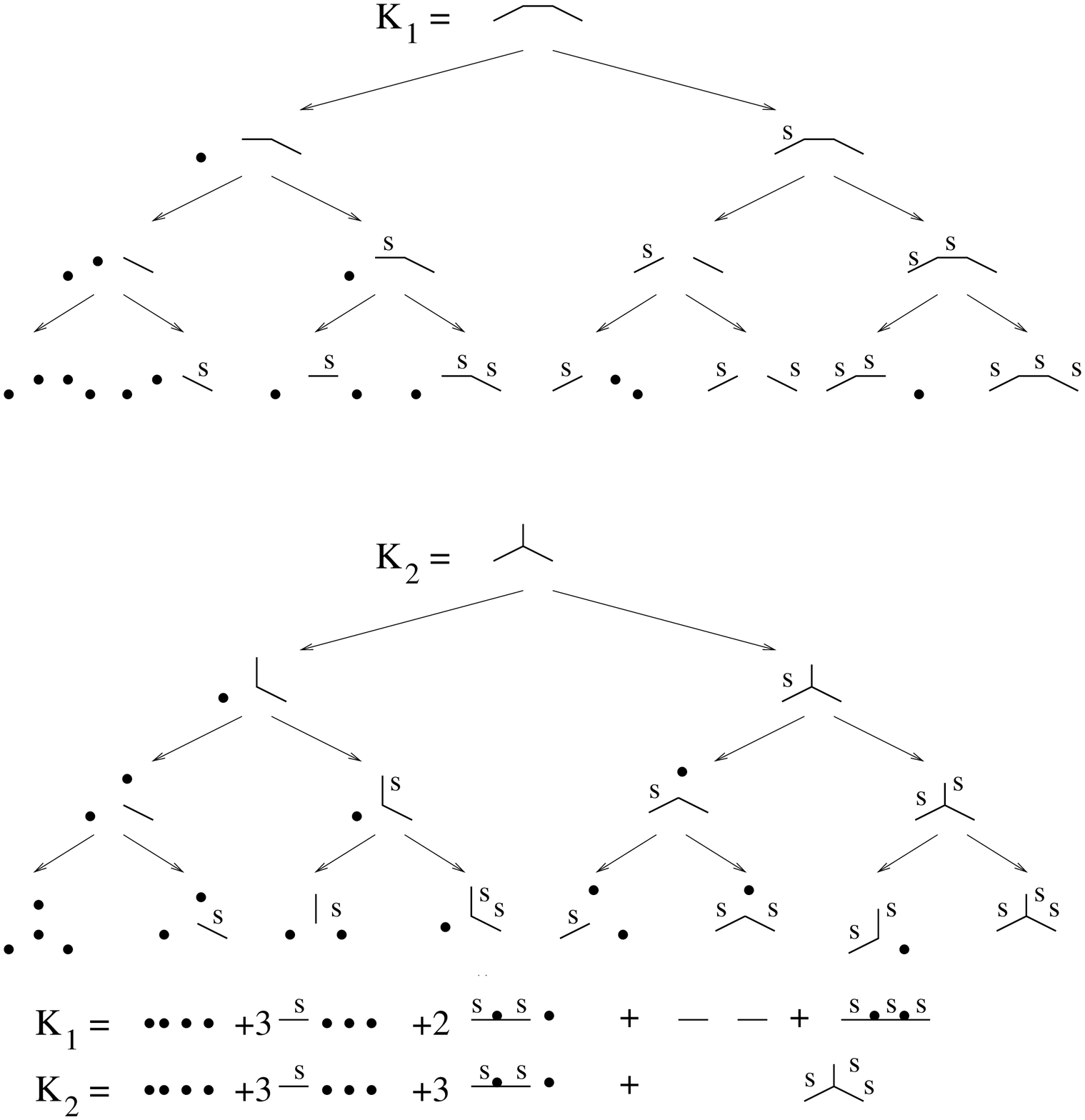,height=10.9cm}}}
%\vspace*{3in}
\begin{center}
Fig. 2.2
\end{center}

\end{example}

\begin{corollary}\label{Corollary 2.9}
$R{\cal G}$ is a ring with the disjoint sum as a product (more of it
after Def. 2.10). All invariants described below are ring homomorphisms
from $R{\cal G}$ to rings of polynomials.
\begin{description}
\item [ (a)] 
Let $<G>'_{\mu,A,B}= \mu B^{-|E(G)|}<G>_{\mu,A,B}$ be a version of 
the Kauffman bracket polynomial
of a graph $G\in {\cal H}$ (see Def. IV.1.x). 
Then the function $<>': R{\cal G}\to R[\mu,A,B]$
is given for a special graph $G^s$ by $<G^s>'_{\mu,A,B}= 
\mu^{p_0(G^s)+ p_1(G^s)}(AB^{-1})^{|E(G^s)|}$.
\item [(b)] 
Let $Q(G;t,z)$ be a dichromatic polynomial of a graph $G$ (see
\cite {Tral-2} for example) then the function $Q: R{\cal G}\to R[t,z]$ is
given for a special graph $G^s$ by $Q(G^s)=t^{p_0(G^s)}z^{p_1(G^s)}$
where $p_0(G^s)$ is the number of components and $p_1(G^s)$ the cyclomatic
number of $G^s$,
\item
[(c)] Let $Z(G)$ be the version of the dichromatic polynomial used for
example in \cite{K-10}, i.e. $Z: R{\cal G}\to R[q,v]$ is $R$-algebras
homomorphism satisfying $Z(\cdot)=q$ and $Z(G)=Z(G-e)+vZ(G/e)$. Then
for a special graph $G^s$ one has $Q(G^s)=q^{p_0(G^s)}v^{|E(G^s)|}.$
>From this we can get the well known formula $Z(G)(tz,z)=v^{|V(G)|}Q(G;t,z)$,
where $|V(G)|$ is the number of vertices of the graph $G$.
\item
[(d)] Let the $R$-algebras homomorphism ${\cal M}:R{\cal G}\to R[y]$, 
associate to $G$ its matching polynomial
\cite {Far-3}, then ${\cal M}$ can be determined
by: ${\cal M}(\cdot)=1, {\cal M}({\cal Y})=y, {\cal M}(G^s)=0$
if $G^s$ is a connected special graph different than one vertex graph and
different than $\cal Y$ where $\cal Y$ is the connected
graph of one edge and two vertices (interval). One can think of the
matching polynomial as a projection of the pattern polynomial
$\sum_{H\in {\cal G}^s} a_H H$ (Lemma IV.2.7(a)). $H$ projects to
$y^n$ if $H$ is composed of $n$ disjoint interval (and, possibly,
isolated vertices). Otherwise $H$ projects to $0$.
\end{description}
\end{corollary}

In the last part of the section we construct Hopf algebras out of graphs.
To have this part self-contained we offer below a short overview of
completions and Hopf algebras.\\
\ \\
{\bf Completions}\\
Let $M$ be a module over a commutative ring with identity $R$.
Consider a filtration $...C_3 \subset C_2 \subset C_1 \subset C_0= M$
of $M$ that is a descending family of submodules of $M$.
We can equip $M$ with a pseudo-metric\footnote{ We do not require
that if $\rho (x,y)=0$ then $x=y$.} $\rho: M\times M \to R_+$, where
$R_+$ denotes  non-negative real numbers, such that 
$\rho (x,y)= \frac{1}{k}$ if $(x-y) \in C_k$ but $(x-y)$ is not in $C_{k+1}$. 
If $(x-y) \in C_k$ for any $k$ then we put $\rho (x,y)=0$.
The pseudo-metric $\rho$ yields a topology on $M$ 
called an adic topology) and it is a Hausdorff topology 
iff $\rho$ is a metric or equivalently $\bigcap_i C_i = \{0\}$.
The pseudo-metric $\rho$ is invariant under addition (i.e. $\rho (x,y) =
\rho (x+a,y+a)$ for any $a \in M$) and under multiplication by 
an invertible scalar $r\in R$ (generally $\rho (rx,ry) \leq \rho (x,y)$). 
One can show that $M$ is a topological module. Using the pseudo-metric 
one can define now a completion $\hat M$ of $M$ by adding to $M$ Cauchy 
sequences with respect to $\rho$ modulo the standard equivalence relation 
on Cauchy sequences.
This relation makes $\hat M$ a metric space (with a metric $\hat\rho$
yielded by $\rho$), and a topological module. We have 
a distance preserving map (homomorphism) 
$i: M \to \hat M$ with $ker\ i = \bigcap_i C_i $ (the set of points with
pseudo-distance zero from $0\in M$).  
We will say that $\hat M$ is a module yielded by a filtration $\{C_i\}$
of $M$. If ${\cal I}$ is an ideal in $R$ then we have ${\cal I}$-adic 
filtration of $R$: $...{\cal I}^3\subset {\cal I}^2 \subset {\cal I} 
\subset R$, and of $M$: $...{\cal I}^3M\subset {\cal I}^2M \subset 
{\cal I}M \subset M$. $\hat M$ yielded by this filtration will be 
called ${\cal I}$-adic completion of $M$. In this
case $\hat M$ can be thought as a topological module over the topological
ring $\hat R$ where $\hat R$ is an ${\cal I}$-adic completion of $R$.

If $D_n$ is defined a a quotient $M/C_{n+1}$ then the completion 
$\hat M$ can be defined as an inverse limit of the 
sequence of $R$-epimorphisms
$...\rightarrow D_k \rightarrow D_{k-1} \rightarrow ... \rightarrow D_0
\rightarrow \{1\}$.

The simplest example of completion is that of polynomials $R[x_1,x_2,...]$
by infinite series $R[[x_1,x_2,...]]$; here filtration is given
by polynomials of degree no less than $i$. If we think of polynomials
as a ring than we have ${\cal I}$-adic filtration where ${\cal I}$ is 
an ideal generated by $x_1,x_2,...$.

We are ready now to describe the dichromatic filtration and completion
of graphs. 

\begin{definition}\label{Definition 1.1} \
\begin{description}
\item
[(a)]  (Dichromatic filtration): Let $C_k $ be the submodule of $R\cal G$
( $=R{\cal G'}/\sim$) generated by graphs with $k$
special edges.  The family $\{C_k\}$ forms a filtration of $R\cal G$
$$R{\cal G}=C_0 \supset C_1 \supset C_2 \supset ...\supset C_k \supset ...$$
and the filtration yields an adic topology on $R\cal G$. In particular
$\{C_i\}$ forms a basis of open sets around $0$.
\item
[(b)] The k'th dichromatic module of graphs is defined to be
$D_k=R{\cal G}/C_{k+1}$.
\item
[(c)] The dichromatic module of graphs, $\widehat{R\cal G} =D_{\infty}$, 
is defined to be the completion of $R\cal G$ yielded by the
filtration $\{C_k\}$.
\end{description}
\end{definition}
$R\cal G$ has has a natural $R$-algebra structure.
Namely, we can introduce a multiplication in $\cal G$ by taking as
$G_1\circ G_2$  the disjoint sum of $G_1$ and $G_2$. We obtain in such
a way a commutative semigroup.
$R\cal G$ is a semigroup ring (it has been introduced by W.T.Tutte
in 1947 \cite{Tut-1}). If we allow the empty graph $T_0=\emptyset$ then we
have also a unit of the multiplication. Then
$R\cal G$ is a semigroup algebra.  It is a filtered algebra 
because $C_i \circ C_j \subset C_{i+j} $ and therefore its 
completion is a (topological) algebra as well.
The following lemma describe some elementary but essential properties
of $R\cal G$ and its completion, $D_{\infty}$.
\begin{lemma}\label{Lemma x.x}\
\begin{description}
\item
[(a)] $R\cal G$ embeds in $D_{\infty}$
\item
[(b)] $D_{\infty}$ is a formal power series algebra in variables: connected
special graphs. $R\cal G$ is its dense subalgebra.
\item
[(c)] If $G_1$ and $G_2$ are two classical graphs with the same
number of vertices then $(G_1 - G_2) \in C_1$. In particular if $T_n$
is a graph with $n$ vertices and no edges an $|V(G)|=n$ then
$(G - T_n)\in C_1$.
[(d)] Let $\bar G_k = T_n + (T_n- G) + (T_n- G)^2 + ... + (T_n- G)^k$.
Then $G\circ \bar G_k = T_n^{k+1} - (T_n- G)^{k+1}$ and
$G\circ \bar G_k - T_n^{k+1} \in C_{k+1}$.
\end{description}
\end{lemma}
\begin{proof}
\begin{description}
\item [(a)] 
$R\cal G$ embeds in $D_{\infty}$ because $\bigcap C_i = \{0\}$. 
The last equality holds because by the change of basis
lemma for any element $a \in R\cal G$  there exists 
$i$ such that $a \notin C_i$
($C_i$ is a free module with basis: special graphs of at least $i$ edges).
\item [(b)] 
It follows from the change of basis lemma that the algebra 
$R\cal G$ can be identified with polynomial algebra in variables: 
connected special graphs. Thus (b) follows 
because the formal powers series algebra is the completion of the
polynomial algebra.
\item [(c)] 
If we construct a computational tree for the "pattern" polynomial
than exactly one leaf has no edges and it is $T_n$. Thus $(G - T_n)\in C_1$
and part (c) follows. 
\item [(d)] 
It is a standard "geometric series" formula. It will be
very useful later in constructing inverse to $G$ in the completion
(assuming $T_n$ invertible).
\end{description}
\end{proof}
\ \\
\ \\
{\bf  Bialgebras and Hopf algebras}
\\
\ \\
Let $A$ be a commutative ring with identity.  Let
  $B$ be an $A$-module with two $A$-module morphisms
  $i:A\to B$ and $\mu :B \otimes_A B\to B$.
  We say that $(B,\mu ,i)$ is an $A$-algebra if
  \begin{enumerate}
  \item[(a)] $\mu $ is associative, i.e.,
  $\mu (1\otimes \mu )=\mu (\mu \otimes 1)\ ;$ see Fig.2.3,
  \item[(b)] the unitary property holds, i.e.,
  $$
  (\mu (i\otimes 1))(a\otimes b) = ab = ba = (\mu
  (1\otimes i)) (b\otimes a)
  $$
  for any $a\in A$ and $b\in B$, where $1=1_B$ is the
  identity morphism on $B$; see Fig.2.4.  $\mu $ is called the
  multiplication map and $i$ the unit map.
  \end{enumerate}

  We define an $A$-coalgebra $B$ dually to an $A$-algebra:\\
Let $B$ be an $A$-module with two $A$-module morphisms
  $\epsilon :B\to A$ and $\nabla : B\to B\otimes_A B $.  
We say that ($B,\nabla ,\epsilon $) is an
  $A$-coalgebra if
\begin{enumerate}
  \item[(a)]  $\nabla $ is coassociative, i.e., $(\nabla
  \otimes 1)\nabla = (1\otimes \nabla )\nabla $ (see Fig.2.3), and
  \item[(b)]  the counitary property holds (see Fig.2.4), i.e.,
  $$
  (\epsilon \otimes 1) \nabla (b)=1\otimes b,\quad
  b\otimes 1=(1\otimes \epsilon )\nabla (b)\qquad
  \ \ for \ \ b\in B\,.
  $$
  $\nabla $ is called the comultiplication map and
  $\epsilon $ the counit map.
  \end{enumerate}
  \vskip2ex

\centerline{{\psfig{figure=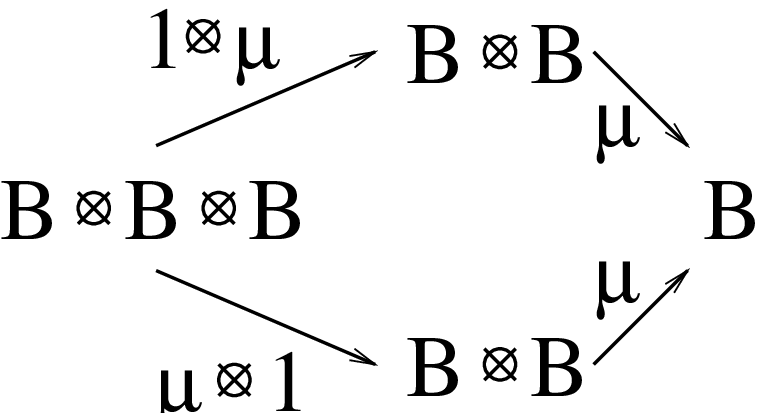,height=2.8cm}}\ \ \ \ \ 
{\psfig{figure=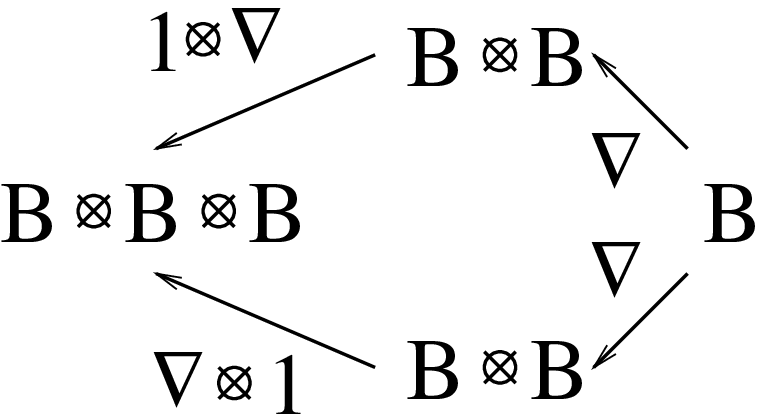,height=2.8cm}}}
\begin{center}
Fig.~2.3.\ Associativity and coassociativity.
\end{center}
\ \\
\centerline{{\psfig{figure=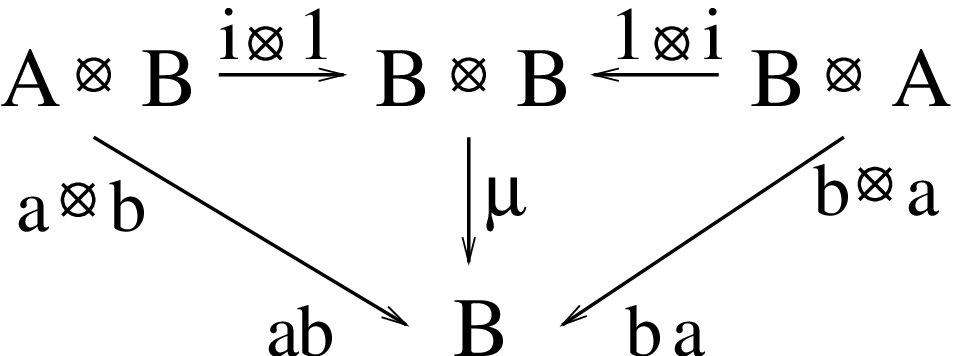,height=2.5cm}}\ \ \ \ \ \  
{\psfig{figure=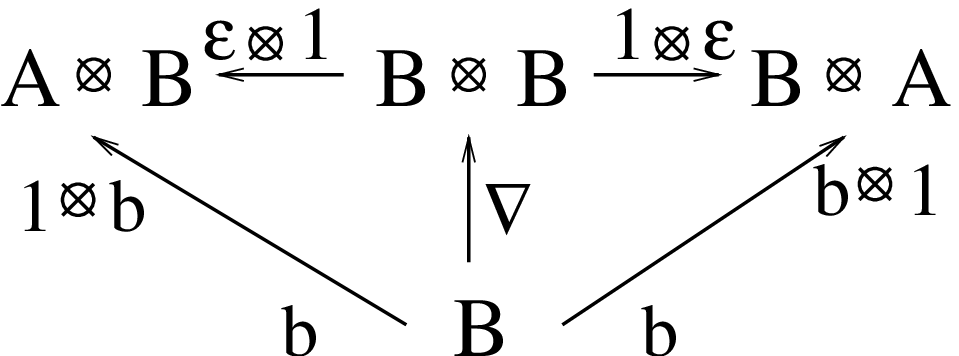,height=2.5cm}}}
\begin{center}
Fig.~2.4.\ Unitary and counitary properties.
\end{center}
  \begin{definition}\label{ 1.1}  Suppose that $(H,\mu ,i)$
  is an $A$-algebra and $(H,\nabla ,\epsilon )$ is an
  $A$-coalgebra.  If $\nabla $ and $\epsilon $ are
  $A$-algebra morphisms then $(H,\mu ,i,\nabla ,\epsilon)$ 
or simply $H$ is called an $A$-bialgebra.  The
  multiplication $\mu _{H\otimes H} : (H\otimes H)\otimes
  (H\otimes H)\to H\otimes H$ is given by the formula
$\mu
  _{H\otimes H}((a\otimes b)\otimes (c\otimes d)=\mu
  (a\otimes c
  )\otimes \mu (b\otimes d))$.
   The condition that $\nabla : H\to H \otimes H$ is an
  $A$-algebra homomorphism can be written as $\nabla \mu
 =
  (\mu \otimes \mu )(1\otimes P \otimes 1)(\nabla \otimes
  \nabla )$ where $P:H\otimes H\to H\otimes H$ is the
  $A$-module isomorphism defined by $P(a\otimes
 b)=b\otimes
  a$ and $\nabla (i(1))= i(1)\otimes i(1)$.  We will
 often
  identify $i(1)$ with 1 in further considerations.
 \end{definition}
  \begin{definition}\label{1.2}
Let $(H,\mu ,i,\nabla,\epsilon )$ be an $A$-bialgebra if there is an
  $A$-module morphism $S:H\to H$ such that $\mu (S\otimes
  1)\nabla =\mu (1\otimes S)\nabla =i\epsilon $ (i.e.,
  diagram 2.5. commutes) then $S$ is called the antipode of
  $H$ and an $A$-bialgebra with an antipodes is called an
  $A$-Hopf algebra.
\end{definition}
\centerline{\psfig{figure=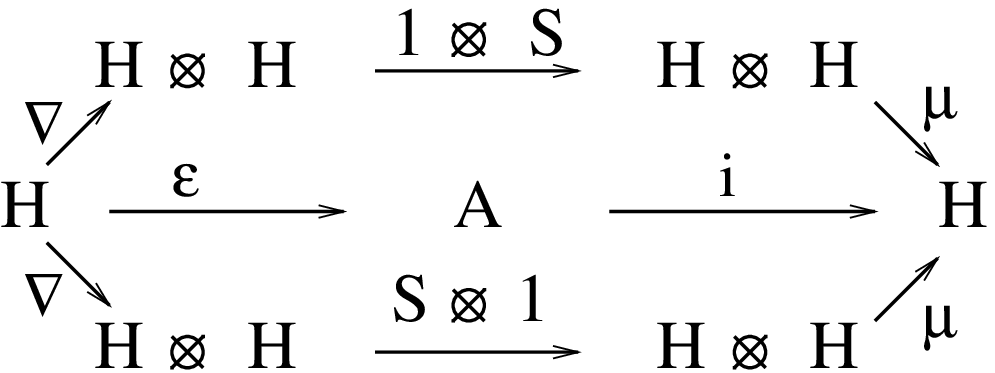,height=4.8cm}}
\begin{center}
Fig.~2.5.\ Antipode property.
\end{center}
  One can easily check that an antipode $S$ is an
  anti-$A$-algebra morphism, that is, $S(i(1)) = i(1)$
 and
$S(\mu (a\otimes b)) = \mu (S(b)\otimes S(a))$.  (See
  \cite{Ab} Thm. 2.1.4).

  An element $h\in H$ such that $\epsilon (h)=1$ and
  $\nabla (h) = h\otimes h$ is called a
  group-like element.

  An element $h\in H$ such that $\nabla (h) = h\otimes
  1+1\otimes h$ is called a primitive element of $H$.  If
  $h$ is a primitive element of $H$ then $\epsilon (h)=0$
  (see \cite{Ab} Thm. 2.1.3). If $S^2=1_H$ then $H$ is called an
involutive Hopf algebra.

The simplest, and relevant, example of a bialgebra is a semigroup
algebra $RG$ for a semigroup with identity $G$ and a commutative
ring $R$. Every element of $G$ is a group like element, that is
$\nabla (g) = g\otimes g$ and $\epsilon (g)=1$. If $G$ is a group
then $RG$ is an involutive Hopf algebra with the antipode $S(g)=g^{-1}$.

We will consider two bialgebra structures on the module of graphs 
$R\cal G$, and show that they are isomorphic. 
Then we extend the structure by the completion
to a Hopf algebra.

$R\cal G$ has a standard semigroup bialgebra structure. It is described by:
\begin{description}
\item
[(a)] unite $i: R \rightarrow R\cal G$ is given by $i(r)=r\emptyset$,
\item
[(b)] multiplication $\mu :R{\cal G} \otimes R{\cal G} \rightarrow R\cal G$
is given by $\mu (G_1,G_2)=G_1\circ G_2$ (disjoint sum),
\item
[(c)] counit $\epsilon : R{\cal G} \to R$ is given by $\epsilon (G)
=1$
\item
[(d)] comultiplication $\nabla : R{\cal G} \to R{\cal G}\otimes R{\cal G}$
is given by $\nabla (G) = G\otimes G$.
\end{description}

Consider $R$ with the discrete topology and $R{\cal G}\otimes R\cal G$ with
the topology yielded by the filtration $\sum_{i=0}^kC_i\otimes C_{k-i}$ (or
equivalently by the filtration
$\{R{\cal G}\otimes C_k +C_k\otimes R\cal G\}$).
\begin{lemma}\label{Lemma 2.4}\
\begin{description}
\item
[(a)]  $R\cal G$ is a topological bialgebra,
\item
[(b)] let $t$ denote the one vertex graph (then $t^n=T_n$ is the n-vertex
graph with no edges). Consider the extension of $R\cal G$ by 
$t^{-1}{\cal G}$ 
(we just make $t$ invertible in the algebra; we do not kill anything 
because $\cal G$
is a semigroup with the unique prime decomposition). 
The completion $D_{\infty}' =
D_{\infty}{\otimes}_{R[t]} R[t^{\pm 1}]$ is a (topological) Hopf algebra.
\end{description}
\end{lemma}

\begin{proof}
(a)
(i) $C_i \circ C_j \subset C_{i+j} $, therefore ${\mu}^{-1}(C_k)\supset
     (C_k \otimes R{\cal G} + R{\cal G} \otimes C_k)$. Thus the multiplication
is continuous,

(ii) $\nabla C_k \subset \sum_{i=0}^{k}C_i\otimes C_{k-i} .$ Thus the
comultiplication is continuous (we will give the exact description of
$\nabla C_k $ in the proof of Theorem 2.17),

(iii) Counit is continuous with discrete topology on the ring $R$ (so
any topology on $R$) because ${\epsilon}^{-1}(0)$ contains $C_1$,

(iv) Unit map is continuous for $R$ with discrete topology.

(b) An element $x$ is invertible in the completion iff it is of the form
$e+c$ where $e$ is invertible and $c\epsilon C_1$. In our case, for an
$n$-vertex graph $G$ we have $(G-t^n) \epsilon C_1$. Therefore in order to
invert any graph $G$ we need $t$ to be invertible. On the other hand
$G$ are group like elements of our bialgebra and they generate it. So
if any $G$ is invertible we can  define the antipode map $S:D_{\infty}' \to
D_{\infty}'$ by $S(G)= G^{-1}= t^{-n}(1+(1 - t^{-n}K)+(1 - t^{-n}K)^2 +
(1 - t^{-n}K)^3+...)$. In such a way we define the antipode map on the dense
subset and then extend it continuously to the whole completion.
\end{proof}
\begin{corollary}\label{Corollary 1.5}\ \\
If $G\in \cal G$ and $G$ has $n$ vertices then for any $k$ there is an element
$d\in D_k$ such that $G\circ d = t^{n(k+1)}$ in $D_k$. In particular $G$ is
invertible in $D_k'=D_k{\otimes}_{R[t]} R[t^{\pm 1}]$;
 compare Lemma IV.2.11(c).
\end{corollary}
Consider another, simpler, filtration $\{E_i\}$ of $R\cal G$, where
$E_i$ is generated by graphs with at least $i$ edges. Let $T_2$ denote
the topology on $R\cal G$ yielded by the filtration. One can immediately
check that the bialgebra $(R{\cal G},i,\mu,\epsilon ,\nabla)$ is a topological
bialgebra with respect to $T_2$. However, its completion does not possess
an antipode map (i.e. is not a Hopf algebra) unless we extend it by inverses
of all graphs in $\cal G$. There is, however, another bialgebra structure
on the semigroup algebra $R\cal G$ (of which I learned from Schmitt 
\cite{Schm-1,Schm-2}),
completion of which (with respect to $T_2$) is a Hopf algebra.
\begin{lemma}(Schmitt)\label{Lemma 1.6}\
\begin{description}
\item
[(a)]
The semigroup algebra $R\cal G$ is a bialgebra with ${\epsilon}'$ and
${\nabla}'$ defined as follows:
${\epsilon}'(G)= 1$ if $G$ has no edges and $0$ otherwise,
${\nabla}'(G)=\sum_{S_1 \cap S_2 = \emptyset}(G-S_2)\otimes (G-S_1)$ 
where the sum is taken over all ordered disjoint pairs of subsets of 
edges of $G$, 
\item
[(b)] The above bialgebra is a topological bialgebra with respect to $T_2$,
\item
[(c)] The completion of $(R{\cal G},T_2)$ is a Hopf algebra (assuming that the
one vertex graph, $t$, is invertible).
\end{description}
\end{lemma}

\begin{proof}
We will show that the bialgebras $(R{\cal G},i,\mu,\epsilon ,\nabla;T_1)$ and
$(R{\cal G},i,\mu,{\epsilon}' ,{\nabla}';T_2)$ are isomorphic by
a homeomorphism.
\end{proof}
\begin{theorem}\label{Theorem 2.17}\ \\
Let $\phi : R{\cal G} \to R{\cal G}$ be an $R$-linear map given by
$\phi(G) = G^s$ where $G\epsilon \cal G$ and
$G^s$ is the special graph obtained from $G$ by changing all its edges to
special edges. \ \ Then $\phi$ is a homeomorphism
and an isomorphism of bialgebras $(R{\cal G},i,\mu,{\epsilon}' ,{\nabla}';T_2)$
and  $(R{\cal G},i,\mu,\epsilon ,\nabla;T_1)$
\end{theorem}

\begin{proof}
$\phi (E_i) = C_i$, thus $\phi$ is a homeomorphism. It is an $R$-algebras
isomorphism, essentially by definition. Also by definition ${\epsilon}'=
\epsilon \phi.$
It remains to analyze $\nabla(G^s)$ where $G^s$
is a special graph with $n$ edges (i.e. $G^s \in C_n$).
First, we illustrate it using the special one edge graph $e_s$:
$\nabla (e_s)=\nabla (e-t^2)= e\otimes e - t^2\otimes t^2=
e\otimes (e-t^2)+(e-t^2)\otimes t^2=
e\otimes e_s + e_s\otimes t^2=t^2\otimes e_s +
e_s\otimes e_s + e_s\otimes t^2$. Inductively we reach the general formula
$$\nabla G^s=\sum_{{\epsilon}_1,{\epsilon}_2,...,{\epsilon}_n}
G^{e_1,e_2,...e_n}_{{\epsilon}_1,{\epsilon}_2,...,{\epsilon}_n} 
\otimes G^{e_1,e_2,...e_n}_ {{\epsilon}_1-1,{\epsilon}_2-1,...,
{\epsilon}_n-1}=$$
$$\sum_{S_1 \cap S_2 = \emptyset}(G^s-S_2)\otimes (G^s-S_1) \in 
\sum_{i}C_i \otimes C_{n-i},$$ 
where $G^s$ is a special graph of $n$ edges $e_1$,...,$e_n$, 
${\epsilon}_i$ is $1$ or $0$
and a sub-index $1$, $0$ or $-1$ under $e_i$ indicates whether 
we deal with a classical edge,
special edge or deleted edge, respectively.
Thus the formula is analogous to the Schmitt co-multiplication and $\phi$ is
a bialgebra isomorphism. 
The inductive step (for the first part of the formula) works as follows:\\
$$\nabla (G^{e_1,e_2,...e_n,e_{n+1}}_{0,\ 0,\ ...\ 0,\ 0})=
\nabla (G^{e_1,e_2,...e_n,e_{n+1}}_{0,\ 0,\ ...\ 0,\ +1}) -
\nabla (G^{e_1,e_2,...e_n,e_{n+1}}_{0,\ 0,\ ...\ 0,\ -1})$$ which
equals by the inductive assumption to:\\
$$ \sum_{{\epsilon}_1,{\epsilon}_2,...,{\epsilon}_n}
(G^{e_1,e_2,...e_n,e_{n+1}}_{{\epsilon}_1,{\epsilon}_2,...,{\epsilon}_n, +1}
\otimes 
G^{e_1,e_2,...e_n,e_{n+1}}_ {{\epsilon}_1-1,{\epsilon}_2-1,...
,{\epsilon}_n-1,+1}
- G^{e_1,e_2,...e_n,e_{n+1}}_{{\epsilon}_1,{\epsilon}_2,...,{\epsilon}_n, -1}
\otimes
G^{e_1,e_2,...e_n,e_{n+1}}_{{\epsilon}_1-1,{\epsilon}_2-1,...
,{\epsilon}_n-1,-1})
=$$
$$\sum_{{\epsilon}_1,{\epsilon}_2,...,{\epsilon}_n}
(G^{e_1,e_2,...e_n,e_{n+1}}_{{\epsilon}_1,{\epsilon}_2,...,{\epsilon}_n, +1}
\otimes 
(G^{e_1,e_2,...e_n,e_{n+1}}_ {{\epsilon}_1-1,{\epsilon}_2-1,...,{\epsilon}_n-1,
+1} - 
G^{e_1,e_2,...e_n,e_{n+1}}_{{\epsilon}_1-1,{\epsilon}_2-1,...,
{\epsilon}_n-1,-1}) \ +$$ 
$$(G^{e_1,e_2,...e_n,e_{n+1}}_{{\epsilon}_1,{\epsilon}_2,...,{\epsilon}_n, +} -
G^{e_1,e_2,...e_n,e_{n+1}}_{{\epsilon}_1,{\epsilon}_2,...,{\epsilon}_n, -1})
\otimes
G^{e_1,e_2,...e_n,e_{n+1}}_{{\epsilon}_1-1,{\epsilon}_2-1,
...,{\epsilon}_n-1,-1})=$$ 
$$\sum_{{\epsilon}_1,{\epsilon}_2,...,{\epsilon}_n}
(G^{e_1,e_2,...e_n,e_{n+1}}_{{\epsilon}_1,{\epsilon}_2,...,{\epsilon}_n, +1} 
\otimes
G^{e_1,e_2,...e_n,e_{n+1}}_ {{\epsilon}_1-1,{\epsilon}_2-1,...
,{\epsilon}_n-1,0} \ +$$
$$ G^{e_1,e_2,...e_n,e_{n+1}}_{{\epsilon}_1,{\epsilon}_2,...,{\epsilon}_n, 0}
\otimes
G^{e_1,e_2,...e_n,e_{n+1}}_ {{\epsilon}_1-1,{\epsilon}_2-1,...
,{\epsilon}_n-1,-1}) =$$
$$\sum_{{\epsilon}_1,{\epsilon}_2,...,{\epsilon}_n,{\epsilon}_{n+1}}
G^{e_1,e_2,...e_n,e_{n+1}}_{{\epsilon}_1,{\epsilon}_2,...,{\epsilon}_n,
{\epsilon}_{n+1}}
\otimes
G^{e_1,e_2,...e_n,e_{n+1}}_ {{\epsilon}_1-1,{\epsilon}_2-1,...
,{\epsilon}_n-1, {\epsilon}_{n+1} -1}$$ where $\epsilon_i \geq 0$. Thus
inductive step is performed. The second part of the equality in
the formula for ${\nabla}(G^s)$ follows by resolving all classical
edges on the left site of the tensor products in the formula.

%$$G= \sum_{S\subset E(G)} (-1)^{|S|}G^{e_1,e_2,...,e_m}_{\epsilon_1,\epsilon_2,
%..., \epsilon_m}$$ where $\epsilon_i = -1$ if $e_i \in S$ and  $1$ otherwise. 
\end{proof}

\begin{remark}\label{Remark 1.8}\ \\
We can consider an involution $\alpha:R{\cal G} \to R\cal G$ given by:
$\alpha (G)=(-1)^{|E(G)|}G$ where $|E(G)|$ is the number of edges of $G$. Then
$\phi \alpha$ is an involution of $R\cal G$; more generally $\alpha$ and $\phi$
generate the group $\{\alpha,\phi : {\alpha}^2=1, \alpha \phi \alpha
={\phi}^{-1}\}$, which is isomorphic to the group of isometries of integers.
\\
Notice that $\alpha(G_{e_s}) = \alpha(G_{e}) + \alpha(G_{e}-e)$.
\end {remark}
%We can consider two bialgebras "dual" to $R\cal G$ first the bialgebra
\begin{corollary}\label{Corollary 2.19}\ \\
Let $\cal I$ be an ideal in $R[q,v]$ generated by $v$, then we have
$\cal I$-adic filtration $\{{\cal I}^k\}$ of $R[q,v]$ and the $\cal I$-adic
completion, $\widehat {R[q,v]}$, of $R[q,v]$ (i.e. formal power series in $v$).
Then the dichromatic polynomial $Z(G)$ yields a filtered $R$-algebras
homomorphism and extends to the completions;
$\hat {Z}:\widehat{R{\cal G}}\to \widehat{R[q,v]}$ (compare 
Corollary 2.9(c)).
\end{corollary}
Our construction of a bialgebra and Hopf algebra can be extended to setoids
$S=\{E,T\}$ by extending edges $E$ by special edges $E^s$ and then
resolving them: $S_{e_s} = S_{e} - (S_{e}-e)$. The construction is
similar to that for the dichromatic module of graphs.
\\
\ \\
\ \\
%\end{document}

%\newpage

\section{Jones polynomials of alternating and adequate diagrams;
Tait conjectures}
\markboth{\hfil{\sc Graphs and links}\hfil} 
{\hfil{\sc Jones polynomial and Tait conjectures}\hfil}

More than a hundred years ago Tait was  
setting his tables of knots \cite{Ta}
using the following working assumptions on alternating
diagrams of links:

\begin{enumerate}
\item[(1)]
\begin{enumerate}
\item[(i)]
 A reduced alternating diagram of a given link has
a minimal number of crossings among all diagrams representing
the link. 
In particular, two  reduced alternating diagrams of the same link 
have the same number of crossings. A diagram is called reduced
if it has no nugatory crossings 
({\psfig{figure=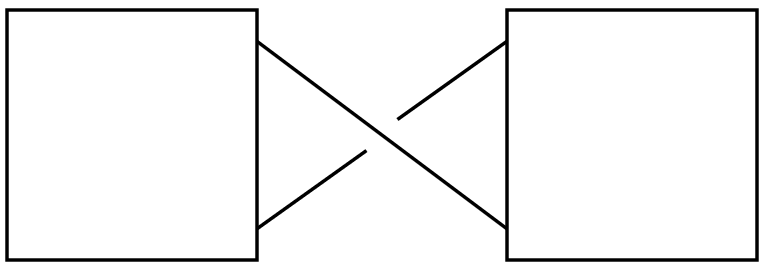,height=0.4cm}}).
\item[(ii)] If we assume additionally that our link is prime and non-split
(i.e. it is not a split or connected sum of links)
then any nonalternating diagram of the link has a non-minimal 
number of crossings.
\end{enumerate}
\item[(2)] Two oriented reduced alternating diagrams of the same link 
have the same Tait(or writhe) number (it is  defined
as a sum of signs of all crossings of the diagram and denoted 
by $Tait(D)$, $w(D)$ or $\tilde n(D)$).
\item[(3)]
 There exist easily recognizable moves
on alternating diagrams, called the Tait moves, or Tait flypes
 (see Fig.~3.1), such that two reduced alternating
diagrams of a given link can be reached one from the other by 
a sequence of such moves.
\ \\
\ \\
%\vspace*{1.5in}\centerline{\Psfig{figure=Rys.2.1}}
\centerline{{\psfig{figure=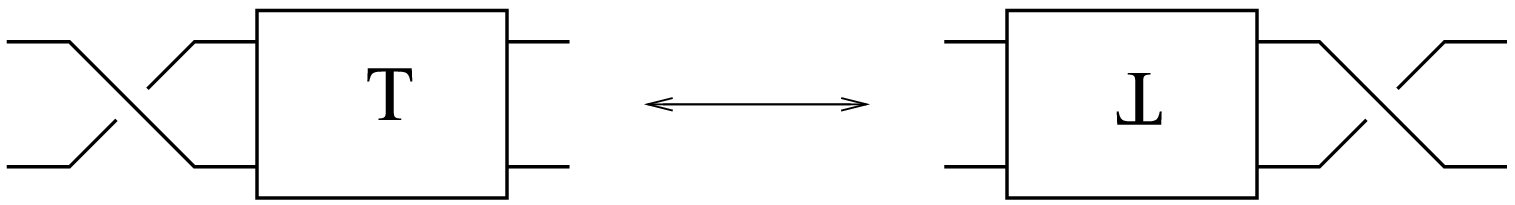,height=1.8cm}}}
\begin{center}
Fig.~3.1
\end{center}

\end{enumerate}

The above three ``assumptions'' are called {\it Tait conjectures}.
One of the most important applications of the Jones polynomial
is a proof of the first two of these conjectures by
Murasugi \cite{M-4,M-5}, Thistlethwaite \cite{This-3} and 
Kauffman \cite{K-6}. All these proofs apply a version of the Jones
polynomial which was discovered by Kauffman in the summer 
of 1985 \cite{K-6}. 
The third Tait conjecture
% is still open.  \footnote{It 
has been proved by Menasco and Thistlethwaite \cite{MT-1,MT-2}. 
The proof combines the use of Jones type polynomials with 
study of incompressible surfaces and goes beyond the scope of this 
book\footnote{Murasugi gave before an elementary proof for 
some special classes of alternating links.}.

Now we will describe the Kauffman's version of the Jones polynomial
and subsequently, following
Lickorish and Thistlethwaite \cite{L-T}, we will 
apply it to study adequate diagrams which are generalization
of alternating diagrams.
In particular, we will prove the first conjecture of Tait and outline 
a proof of the second.

The
Kauffman bracket polynomial $\langle L\rangle$ 
was defined by Kauffman in the summer of 1985 independently on the 
Jones polynomial and without relation to Tutte polynomial. Kauffman 
was investigating possibility that three diagrams 
\parbox{0.5cm}{\psfig{figure=L+nmaly.eps}}\ ,\ 
\parbox{0.5cm}{\psfig{figure=L+nmaly.eps}}\ and \ 
\parbox{0.5cm}{\psfig{figure=Linftynmaly.eps}} 
can be linked by a linear relation 
leading to a link invariant. Only later he realized that he constructed 
a variant of the Jones polynomial.

\begin{definition}\label{4:3.1}
Let $D$ be an unoriented diagram of a link. Then the Kauffman 
bracket polynomial
$\langle D \rangle \in Z[A^{\mp 1}]$ is defined by the following
properties:
\begin{enumerate}
\item[(i)] $\langle \bigcirc\rangle = 1$

\item[(ii)] $\langle\bigcirc\sqcup D\rangle = -(A^2+A^{-2})\langle L\rangle$

\item[(iii)] $\langle$ {{\psfig{figure=L+nmaly.eps}}}$\rangle =
A\langle {\mbox{{\psfig{figure=L0nmaly.eps}}}}
\rangle
+
A^{-1}\langle
{{\psfig{figure=Linftynmaly.eps}}}
\rangle$
\end{enumerate}
\end{definition}

\begin{proposition}\label{V.3.2}
\begin{enumerate}
\item[(i)] The Kauffman bracket polynomial is well defined, that is 
conditions (i)-(iii) define the unique function from the set of all 
diagrams, $\cal D$, to the ring of Laurent polynomials,
 $\langle \rangle : {\cal D} \to Z[A^{\pm}]$. 
\item[(ii)] The bracket $\langle D\rangle$ is an invariant of a 
regular isotopy of link diagrams, that is it is preserved 
by the second and third Reidemeister moves. Furthermore,
it is an invariant of a weak regular isotopy\footnote{Recall 
that two link diagrams are weak regular isotopic if they are 
related by the second and third Reidemeister moves and the weak 
first Reidemeister move in with two opposite Reidemeister moves,
$R_{+1}$ and $R_{-1}$ can be canceled 
(\parbox{0.5cm}{\psfig{figure=R+maly.eps,height=0.4cm}}
\parbox{0.4cm}{\psfig{figure=R-maly.eps,height=0.4cm}}
$\leftrightarrow$
%{\psfig{figure=R2-p.eps,height=0.4cm}}), 
\parbox{0.4cm}{\psfig{figure=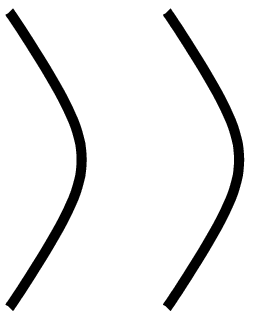,height=0.4cm}}), 
compare Lemma I.5.6. 
In \cite{A-P-R} the weak regular
isotopy is called balanced isotopy.}.
\end{enumerate}
\end{proposition}

\begin{proof}
\begin{enumerate}
\item[(i)] First we note that the value of  
$\langle D\rangle$ for a given diagram of a link does not depend
on the way we compute it.
This is a special case of Lemma \ref{4:1.4}, 
we repeat its proof once more, now
without using the graph associated to the diagram.

Let $\row{c}{n}$ denote crossings of the diagram $D$.

By the Kauffman state of $D$ we understand a function 
$s:\{i:1\leq i\leq n\}\rightarrow\{
-1,1\}$, that is every crossing has associated $+1$ or $-1$ and it will 
be treated depending on the sign. Let $D_s$ (or $sD$) denote 
the diagram obtained from 
$D$ according to the following rules:
in the crossing 
$c_i$ the diagram {\psfig{figure=L+nmaly.eps}}
is changed to
{\psfig{figure=L0nmaly.eps}}
if $s(i)=1$ and it is changed to
{\psfig{figure=Linftynmaly.eps}}
if $s(i) = -1$. 
Let $|D_s|$ (or simply $|s|$) denote the number of components of $D_s$, then
\begin{formulla}\label{4:2.2}
$$\langle D\rangle = \sum_{s\in 2^n} A^{\sum s(i)} (-A^2-A^{-2})^{|s|-1}$$
\end{formulla}

The formula  \ref{4:2.2} follows immediately from conditions
(i)---(iii) of Definition V.2.1 and it is a special case of 
Lemma V.1.4.
\item[(ii)]
We have yet to prove that $\langle D\rangle$ is not changed when we 
apply the second and the third Reidemeister moves --- this
concerns regular isotopy --- and also when we apply first weak 
Reidemeister move --- this in the case of a weak regular isotopy.

\begin{enumerate}
\item The second Reidemeister move.\\
\centerline{ $\langle \parbox{0.6cm}{\psfig{figure=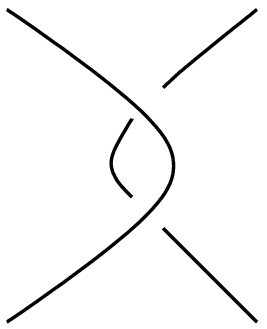,height=0.4cm}} 
\rangle=$}
\centerline{$A\langle\parbox{0.6cm}{\psfig{figure=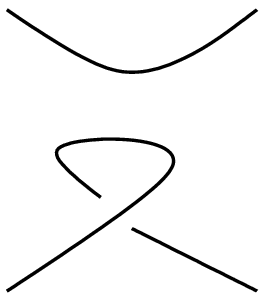,height=0.4cm}}
\rangle + A^{-1}\langle 
\parbox{0.6cm}{\psfig{figure=L-nmaly.eps,height=0.4cm}} \rangle =$}
\centerline{$A(A\langle {\psfig{figure=L0nmaly.eps,height=0.4cm}} \rangle +
A^{-1}\langle{\psfig{figure=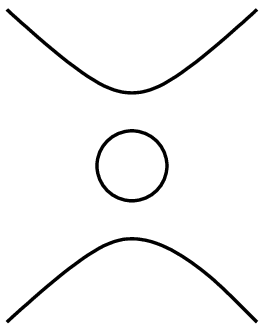,height=0.4cm}} \rangle)
+A^{-1}(A\langle{\psfig{figure=Linftynmaly.eps,height=0.4cm}} \rangle
+A^{-1}\langle {\psfig{figure=L0nmaly.eps,height=0.4cm}} \rangle)=$}
\centerline{$(A^2+AA^{-1}(-A^2-A^{-2})+A^{-2})\langle
{\psfig{figure=L0nmaly.eps,height=0.4cm}} \rangle
+\langle {\psfig{figure=Linftynmaly.eps,height=0.4cm}} \rangle=$}
\centerline{$\langle 
\parbox{0.5cm}{\psfig{figure=Linftynmaly.eps,height=0.4cm}} \rangle.$}

\item The third Reidemeister move.
\begin{eqnarray*}
&\langle {\psfig{figure=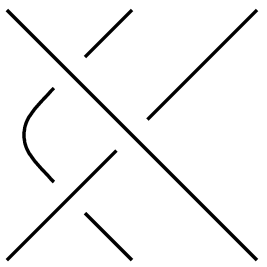,height=0.4cm}} \rangle& =\\
&A\langle
{\psfig{figure=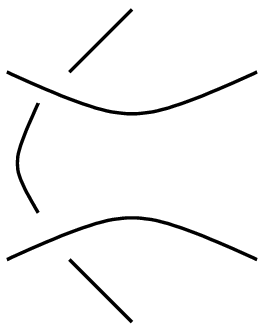,height=0.4cm}}
\rangle+A^{-1}\langle
{\psfig{figure=R2-k.eps,height=0.4cm}}\ {\psfig{figure=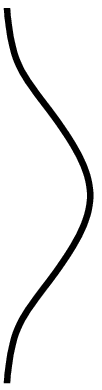,height=0.4cm}}
\rangle&=\\
=&A\langle
{\psfig{figure=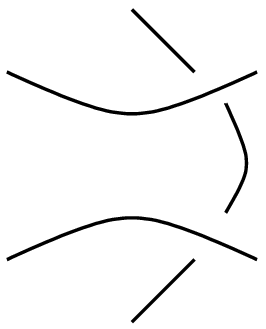,height=0.4cm}}
\rangle
+A^{-1}\langle
{\psfig{figure=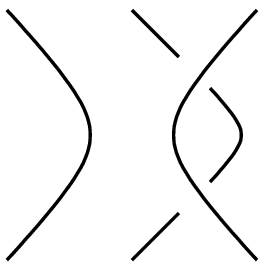,height=0.4cm}}
\rangle&=\\
=&\langle
{\psfig{figure=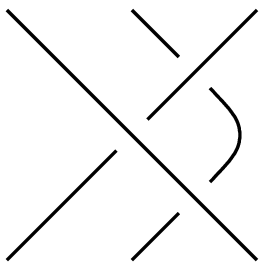,height=0.4cm}}
\rangle .&\\
\end{eqnarray*}
here we have used twice the invariance under the second
Reidemeister move.

\item First weak Reidemeister move.

Let us check at the beginning how first Reidemeister moves changes
$\langle D\rangle$:
\begin{formulla}\label{4:4:2.3}
\begin{eqnarray*}
&\langle
\parbox{0.6cm}{\psfig{figure=R+maly.eps}}
\rangle&=\\
&A\langle  {\psfig{figure=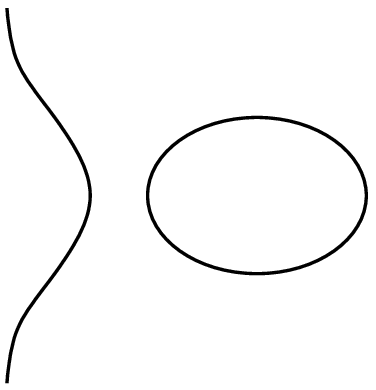,height=0.4cm}}
\rangle+A^{-1}
\langle  {\psfig{figure=ver.eps,height=0.4cm}}
\rangle&=\\
&(A(-A^2-A^{-2})+A^{-1})\langle {\psfig{figure=ver.eps,height=0.4cm}}
\rangle&=\\
&-A^{3}\langle  {\psfig{figure=ver.eps,height=0.4cm}}
\rangle&\\
&\langle \parbox{0.6cm}{\psfig{figure=R-maly.eps}}
\rangle&=\\
&A\langle  {\psfig{figure=ver.eps,height=0.4cm}} \rangle
+A^{-1}\langle  {\psfig{figure=v-o.eps,height=0.4cm}} \rangle
&=\nonumber\\
&(A+A^{-1}(-A^2-A^{-2}))\langle  {\psfig{figure=ver.eps,height=0.4cm}}
\rangle&=\\
&-A^{-3}\langle  {\psfig{figure=ver.eps,height=0.4cm}}
\rangle&\\
\end{eqnarray*}
\end{formulla}

Therefore for the first weak Reidemeister move if follows that:
$$\langle
\parbox{0.6cm}{\psfig{figure=R+maly.eps}},
\parbox{0.6cm}{\psfig{figure=R-maly.eps}}
\rangle = -A^3\langle {\psfig{figure=ver.eps,height=0.4cm}},
 \parbox{0.6cm}{\psfig{figure=R-maly.eps}} \rangle
=-A^3(-A^{-3})\langle 
{\psfig{figure=ver.eps,height=0.4cm}}, {\psfig{figure=ver.eps,height=0.4cm}}
\rangle=\langle 
{\psfig{figure=ver.eps,height=0.4cm}}, {\psfig{figure=ver.eps,height=0.4cm}}
\rangle$$
which completes the proof of Proposition 3.2.
\end{enumerate} 
\end{enumerate}
\end{proof}

After a slight modification the Kauffman bracket 
will give an invariant of global isotopy:

\begin{lemma}\label{4:2.4}
%\begin{lemma}\label{V.3.5}
\ \\
\begin{enumerate}
\item Let $\mbox{sw}(D)$ denote the algebraic self-crossing number
of the diagram $D$, i.e.~$\mbox{sw}(D)$ is equal to  
to the sum of signs of the self-crossings of $D$.
Then the polynomial
$\hat{f}_D(A) = (-A^3)^{-\mbox{sw}(D)} \langle D\rangle$ 
is an invariant of global isotopy of the unoriented link 
determined by the diagram $D$.

\item The polynomial
$f_D(A) = (-A^3)^{-Tait(D)}\langle A\rangle =
(-A^3)^{-2\mbox{lk}(D)}\hat{f}_D(A)$ 
is invariant of global isotopy of oriented link determined
by the oriented diagram $D$.
\end{enumerate}
\end{lemma}

Proof. Both, $Tait(D)$ and $\mbox{sw}(D)$, are invariants 
of regular isotopy  (note that for $Tait(D)$ the diagram
$D$ must be oriented) and therefore $\hat{f}_D(A)$ and $f_D(A)$ 
are invariants of regular isotopy. Now, it is sufficient to prove
that both are preserved by the first Reidemeister move.
Since
$Tait({\psfig{figure=R+maly.eps}}) =
Tait({\psfig{figure=ver.eps,height=0.4cm}})
+1$ and
$Tait({\psfig{figure=R-maly.eps}})
=
Tait({\psfig{figure=ver.eps,height=0.4cm}})) -1$ then from V.3.4
%\ref{IV:2.3
%} 
and by the definition of $f_D$ and
$\hat{f}_D$ it follows that
$$f_{{\psfig{figure=R+maly.eps}} } (A) =
f_{{\psfig{figure=ver.eps,height=0.4cm}}} (A)
=
f_{{\psfig{figure=R-maly.eps}}} (A)$$
and
$$\hat{f}_{{\psfig{figure=R+maly.eps}}}
(A)
=
\hat{f}_{{\psfig{figure=ver.eps,height=0.4cm}}}
(A)
=
\hat{f}_{{\psfig{figure=R-maly.eps}}}
(A)$$
which completes the proof of Lemma 3.5.

\begin{theorem}\label{4:2.5}
The polynomial $f_L(A)$ is equal to the Jones polynomial $V_L(t)$ for
$A = t^{-\frac{1}{4}}$, that is
$$V_L(t) = f_L(t^{-\frac{1}{4}})$$
\end{theorem}

Proof. 
$$\langle
\parbox{0.6cm}{\psfig{figure=L+nmaly.eps}}
\rangle = A\langle
{\psfig{figure=L0nmaly.eps}}
\rangle
+
A^{-1}\langle
{\psfig{figure=Linftynmaly.eps}}
\rangle$$
and
$$\langle
\parbox{0.6cm}{\psfig{figure=L-nmaly.eps}}
\rangle = A^{-1}\langle
{\psfig{figure=L0nmaly.eps}}
\rangle+A\langle
{\psfig{figure=Linftynmaly.eps}}
\rangle$$
thus
$$A\langle
\parbox{0.6cm}{\psfig{figure=L+nmaly.eps}}
\rangle - A^{-1}\langle
\parbox{0.6cm}{\psfig{figure=L-nmaly.eps}}
\rangle = (A^2 - A^{-2})\langle
{\psfig{figure=L0nmaly.eps}}
\rangle$$
Assuming now that $L$ is oriented as in the diagram 
\parbox{0.6cm}{\psfig{figure=L+maly.eps}}, we will get:
$$A(-A^3)^{Tait({\psfig{figure=L+maly.eps}})}
%\Psfig{figure=sfo.eps} } + 1)}
f_{{\psfig{figure=L+maly.eps}}
%\Psfig{figure=skrzyzowanie2.eps}
} (A) -
A^{-1}(-A^3)^{Tait({\psfig{figure=L-maly.eps}})}
f_{{\psfig{figure=L-maly.eps}}} (A) =$$
$$(A^2 - A^{-2})(-A^3)^{Tait({\psfig{figure=L0maly.eps}})} 
f_{\psfig{figure=L0maly.eps}}(A)$$
thus 
$$A^4f_{{\psfig{figure=L+maly.eps}}}(A) -
A^{-4}f_{{\psfig{figure=L-maly.eps}}}(A) = (A^{-2} -
A^2)f_{{\psfig{figure=L0maly.eps}}}(A),$$
which for  $A = t^{-\frac{1}{4}}$ yields
$$ t^{-1}f_{{\psfig{figure=L+maly.eps}}}( t^{-\frac{1}{4}}) -
tf_{{\psfig{figure=L-maly.eps}}}( t^{-\frac{1}{4}}) =
(t^{\frac{1}{2}} - t^{-\frac{1}{2}})f_{{\psfig{figure=L0maly.eps}}}(
t^{-\frac{1}{4}})$$ 
The last formula is equivalent to the standard
skein relation satisfied by the Jones polynomial
$$ t^{-1}V_{{\psfig{figure=L+maly.eps}}}(t) -
tV_{{\psfig{figure=L-maly.eps}}}( t) =
(t^{\frac{1}{2}} - t^{-\frac{1}{2}})
V_{{\psfig{figure=L0maly.eps}}}(
t)$$
(for the other orientation of 
\parbox{0.6cm}{{\psfig{figure=L+nmaly.eps}}} we will
get a similar equation).

Moreover, for a trivial knot we get $f_\bigcirc(A) = 1 =V_\bigcirc (t)$,
which concludes the proof of Theorem \ref{4:2.5} (existence of the Jones
polynomial). Uniqueness is an easy exercise.

Theorem \ref{4:2.5} provides a short proof of
the Jones reversing result (Lemma III.5.15).

%(c.f.~Lemma \ref{l:5.15}).

\begin{corollary}\label{4:2.6}
Let us assume that $L_i$ is a component of an oriented
link $L$ and let us set $\lambda = \lk (L_i,L-L_i)$. 
Suppose that $L'$ is an oriented link obtained from $L$
by reversing the orientation of the component $L_i$.
Then
$V_{L'}(t) = t^{-3\lambda}V_L(t)$.
\end{corollary}

Proof. Let $D$ (resp. D') denotes a diagram of the link $L$ (resp. L'). 
The Kauffman bracket polynomial does not depend on the orientation of 
the diagram and therefore
$\langle D \rangle = \langle D'\rangle$. Subsequently
$$f_{D'}(A) = (-A^3)^{-Tait(D') + Tait(D)} f_D (A) = (-A^3)^{4\lambda} f_D(A)
= (A^4)^{3\lambda} f_D (A)$$
and thus $V_{L'} (t) = t^{-3\lambda} V_L (t)$.

The Kauffman interpretation of the Jones polynomial
provides also a short proof of a theorem of Lickorish,
Theorem III.5.13, that the Jones polynomial is a
specialization of the Kauffman polynomial.

\begin{theorem}\label{4:2.7}
%\begin{theorem}\label{V.3.8}
\ \\
\begin{enumerate}
\item[(1)]
 If $L$ is a diagram of an unoriented link
then  $\langle L\rangle = \Lambda_L(a,x)$ for $a= -A^3$, $x= A+A^{-1}$, 
that is
$$\langle L\rangle = \Lambda_L (-A^3, A+A^{-1}).$$

\item[(2)]
 If $L$ is an oriented link then
$$V_L(t) = F_L(-t^{-\frac{3}{4}},t^{\frac{1}{4}} + 
t^{-\frac{1}{4}}) =F_L(t^{-\frac{3}{4}},-(t^{\frac{1}{4}} + 
t^{-\frac{1}{4}})).$$
\end{enumerate}
\end{theorem}

Proof of (1). If $\bigcirc$ is a trivial  diagram of a knot then
$\langle \bigcirc \rangle = 1 = \Lambda_{\bigcirc}(a,x)$
and moreover $\langle 
{\psfig{figure=R+maly.eps}} 
\rangle = -A^3\langle {\psfig{figure=ver.eps,height=0.4cm}}
\rangle$ 
and $\Lambda_{
{\psfig{figure=R+maly.eps}}}
= -A^3\Lambda_{ {\psfig{figure=ver.eps,height=0.4cm}}
} (-A^3,A+A^{-1})$. Similarly 
$\langle {\psfig{figure=R-maly.eps}} \rangle = -A^{-3}\langle
{\Psfig{figure=).eps}}
{\psfig{figure=ver.eps,height=0.4cm}} \rangle$ and 
$\Lambda_{{\psfig{figure=R-maly.eps}} }(-A^3,A+A^{-1})
= -A^{-3}\Lambda_{ {\psfig{figure=ver.eps,height=0.4cm}}
} (-A^3,A+A^{-1})$.

Let us add the sides of the following two equations
$$\langle 
\parbox{0.6cm}{\psfig{figure=L+nmaly.eps}}\rangle = A\langle
{\psfig{figure=L0nmaly.eps}}\rangle + 
A^{-1}\langle {\psfig{figure=Linftynmaly.eps}}\rangle$$

$$\langle \parbox{0.6cm}{\psfig{figure=L-nmaly.eps}}\rangle = 
A^{-1}\langle
{\psfig{figure=L0nmaly.eps}}\rangle + 
A\langle {\psfig{figure=Linftynmaly.eps}}\rangle$$
obtaining 
$$\langle \parbox{0.6cm}{\psfig{figure=L+nmaly.eps}}\rangle + 
\langle
\parbox{0.6cm}{\psfig{figure=L-nmaly.eps}}\rangle = 
(A+A^{-1})(\langle
{\psfig{figure=L0nmaly.eps}}\rangle+\langle 
{\psfig{figure=Linftynmaly.eps}}\rangle),$$
which is equivalent to the equation\\
$\Lambda_{{\psfig{figure=L+nmaly.eps}}} (-A^3,A+A^{-1}) +
\Lambda_{{\psfig{figure=L-nmaly.eps}}} (-A^3,A+A^{-1}) =$\\
$(A+A^{-1})(\Lambda_{{\psfig{figure=L0nmaly.eps}}}(-A^3,A+A^{-1})
+\Lambda_{{\psfig{figure=Linftynmaly.eps}}}(-A^3,A+A^{-1}))$\\
for the Kauffman polynomial. This concludes the proof of (1).

Part (2) of the Theorem follows from (1) and Theorem
\ref{4:2.5}.

%\begin{definition}

Now let  $s_+$ (respectively, $s_-$) be a state of a diagram $D$ 
such that $s_+(i) = 1$ 
(respectively, $s_-(i) = -1$) for any crossing $c_i$ of $D$.
Then the  diagram $D$ will be called $+$-adequate if after changing it
to $s_+D$ by replacing 
\parbox{0.6cm}{\psfig{figure=L+nmaly.eps}} by 
\parbox{0.6cm}{\psfig{figure=L0nmaly.eps}} the newly
created arcs (for every crossing)
of the new diagram, $s_+D$, are in the different components of $s_+D$.
Similarly --- replacing $s_+D$ by $s_-D$ --- we define 
$-$-adequate diagram. A diagram is called adequate if it is $+$ and $-$
adequate (compare  Section 5; 
here we do not assume that every component of the 
diagram has a crossing). Equivalently, $+$ (resp. $-$) adequate
diagrams can be characterize as follows:
If a state $s$ differs from $s_+$ (resp. $s_-$) at one crossing only,
then for a $+$ adequate diagram (resp. $-$ adequate 
diagram) we have $|sD|<|s_+D|$ (resp. $|sD|<|s_-D|$). 
We will use this crucial
property of adequate diagrams in later considerations. 
%\end{definition}
%\begin{lemma}\label{V.3.9}
\begin{lemma}\label{4:2.8}
\begin{enumerate}
\item[(i)]
Reduced alternating diagrams 
%with no nugatory crossings
are adequate.
\item[(ii)] An alternating diagram is $+$-adequate if and only if 
all its nugatory crossings are positive (Fig. 3.2(a)).
\item[(iii)] An alternating diagram is $-$-adequate if and only if
all its nugatory crossings are negative (Fig. 3.2(b)).
\end{enumerate}
\end{lemma} 

\begin{proof} 
(i) For simplicity let us assume that $D$ is a connected 
alternating diagram. We color the components (regions) of the complement
of $D$ in the plane in black and white (checkerboard coloring).
>From the fact that the diagram is alternating it follows that
either all crossing look like 
%{\psfig{figure=L+nciemn.eps,height=0.4cm}}?
{\psfig{figure=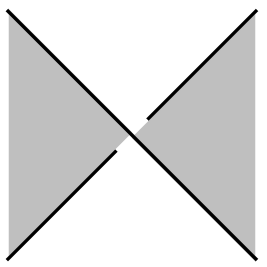,height=0.4cm}}
(and then the associated graph
$G(D)$ has only black edges) or all crossings look like 
%{\psfig{figure=L-ncien.eps,height=0.4cm}}?
{\psfig{figure=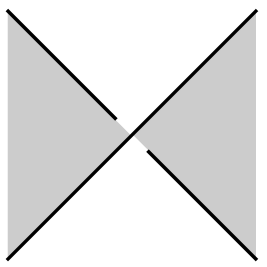,height=0.4cm}}
(and $G(D)$ has only white edges).
Now if $D$ has no nugatory crossing 
then no region is joined with itself by a crossing.
Therefore $D$ is an adequate diagram. We note that if all crossings
are of the type {\psfig{figure=L-nciemn.eps,height=0.4cm}} 
%{\psfig{figure=L+ncien.eps,height=0.4cm}}?
then 
$|s_+D|$ is equal to the number of black regions of the divided plane,
and $|s_-D|$ is equal to the number of white regions.
Therefore $|s_+ D| + |s_-D| = n(D)+2$ (we use an easy Euler 
characteristic argument).\\
(i)-(ii) It is explained in Fig. 3.2. For a nugatory crossing 
the positive marker agrees with the with the (orientation preserving) 
smoothing of the crossing. Recall that for a selfcrossing  the 
orientation preserving smoothing does not depend on an orientation of 
the link diagram.
\end{proof}
\centerline{{\psfig{figure=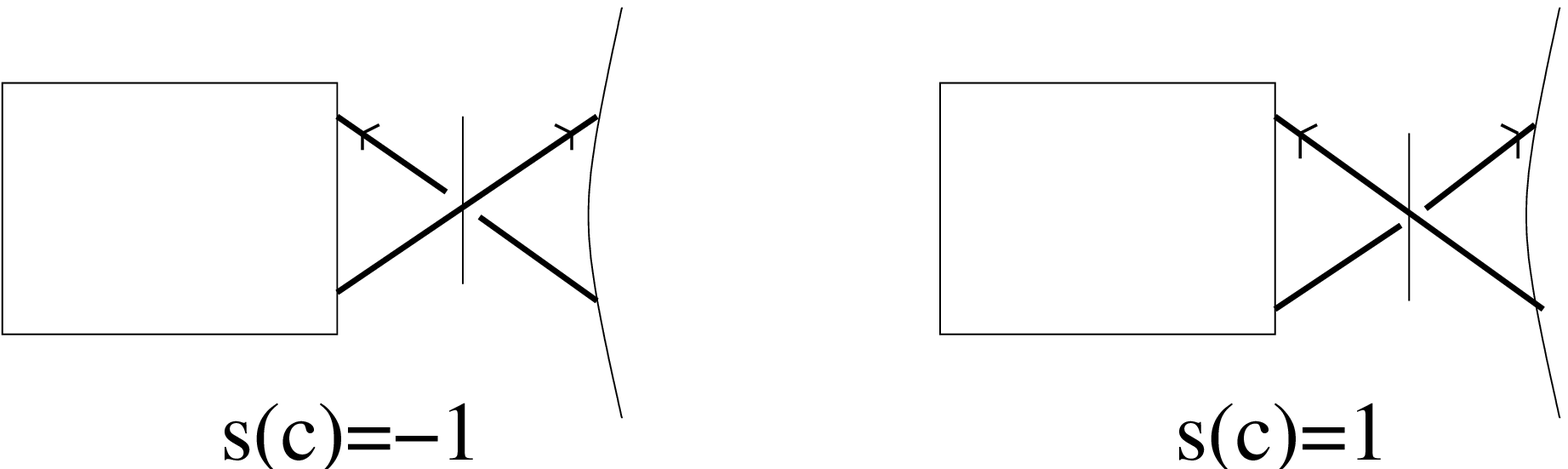,height=2.1cm}}}
\begin{center}
Fig. 3.2
\end{center}

We will show in Lemma 6.17 that every positive diagram is $+$-adequate.

Now let $\max\langle D\rangle$ and $\min\langle D\rangle$ denote maximal
and, respectively, minimal exponent of the variable 
$A$ in $\langle D\rangle$ and let
$\mbox{span }(D)$ be equal to
$\mbox{max }\langle D\rangle - \mbox{min }\langle D\rangle$. 

%\begin{lemma}\label{V.3.10}
\begin{lemma}\label{4:2.9}
Let $D$ be an unoriented $+$-adequate (resp. $-$-adequate) diagram 
of a link. Then the monomials
of $\langle D\rangle$ of maximal (respectively, minimal)
degree in $A$ are of the form
$$(-1)^{|s_+ D| - 1} A^{n+2|s_+ D|-2}$$ 
and, respectively,
$$(-1)^{|s_- D|-1} A^{-n-2|s_- D|+2}.$$
In particular $\max\langle D\rangle = n+2|s_+ D| - 2$, respectively 
$\min\langle D\rangle = -n-2|s_- D|+2$ and  if $D$ is adequate diagram then 
$\mbox{span }\langle D\rangle = 2n+2(|s_+
D| + |s_- D|) -4$. 
\end{lemma}

Proof. From the identity
$$A^{\sum s_+(i)}\langle s_+ D\rangle = A^n(-A^{-2}-A^2)^{|s_+
D| -1}$$ 
it follows that 
$$\max A^{\sum s_+(i)}\langle s_+ D\rangle  = n+2|s_+ D| - 2.$$ 
If now $s$ is any other state then there exists a sequence
of states $s_+ =\row{s}{k} = s$ such that any two subsequent states
$s_{r-1}$ and $s_r$ agree on all elements $i$ different than $i_r$ 
from the set
$\{ i:1\leq i\leq n\}$ and $s_{r-1}(i_r) = 1, s_r(i_r) = -1$. 
Thus, for $r\leq k$ we have $\sum_{i} s_r(i) = n-2r$ 
and $|s_r D| = |s_{r-1} D|\mp 1$.
It follows now that if $r$ increases to $r+1$ then $\max A^{\sum
s_r(i)}\langle s_r D\rangle$ decreases 
by 4 or is unchanged. Furthermore from the fact that $D$ is $+$ adequate
follows that $|s_1 D| = |s_{2} D| +1$ so in the first step ($s_1$ changed to
$s_2$)
$\max A^{\sum
s_1(i)}\langle s_1 D\rangle >\max A^{\sum
s_2(i)}\langle s_2 D\rangle$ therefore $\max A^{\sum
s_+(i)}\langle s_+ D\rangle  > A^{\sum
s_r(i)}\langle s_r D\rangle$ for all $r\geq 2$. 
Now from \ref{4:2.2} it follows that $\max\langle D\rangle = \max A^{\sum
s_+(i)}\langle s_+ D\rangle$ and thus we have the first part of \ref{4:2.9}.
The second part of the lemma, about the monomial of minimal degree in
$A$, can be proved similarly.

%\begin{corollary}\label{4:2.10}
\begin{corollary}\label{V.3.11}
\begin{enumerate}
\item[(i)] If the coefficient of the monomial of the maximal (resp. minimal) 
degree in $V_L(t)$ is not equal to $\pm 1$ then $L$ 
cannot be represented by a $+$ (resp. $-$) 
adequate diagram.
\item[(ii)]
Let $L$ be a connected alternating diagram of a link. If $L$ has
$n(L)$ crossings and none of them is nugatory then
$\mbox{span }\langle L\rangle = 4\mbox{span }V_L(t) = 4n$. 
\end{enumerate}
\end{corollary}

Proof. (i) It follows straight from Lemma 3.10.\\
(ii) We apply Lemma 3.10 and the fact that connected alternating 
diagram with no nugatory crossing is adequate and 
$|s_+ L| + |s_- L| = n(L)+2$. Clearly $4\mbox{span }V_L(t) =
\mbox{span }\langle L\rangle$ (Theorem \ref{4:2.5}).

%\begin{corollary}\label{V.3.12}
\begin{corollary}\label{4:2.11}
If $L$ is a prime non-split link  
%(i.e. $L$ is not a split or connected sum of links)
then for any nonalternating diagram
of $L$ we have $\spn V_L(t) < n(L)$. 
\end{corollary}

Corollaries 2.11 and 2.12 
%\ref{4:2.10} and \ref{4:2.11} 
imply the first Tait conjecture.

Corollary \ref{4:2.11} follows from the subsequent observation,
which has an easy proof when translated to the fact on 2-color graphs.

\begin{lemma}[On dual states \cite{K-6,Wu}.]\label{4:2.12}
\ \\
\begin{enumerate}
\item[(1)]
Let $D$ be a connected diagram of a link and let $s$ be its state.
let $s^\star$ denote the dual state $s^\star (i) = -s(i)$. 
Then $|s^\star D| + |sD|\leq n(D)+2$. 
\item[(2)] If $D$ is not a connected sum of connected alternating diagrams
then  $|s_+D | +|s_-D|<n(L) +2$. 
%More generally, the above inequality is true if $L$ is not a connected
%sum of alternating links.
\end{enumerate}
\end{lemma}
Hint.\ \  A simple way to show Lemma 3.13 is to translate it
to the language of graphs (as shown in Fig.~1.4). The respective
lemma for graphs is true also for non-planar graphs (hence more general
than we need to prove Lemma 3.13).

%\begin{lemma}\label{4:2.12a}\ \\
\begin{lemma}\label{V.3.14}\ \\
Let $G$ be a connected 2-color graph with edges colored in $b$ and $w$.
Let $G_b$ (respectively, $G_w$) be a graph which has the vertices 
of the graph $G$ and edges of color $b$ (respectively, $w$) --- taken 
from $G$.  Then
$$p_0(G_b)+p_1(G_b)+ p_0(G_w)+p_1(G_w)\leq E(G)+2,$$
and the inequality becomes equality if and only if
$G$ has no 2-color cycle (i.e every cycle is compose of only $b$ edges or 
only $w$ edges). 
\end{lemma}
\begin{proof}
Since for any graph $H$ we have $p_1(H)=E(H)-V(H)+p_0(H)$ 
then the inequality from  Lemma 3.14 can be reduced to
$2(p_0(G_b)+p_0(G_w))\leq 2(V(G)+1)$. We leave the rest of the proof
in the form of the following easy exercise.
\end{proof}
\begin{exercise}\label{4:2.12b}
Let $G$ be a 2-color graph as in Lemma 3.14,
but not necessarily connected. Then $p_0(G_b)+p_0(G_w)\leq V(G)+p_0(G)$
and the inequality becomes equality if and only if $G$ does not 
contain 2-color cycles.
\end{exercise}
Hint.\ Solve the exercise first for $G$ being a forest.

In the fifth section we will show (following Thistlethwaite) how to use
Kauffman polynomial to prove second Tait conjecture.
Murasugi \cite{M-4,M-5} proved this conjecture
applying Jones polynomial and signature.
Namely, he showed the following

%\begin{theorem}[\cite{M-5}.]\label{4:2.13}
\begin{theorem}[\cite{M-5}.]\label{V.3.16}
For any connected diagram of a link, $L$, we have:
\begin{enumerate}

\item $\max V_L (t) \leq n^+(L) - \frac{1}{2}\sigma (L)$

\item  $\min V_L (t) \geq -n^-(L) - \frac{1}{2}\sigma (L)$
\end{enumerate}
where $n^+(L)$ (respectively, $n^-(L)$) denotes the number of positive 
(respectively, negative) crossings of the diagram $L$.

\end{theorem}

Both inequalities become simultaneously equalities if and only if 
$L$ is an alternating diagram without nugatory crossings or it is 
a connected sum  of such diagrams. We demonstrate, after Traczyk, 
the equality for alternating diagrams in Chapter IV.

\begin{corollary}\label{4:2.14}
If $L$ is either an alternating diagram with no nugatory crossing
or a connected sum  of such diagrams then the Tait or writhe number 
$Tait(L)$ (denoted also by $\tilde n(L)$) 
is equal to
$$n^+(L) - n^-(L) = \max V_L (t) + \min V_L(t) + \sigma (L).$$
\end{corollary}
We can obtain additional properties of Kauffman bracket and Jones 
polynomial of alternating link diagrams from corresponding properties 
of Tutte polynomial (Corollary 1.9 and Exercise 1.11), in particular 
we prove that alternating links have alternating 
Jones polynomial \cite{This-3}.

\begin{theorem}[Thistlethwaite]\label{V.3.18} 
\begin{enumerate}
\item[(a)]
If $L$ is a nonsplit alternating link, then the coefficients 
of the Jones polynomial of $L$ are alternating (we allow $0$).
\item[(b)] 
If $L$ is a nonsplit prime alternating link different from a 
$(2,k)$ torus link then every coefficient of $V_L(t)$ between 
$\max V_L(t)$ and $\min V_L(t)$ is different from zero. For 
a positive $(2,k)$ torus link, $T_{2,k}$, we have 
$V_{T_{2,k}}(t)= 
-t^{\frac{1}{2}(n-1)}(t^{n}- t^{n-1} +...+(-1)^{n-2}t^2 + (-1)^n)$
\item[(c)] If $L$ is a reduced diagram of a
 nonsplit prime alternating link different from a $(2,k)$ torus link then 
\begin{enumerate}
\item[(i)] the Kauffman bracket polynomial satisfies
$$<L> =\sum_{i=0}^n(-1)^{i+|s_-L|-1}a_iA^{4i-n-2|s_-L|+2}
\ \textup{with}\ a_i>0,\ a_0=a_n=1.$$ 
\item[(ii)] the Jones polynomial of an  oriented $\vec L$ satisfies
$$V_{\vec L}(t)= \sum_{i=0}^n (-1)^{(i-|s_-L|+1)}
a_{n-i}t^{i-n^- - \frac{1}{2}\sigma (\vec L)},$$ 
where $\vec L$ is 
an oriented link diagram with underlining unoriented link diagram $L$ 
\end{enumerate}
\end{enumerate}
\end{theorem}
\begin{proof} We apply Theorem 1.6 for $B=A^{-1}$ and $\mu =-A^2-A^{-2}$.
Then $x=-A^{-3}$ and $y=-A^3$, and the Kauffman bracket polynomial of 
a nonsplit alternating link obtained from a (black edged) graph is 
up to $\pm A^i$ obtained by this substitution from the Tutte polynomial.
Now Theorem 3.18 follows from Corollary 1.9, Exercise 1.11 and Theorem 3.6.
Part (c) is a combination of (a),(b), Lemma 3.10 and Theorem 3.16.
\end{proof}

We will finish this section by describing one more possible
generalization of the Kauffman bracket polynomial of a diagram of a link $<D>$,
leading to a 2-variable polynomial invariant of links.
We will show however that we do not gain any new information in
this approach. Following our definition for polynomial of graphs,
we can define the Kauffman bracket of link diagrams as a polynomial 
in three variables
$<D>_{A,B,\mu}\in Z[A,B,\mu]$ which satisfies the following conditions:
\begin {enumerate}
\item
[(a)] $<T_n>_{A,B,\mu}={\mu}^{n-1},$
\item
[(b)] $<
\parbox{0.6cm}{\psfig{figure=L+nmaly.eps}}
>_{A,B,\mu} = A<
\parbox{0.6cm}{\psfig{figure=L0nmaly.eps}}
 >_{A,B,\mu}  + B<
\parbox{0.6cm}{\psfig{figure=Linftynmaly.eps}}
 >_{A,B,\mu}.$
\end{enumerate}
A direct induction with respect to the number of crossings of $D$, 
denoted $n(D)$,
provides that $<D\sqcup  \bigcirc>_{A,B,\mu}=\mu<D>_{A,B,\mu}$.
Considering the second Reidemeister move we get:
$<
\parbox{0.6cm}{\psfig{figure=R2-k.eps,height=0.5cm}}
%{\psfig{figure=R2maly.eps}}
>_{A,B,\mu} =AB<
\parbox{0.6cm}{\psfig{figure=Linftynmaly.eps}}
 >_{A,B,\mu} +
(A^2+B^2+\mu AB)<
\parbox{0.6cm}{\psfig{figure=L0nmaly.eps}}
 >_{A,B,\mu}.$  
If we assume that 
$A^2+B^2+\mu AB=0$ then $\mu=-\frac{A^2+B^2}{AB}=-\frac{A}{B}-
\frac{B}{A}$ and $<
\parbox{0.6cm}{\psfig{figure=R2-k.eps,height=0.5cm}}
%{\psfig{figure=R2maly.eps}}
 >_{A,B,\mu} =AB<
\parbox{0.6cm}{\psfig{figure=Linftynmaly.eps}}
>_{A,B,\mu}.$
%\begin{exercise}\label{V.3.19}
\begin{exercise}\label{IV:2.17}
 Prove that, for  $\mu =-(\frac{A}{B}+\frac{B}{A}),$
the (generalized) Kauffman bracket polynomial
$<D>_{A,B}\in Z[A^{\pm 1},B^{\pm 1}]$ has the following properties:
\begin{enumerate}
\item
[(a)] $<D>_{A,B}$ is preserved by the third
Reidemeister move.
\item
[(b)] $(AB)^{-n(D)/2}<D>_{A,B}$ is an invariant of regular isotopy.
\item
[(c)]  $<
\parbox{0.6cm}{\psfig{figure=R+maly.eps}} >_{A,B} =-A^2B^{-1}<
\parbox{0.1cm}{\psfig{figure=ver.eps,height=0.4cm}}
 >_{A,B}$, $<
\parbox{0.6cm}{\psfig{figure=R-maly.eps}}
 >_{A,B}= -B^2A^{-1} <
\parbox{0.2cm}{\psfig{figure=ver.eps,height=0.4cm}} 
>_{A,B}$,
\item[(c')]  $(AB)^{-n(\parbox{0.5cm}{\psfig{figure=R+maly.eps}})/2}
<\parbox{0.5cm}{\psfig{figure=R+maly.eps}} >_{A,B} = -A^{3/2}B^{-3/2}
<\parbox{0.1cm}{\psfig{figure=ver.eps,height=0.4cm}}> = -(AB^{-1})^{3/2}
<\parbox{0.1cm}{\psfig{figure=ver.eps,height=0.4cm}}>,$ 
\item
[(d)]
${\hat f}_D(A,B)= (AB)^{-n(D)/2}(-(AB^{-1})^{3/2})^{-sw(D)}
<D>_{A,B}$ is an invariant of unoriented links.
\item
[(e)] $f_{\vec{D}}(A,B)=(AB)^{-n(D)/2}
(-(AB^{-1})^{3/2})^{-Tait(\vec{D})}<D>_{A,B},$
where $\vec{D}$ is an oriented diagram obtained from $D$ by
equipping it with an orientation, is an invariant of oriented links.
\end {enumerate}
\end{exercise}

Consequently, one may think that $ <D>_{A,B}$ provides a better
knot invariant than the usual Kauffman bracket of one variable.
This is, however, not the case, as we see from the following exercise
%\begin{exercise}\label{V.3.20}
\begin{exercise}\label{IV:2.18}
Prove that, if $<D>=\Sigma a_iA^i$ then
$$<D>_{A,B}=\Sigma a_iA^{(n(D)+i)/2}B^{(n(D)-i)/2}$$
$$ and \ \ (AB)^{-n(D)/2}<D>_{A,B}=\Sigma a_i(AB^{-1})^{i/2}.$$
\end{exercise}

%%%%%%%%%%%%%%%%%%%%%%%%%%%%%%%%%%%%%%%%%%%%%%%%%%%%%%%%%%%%%%%%

\section{Application of Kauffman polynomial 
to alternating links}\label{V.4}
\markboth{\hfil{\sc Graphs and links}\hfil} 
{\hfil{\sc Kauffman polynomial of alternating links}\hfil}

We start with a special version of the Kauffman polynomial, which is
the polynomial of Brandt-Lickorish-Millett  and Ho \cite{B-L-M, Ho}, 
compare Chapter II for historical remarks.
We denote this polynomial by $Q_L(x)$ and we get it by setting
$a=1$ in the Kauffman polynomial. Therefore, $Q_L(x)$ is uniquely
defined by the conditions:

\begin{formulla}\label{4:3.1}
$$
\left\{
\begin{array}{lc}
(i)&Q_{T_1}(x) = 1 \mbox{ for the trivial knot } T_1\\
(ii)&Q_{L_{{\psfig{figure=L+nmaly.eps,height=0.4cm}}}}(x) +
Q_{L_{{\psfig{figure=L-nmaly.eps,height=0.4cm}}}}(x) =
xQ_{L_{{\psfig{figure=L0nmaly.eps,height=0.4cm}}}}(x) + 
xQ_{L_{{\psfig{figure=Linftynmaly.eps,height=0.4cm}}}}(x).\\
\end{array}
\right.
$$
\end{formulla}

As a consequence, for the trivial link of $n$ components, $T_n$, 
we have $Q_{T_n}(x) = (\frac{2-x}{x})^{n-1}$, or more generally 
$Q_{L\sqcup O}(x)= (\frac{2-x}{x})Q_{L}$.

The following two theorems are from a paper of M.~Kidwell, \cite{Kid}.

\begin{theorem}\label{4:3.2}
%\begin{theorem}\label{V.4.2}
Let $L$ be a diagram of a link with $n(L)$ crossings.
Let $b(L)$ denote the length of the longest bridge in $L$. 
%(???czy to zostalo zdefiniowane???czy dobrze przetlumaczone???)
Then
$\deg Q_L\leq n(L) - b(L)$, where by $\deg Q_L = max\ deg\ Q_L$ 
we understand the highest
degree of $x$ in $Q_L(x)$ (the polynomial $Q$ can have terms of negative
degree as well\footnote{If $L$ has $com(L)$ components then 
$min\ deg\ Q_L(x) = 1-com(L)$, as can be checked by induction.}).
\end{theorem}

\begin{theorem}\label{4:3.3}
Let $L$ be a connected prime alternating diagram of a link with
$n(L)>0$ crossings. Then the coefficient of the monomial $x^{n(L)-1}$ in $Q(L)$ 
is positive. A connected diagram of a link is called prime 
if there is no a simple closed curve $C$ on the plane which 
meets $L$ transversally in two points
and each of the two components of the complement of $C$ contains 
a crossing of $L$. 
Figure 4.1 presents examples of prime diagrams of links.
\\
%\vspace*{1.5in}\centerline{\Psfig{figure=Rys.3.1}}
\ \\
\centerline{{\psfig{figure=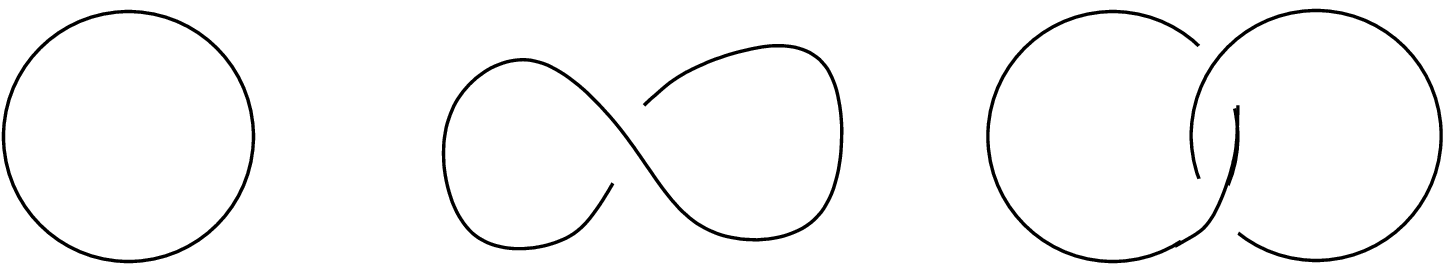,height=2.1cm}}}
\begin{center}
Fig.~4.1
\end{center}

\end{theorem}

Proof of Theorem \ref{4:3.2}.

Assume that the there exists a counterexample
to our theorem. Let $L$ be the counterexample with the smallest number
of crossings and the longest bridge (among these diagrams which have 
the smallest number of crossings). 
Therefore  $\deg Q_L > n(L) - b(L)$. 
Let $B$ be a bridge of length $b(L)$ in $L$. Now we have two
possibilities: either (1) the bridge is not proper or (2) it is proper.
Let us explain both.
\begin{enumerate}
\item The bridge is not proper which means that either

\begin{enumerate}

\item $B$ is a simple closed curve and $b(L)>0$ (Fig.~4.2(i)), or

\item $B$ ends, at least at one side, with a tunnel passing under itself
(Fig.~4.2(ii)), or

\item $B$ ends from both sides with the same tunnel and $b(L)>1$
(Fig.~4.2(iii)).
\end{enumerate}
%\vspace*{2in}\centerline{\Psfig{figure=Rys.3.2}}
\ \\
\centerline{{\psfig{figure=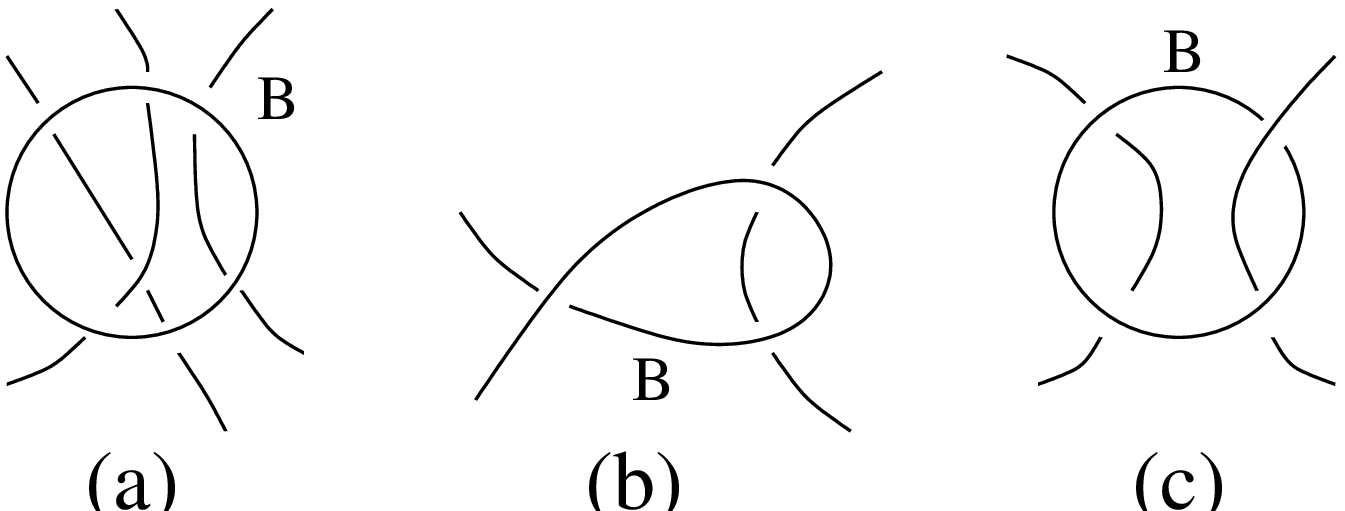,height=3.1cm}}}
\begin{center}
Fig.~4.2
\end{center}

In cases 1(a) and 1(b) we can change $L$ by using isotopy to get a diagram
$L'$ with $n(L') = n(L) - b(L)$ crossings. In the case 1(c) we can reach
(via isotopy) a diagram $L'$ which has at least one crossing 
(i.e.~$b(L')\geq 1$) and such that $n(L') = n(L) - b(L) +1$. 
In any case we get $n(L') - b(L')\leq n(L) -b(L)$. 
Moreover $n(L')< n(L)$ which contradicts our assumption on $L$ 
(note that $\deg Q_L(x) = \deg Q_{L'}(x)$ since $L$ and $L'$ are isotopic).

\item The bridge $B$ is proper if, by definition,  
none of the above situations 
(neither (a), (b) nor (c)) is true.
For diagrams with no crossing $\deg Q_L(x) = 0$, and 
\ref{4:3.2} is true. Therefore in our counterexample 
$L$  we have $n(L)\geq 1$ and $b(L)\geq 1$. 
Since $B$ is a proper bridge, the crossing which ends it --- call it $p$ ---
is not a part of $B$  (see Fig.~4.3)
\\
%\vspace*{1.5in}\centerline{\Psfig{figure=Rys.3.3}}
\ \\
\centerline{{\psfig{figure=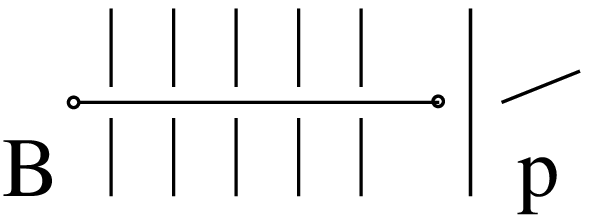,height=2.1cm}}}
\begin{center}
Fig.~4.3
\end{center}

Now let us consider the crossing $p$. The diagrams 
$L_{{\psfig{figure=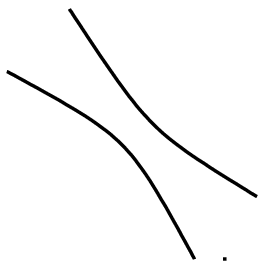,height=0.4cm}}}$  and
$L_{{\psfig{figure=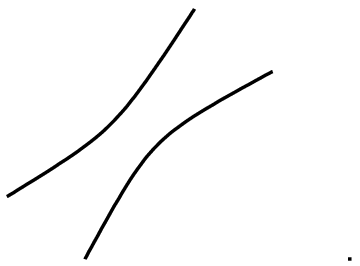,height=0.4cm}}}$ have 
one crossing less than $L$ but 
their longest bridges are not shorter than $B$. 
Because of our assumption on $L$,  the diagrams 
$L_{{\psfig{figure=nw0.eps,height=0.4cm}}}$  and
$L_{{\psfig{figure=ne0.eps,height=0.4cm}} }$ satisfy 
the following inequalities
$$\deg Q_
{L_{{\psfig{figure=nw0.eps,height=0.4cm}} }} \leq 
n(L_{{\psfig{figure=nw0.eps,height=0.4cm}} }) - 
b(L_{{\psfig{figure=nw0.eps,height=0.4cm}} })$$ and
$$\deg Q_{L_{{\psfig{figure=ne0.eps,height=0.4cm}} }}\leq
n(L_{{\psfig{figure=ne0.eps,height=0.4cm}} }) - 
b(L_{{\psfig{figure=ne0.eps,height=0.4cm}} })$$
and therefore
$$\deg(xQ_{L_{{\psfig{figure=nw0.eps,height=0.4cm}}}} 
+xQ_{L_{{\psfig{figure=ne0.eps,height=0.4cm}} }})
\leq\max((n(L_{{\psfig{figure=nw0.eps,height=0.4cm}} } ) - 
b(L_{{\psfig{figure=nw0.eps,height=0.4cm}} }),$$ 
$$(n(L_{{\psfig{figure=ne0.eps,height=0.4cm}}}) - 
b(L_{{\psfig{figure=ne0.eps,height=0.4cm}} } ))+1\leq n(L) - b(L).$$

Moreover $n(L_{{\psfig{figure=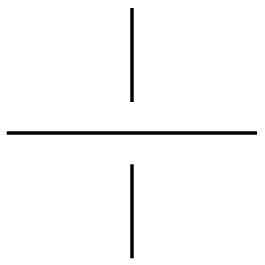,height=0.4cm}}}) = n(L)$ and
$b(L_{{\psfig{figure=h-v-cr.eps,height=0.4cm}}}) > b(L)$ hence, because
of the assumption that $L$ is minimal among counterexamples, 
it follows that the theorem is true for 
$L_{{\psfig{figure=h-v-cr.eps,height=0.4cm}}}$  and thus
$$\deg Q_{L_{{\psfig{figure=h-v-cr.eps,height=0.4cm}}}}(x)\leq
n(L_{{\psfig{figure=h-v-cr.eps,height=0.4cm}}}) - 
b(L_{{\psfig{figure=h-v-cr.eps,height=0.4cm}}})
< n(L) -b(L).$$ 
Hence
$$\deg Q_{L} = \deg (-Q_{L_{{\psfig{figure=h-v-cr.eps,height=0.4cm}}}}
+x(Q_{L_{{\psfig{figure=nw0.eps,height=0.4cm}}}} 
+ Q_{L_{{\psfig{figure=ne0.eps,height=0.4cm}}}})) \leq n(L) -b(L)$$ 
and $L$ can not be a counterexample to our theorem, which concludes
the proof of \ref{4:3.2}
\end{enumerate}

Proof of Theorem \ref{4:3.3}

The theorem is true for a diagram with one crossing:
{\psfig{figure= 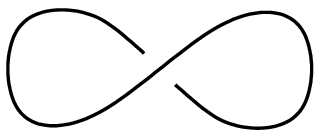,height=0.4cm}} and
{\psfig{figure=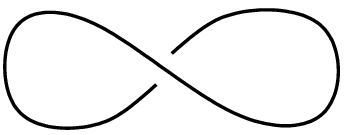,height=0.4cm}}.

Let us assume that the theorem is true for diagrams with less than $n(L)$
crossings ($n(L)\geq 2$).  
Now let us consider an arbitrary crossing $p$ of 
the diagram $L = L_{{\psfig{figure=L+nmaly.eps,height=0.4cm}}}$.
The diagram $ L_{{\psfig{figure=L-nmaly.eps,height=0.4cm}}}$ 
has a bridge of length at least two. Hence
$$\deg Q_{ L_{{\psfig{figure=L-nmaly.eps,height=0.4cm}}}} 
\leq n( L_{{\psfig{figure=L-nmaly.eps,height=0.4cm}}}) -2 
= n( L_{{\psfig{figure=L+nmaly.eps,height=0.4cm}}}) - 2
= n(L) -2.$$

Therefore the coefficient of $x^{n(L) -1}$ in $Q_L$ is equal to
the coefficient of $x^{n(L) -2}$ 
in $Q_{L_{{\psfig{figure=L0nmaly.eps,height=0.4cm}}}} 
+ Q_{L_{{\psfig{figure=Linftynmaly.eps,height=0.4cm}}}}$. 
%(?I tutaj sie zgubilem?co to jest \L ? czy bylo zdefiniowane ?) 
The diagrams $L_{{\psfig{figure=L0nmaly.eps,height=0.4cm}} }$ 
and $L_{{\psfig{figure=Linftynmaly.eps,height=0.4cm}}}$ are alternating
with $n(L) -1$ crossings. We will be done if we prove that 
either $L_{{\psfig{figure=L0nmaly.eps,height=0.4cm}}}$ 
or $L_{{\psfig{figure=Linftynmaly.eps,height=0.4cm}}}$ 
is a prime connected
diagram (by inductive assumption). Note also that if $D$ is not a prime
diagram than $\deg Q_D < n(D) -1$ ((c.f. Exercise 3.11). 
%(c.f.~exercise).  (???caly ten ostatni paragraf jest dla mnie niejasny???)

%\begin{lemma}\label{V.4.4}
\begin{lemma}\label{4:3.4}
If $L$ is a prime connected diagram and  $p$ is an arbitrary crossing
of $L$ then either $L_{{\psfig{figure=L0nmaly.eps,height=0.4cm}} }$ 
or $L_{{\psfig{figure=Linftynmaly.eps,height=0.4cm}} }$ 
is a prime and connected diagram.
\end{lemma}

With the exception of $L=${\psfig{figure= 1-cr+fig.eps,height=0.4cm}} or
{\psfig{figure=1-cr-.eps,height=0.4cm}} the diagrams
 $L_{{\psfig{figure=L0nmaly.eps,height=0.4cm}} }$ 
and $L_{{\psfig{figure=Linftynmaly.eps,height=0.4cm}} }$ 
are connected. 

We will prove that one of them is prime. In fact we show that
this statement follows from  the fact, proved in Lemma 1.8, 
that if a graph $G$ is 2-connected 
 then $G-e$ or $G/e$ is 2-connected for every edge $e$ of $G$.
%\begin{lemma}\label{IV.3.5}
\begin{lemma}\label{V.4.5}
A connected diagram, $D$, of a link is prime if and only if an
associated graph $G(D)$ is 2-connected. The lemma holds for 
any checkerboard coloring of regions of the diagram complement and 
any decoration of the graph.
\end{lemma}
\begin{proof}
If $D$ is composite then there is a closed curve cutting $D$ in two
points into $D_1\# D_2$. The same curve divides the graph $G(D)$
into $G(D_1)*G(D_2)$. The same reasoning also shows that if
$G(D)$ is not 2-connected than $D$ is not a prime connected diagram.
\end{proof}

%{Let us assume the contrary, that is, suppose that both  
%$L_{{\psfig{figure=L0nmaly.eps,height=0.4cm}} }$ 
%and $\L_{{\psfig{figure=Linftynmaly.eps,height=0.4cm}} }$ 
%are not prime.
%\centerline{{\psfig{figure=4-3-4.eps,height=3.2cm}}} 
%\begin{center} Fig.~3.4 \end{center}
%\centerline{{\psfig{figure=4-3-5.eps,height=3.2cm}}}
%\begin{center} Fig.~3.5 \end{center}
%$g$, $d$, $\pi$ and $e$.
Lemma 4.4 follows from Lemma 4.5 and Lemma 1.8.

\begin{corollary}\label{4:3.5}
Consider a connected prime alternating
diagram of a link $L$ with $n(L)\geq 3$ and which contains a clasp
as pictured in Fig.~4.4(a). 
Let $L^p_{{\psfig{figure=L0nmaly.eps,height=0.4cm}}}=
L_{{\psfig{figure=L0nmaly.eps,height=0.4cm}}}$ denote, as before, 
the diagram obtained from $L$ by smoothing the crossing $p$ horizontally 
(Fig. 4.4(d)). Similarly, the meaning of 
$L_{{\psfig{figure=Linftynmaly.eps,height=0.4cm}}}$ and 
$L_{{\psfig{figure=L+nmaly.eps,height=0.4cm}}}$ is illustrated 
in Fig. 4.4.
Then the polynomials $Q_{L_{{\psfig{figure=L0nmaly.eps,height=0.4cm}}}}$ 
and $Q_L$
have the same coefficient at the term of the highest degree in $x$.
The modification $L\rightarrow 
L_{{\psfig{figure=L0nmaly.eps,height=0.4cm}}}$ is called
elimination of a clasp.

%\vspace*{2in}\centerline{\Psfig{figure=Rys.3.5}}
\ \\
\centerline{{\psfig{figure=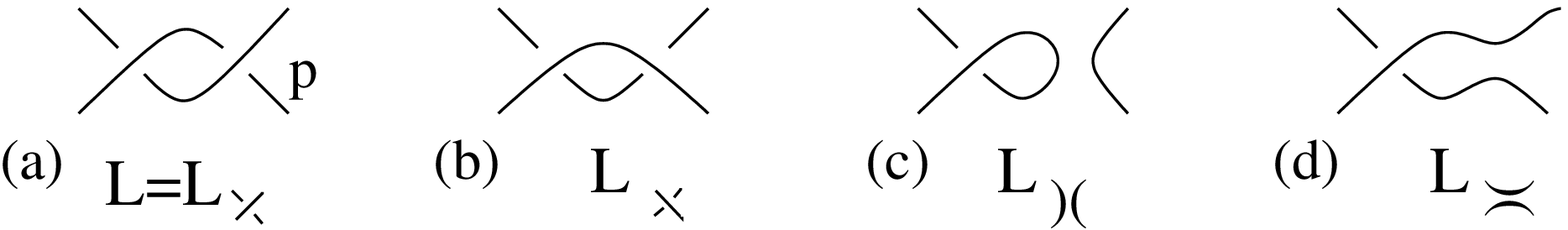,height=2.2cm}}}
\begin{center}
Fig.~4.4
\end{center}

\end{corollary}

Proof. 
Since the diagram  $L_{{\psfig{figure=Linftynmaly.eps,height=0.4cm}}}$ 
is not prime hence by  Lemma \ref{4:3.4}, it follows
that $L_{{\psfig{figure=L0nmaly.eps,height=0.4cm}}}$ 
is connected, prime and alternating.
Now Theorem 4:3 yields that the highest degree terms
of $Q_L$ and $Q_{L_{{\psfig{figure=L0nmaly.eps,height=0.4cm}}}}$ 
are $ax^{n(L)-1}$ and $bx^{n(L)-2}$, 
respectively. From Figure 4.4. it follows that the highest exponent of
$x$ in $Q_{L_{{\psfig{figure=Linftynmaly.eps,height=0.4cm}}}}$ and 
$Q_{L_{{\psfig{figure=L+nmaly.eps,height=0.4cm}}}}$ is at most
${n(L)-3}$ and thus, because of the recursive
 definition of $Q$, we get $a=b$. 
%(???troche zmienilem, bo w polskiej wersji jednomiany, wspolczynniki i

\begin{corollary}\label{4:3.6}
If a given connected prime alternating
diagram of a link $L$ can be reduced by an elimination of clasps
to a Hopf
diagram (\parbox{0.8cm}{\psfig{figure=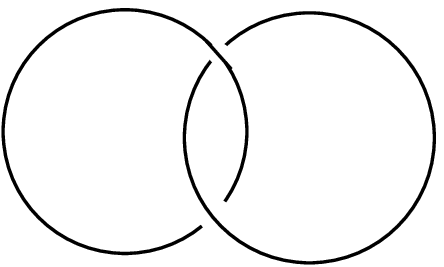,height=0.5cm}}) 
then the leading coefficient
(i.e.~the coefficient of the highest degree monomial) of $Q_L(x)$ is equal 
to 2.
\end{corollary}

Proof. We note that $Q_{{\psfig{figure=Hopfn.eps,height=0.4cm}}}(x) 
= -2x^{-1} +1 +2x$ 
and next we apply Corollary \ref{4:3.5}.

\begin{exercise}\label{4:3.7}
Let us consider the following family of tangles
which Conway called rational tangles, defined inductively as follows: 
\begin{enumerate}
\item[(1)]
\  {\psfig{figure=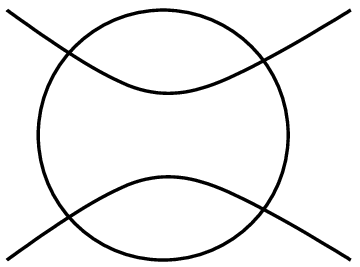,height=0.4cm}}
and {\psfig{figure=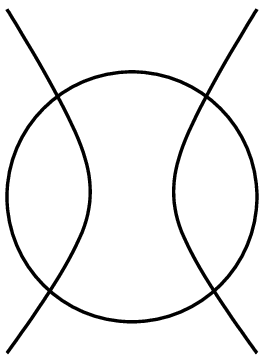,height=0.4cm}}
are rational tangles.

\item[(2)] 
If {\psfig{figure=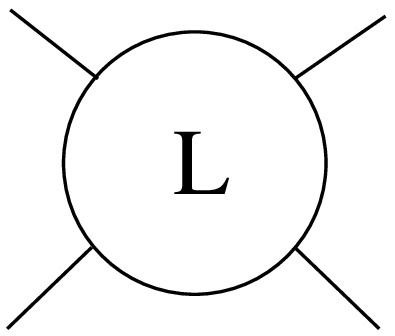,height=0.4cm}} 
is a rational tangle then
{\psfig{figure=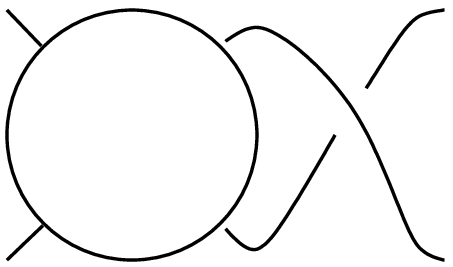,height=0.4cm}}\ \ , 
{\psfig{figure=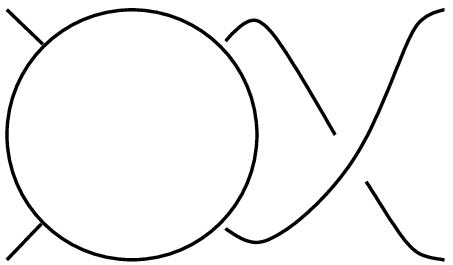,height=0.4cm}}\ \ ,
{\psfig{figure=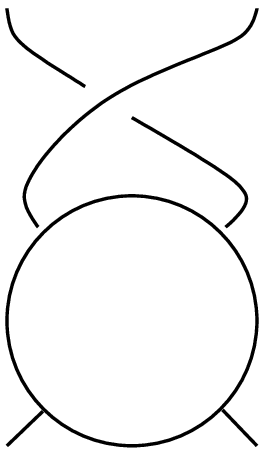,height=0.5cm}}
and {\psfig{figure=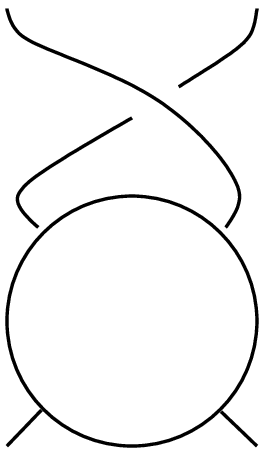,height=0.5cm}} are rational
tangles.
\end{enumerate}

The numerator {\psfig{figure=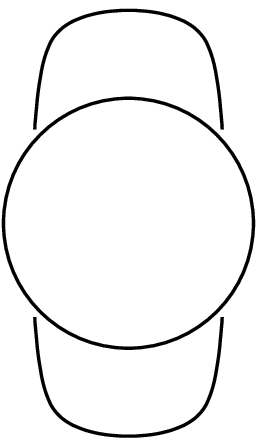,height=0.5cm}}
and the denominator 
{\psfig{figure=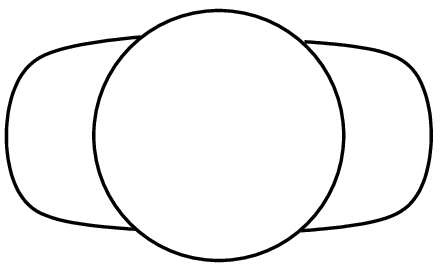,height=0.4cm}}
of a rational tangle are called rational links.
Show that the leading coefficient in the
polynomial $Q(x)$ of a rational link is equal 2.

Hint. Apply Corollary \ref{4:3.6} and use the fact that rational
tangles can be associated with rational numbers and diagrams of 
rational tangles links can be described by continuous 
fractions.
In particular, we can associate 
to any alternating diagram of a rational tangle
a positive fraction (c.f.~\cite{Co-1} and \cite{B-Z}; compare also 
Chapter XX).
\end{exercise}
The family of rational links was extended by Conway to the 
algebraic links (Chapter XX). Formulate the version of Exercise 4.8 
for prime links which have an alternating algebraic diagram.
\begin{exercise}[The first Tait Conjecture]\label{4:3.8}\ \\
Prove that all connected prime alternating
diagrams (of more than one crossing) 
of a given link have the same number of crossings
which is smaller than the number of crossings in any non-alternating
diagram of the link in question.

Hint. Apply theorems \ref{4:3.2} and \ref{4:3.3}.
\end{exercise}

A link is called alternating if it admits an alternating diagram.
\begin{exercise}\label{4:3.9}
Prove that if an alternating link is not prime (i.e.~it decomposes to
a connected sum) then any alternating diagram of it
is composed (this was first proved by  Menasco \cite{Men}).
We summarize this fact by saying that alternating diagram of 
a composite alternating link is visibly composite.
 
Hint. If $L=L_1\kwad L_2$ then $Q_L = Q_{L_1}\cdot Q_{L_2}$ 
and therefore $\deg Q_L < n(L)- 2$. On the other hand if the diagram 
was not composite (but alternating) than  
$\deg Q_L = n(L) - 1$.
\end{exercise}

\begin{exercise}\label{4:3.10}
Let us define a generalized bridge of a diagram $L$
to be the part of the diagram which is descending. That is,
moving along a generalized bridge, any crossing which we meet
for the first time is passed by overcrossing.
(c.f.~Fig.~4.5).

% \vspace*{1.5in}\centerline{\Psfig{figure=Rys.3.6}}
\ \\
\centerline{{\psfig{figure=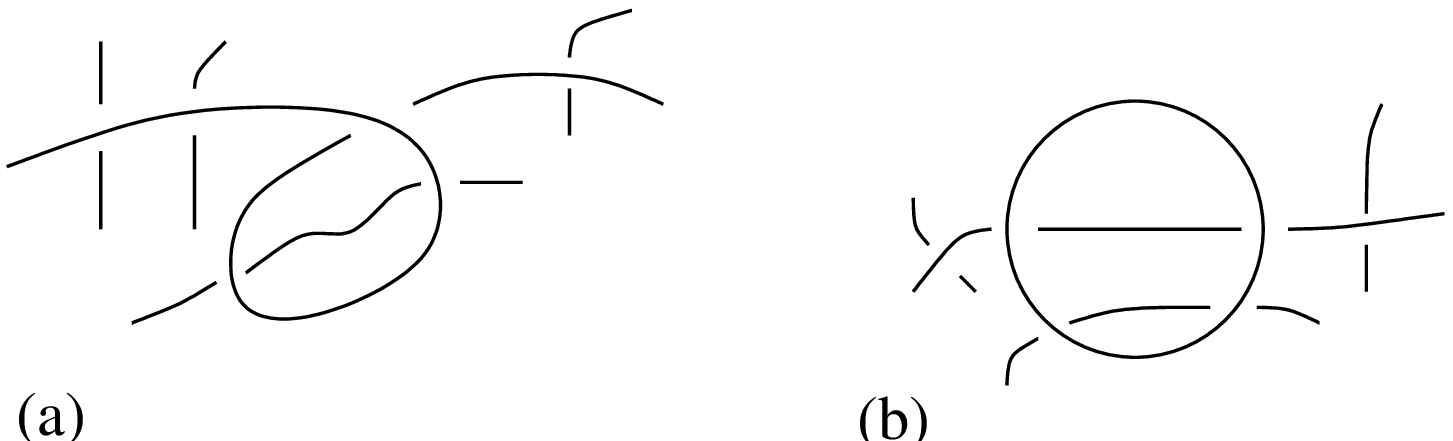,height=2.2cm}}}

\begin{center}
Fig.~4.5 \ Generalized bridges of length 6.
\end{center}

A generalized bridge does not have to be connected.
If a component belongs to the bridge then we can extend the bridge by moving 
along another component (see Fig.~4.5(b)).
The length of a generalized bridge is equal to the number of crossings
in the bridge.

Let $b'(L)$ be the maximal length of a generalized bridge in
a diagram $L$. Prove that $\deg Q_L \leq n(L) - b'(L)$.
\end{exercise}

\begin{exercise}\label{4:3.11}
Assume that a connected diagram $L$ decomposes into a connected sum 
of prime diagrams $\row{L}{k}$. Prove that $\deg Q_L\leq n(L) - \sum^{k}_{i=1} b(L_i)$. 
Generalize this claim for the case of disconnected diagrams.

Hint. Apply formulas: $Q_{L_1\sqcup L_2} = ( {\frac{2}{x}}
- 1)Q_{L_1}Q_{L_2}$ and $Q_{L_1 \kwad L_2} = Q_{L_1}Q_{L_2}$. 
\end{exercise}

The Kidwell result was generalized to the Kauffman polynomial,
by Thistlethwaite \cite{This-4}.

\begin{theorem}\label{4:3.12}
Let $L$ be a diagram of a link and let  $\Lambda_L (a,x) =
\sum u_{r,s}a^r x^s$ be its Kauffman polynomial 
If $u_{r,s}\neq 0$ then $|r|+s\leq n(L)$ and $s\leq n(L)-b(L)$. 
\end{theorem}

Proof. The proof of the second inequality is similar to that of
\ref{4:3.2} for $\Lambda_L$. In order to prove the first one
we note that it holds for descending diagrams. If 
$L$ is a descending diagram of a link with $c$ components
then  $\Lambda_L (a,x) = a^{ Tait(L)} ( \frac{a+a^{-1}}{x} - 1)^{c-1}$. 
Since $|Tait(L)|\leq n(L)$ it follows that $|r|+s\leq n(L)$. 
Now to conclude the proof we apply induction with respect
to the number of crossings and the number of ``bad'' crossings;
we note that, if the theorem is true for 
$L_{\psfig{figure=L+nmaly.eps}}$, $L_{\psfig{figure=L0nmaly.eps}}$ and
$L_{\psfig{figure=Linftynmaly.eps}}$ then it is true for
$L_{\psfig{figure=L-nmaly.eps}}$.

\begin{theorem}\label{4:3.13}
%\begin{theorem}\label{V.4.14}
If $L$ is a connected prime alternating diagram (with $n(L)\neq 1$)
then the coefficient of $x^{n-1}$ in $\Lambda_L$ is equal to $\alpha
(a+a^{-1})$ where $\alpha \geq 1$. 
\end{theorem}

Proof. Let us consider the formula 
$\Lambda_{L_{\psfig{figure=L+nmaly.eps}}} =
x(\Lambda_{L_{\psfig{figure=L0nmaly.eps}}} +
\Lambda_{L_{\psfig{figure=Linftynmaly.eps}}}) -
\Lambda_{L_{\psfig{figure=L-nmaly.eps}}}.$
Following Kidwell we note that to get the coefficient of
$x^{n-1}$ we can ignore $\Lambda_{L_{\psfig{figure=L-nmaly.eps}}}$.
Therefore, the equality is similar to that of Tutte polynomial
or Kauffman bracket, see Section 1.

Now if we present the diagram as a positive graph
with Tutte polynomial 
$\chi _{G(L)} = \sum v_{ij}x^i y^j$ (note that the variable $x$ is not the same
as in $\Lambda$) then we get $u_{1,n(L) - 1} = v_{0,1}\geq 1$ which concludes
the proof of Theorem \ref{4:3.13}.

%\begin{corollary}\label{4:3.14}
\begin{corollary}[The second Tait Conjecture]\label{V.4.15}\ \\
The writhe number $(Tait(L))$ of a prime connected
alternating diagram ($n(L)\neq 1$) is an invariant of isotopy of the link.
\end{corollary}

Proof. We have $F_L (a,x) = a^{-Tait (L)}\Lambda_L (a,x)$ and therefore
the coefficient of $x^{n(L) -1}$ in $F_L$ is equal to
$\alpha (a+a^{-1})a^{-Tait (L)}$. Hence $Tait (L)$ 
is an invariant of isotopy of the link.

\begin{exercise}\label{V.4.16}
Recall that a diagram of a link is called reduced if it contains 
no nugatory crossing ({\psfig{figure=nugatory.eps,height=0.4cm}}).

Prove that Corollary 4.15 holds for reduced
alternating diagrams of links.
\end{exercise}

\begin{exercise}\label{V.4.17}
Prove a version of Theorem \ref{4:3.2} for 
Jones-Conway polynomial
(i.e.~$\deg_Z P_L (a,z)\leq n(L) - b(L)$).
\end{exercise}

In the subsequent section we will describe further 
applications of the Kauffman polynomial. In particular, 
we will deal with a class
of diagrams generalizing alternating diagrams.
We will also describe relations of ``boundary'' coefficients
of the Kauffman polynomial with coefficients of Tutte polynomial
associated to a diagram of a link, and deduce from this that 
computing Kauffman polynomial is NP-hard.

%%%%%%%%%%%%%%%%%%%%%%%%%%%%%%%%%%%%%%%%%%%%%%%%%%%%%%%%%%%%%%%%%%%%%%

\section{Kauffman polynomial of adequate links.}\label{V.5}

In this section we consider relations between the Kauffman 
and Tutte polynomials. 
In the previous section we refrained  from extensive use of properties
of graphs related to links. Graph theory allows, however, a substantial
extention of previous results. Our approach is based on the paper of
Thistlethwaite \cite{This-5}. 

First we have to extend the language of graph theory which we
introduced in Section 1. For a 2-color graph $G$ let $G_b$ denote
a subgraph of $G$ which consists of vertices and
black edges of $G$. Similarly, $G_w$ denotes a graph obtained from
$G$ by removing black edges.
Then a quotient graph $\overline{G}_b$ is obtained by identifying
vertices of $G_b$ which were connected by white edges of $G$.
In other words, $\overline{G}_b$ is obtained from $G$ by removing
white loops and collapsing other white edges of $G$.
Similarly we define $\overline{G}_w$. 

An edge of a graph is called proper if it is neither a loop nor an isthmus.
A crossing in a diagram of a link to which we associate an improper
edge of the associated planar graph is called nugatory (it agrees with
the previous definition of a nugatory crossing) (see Fig. 5.1). 

%\end{document}

% \vspace*{2in}\centerline{\Psfig{figure=Rys.4.1}}
\centerline{{\psfig{figure=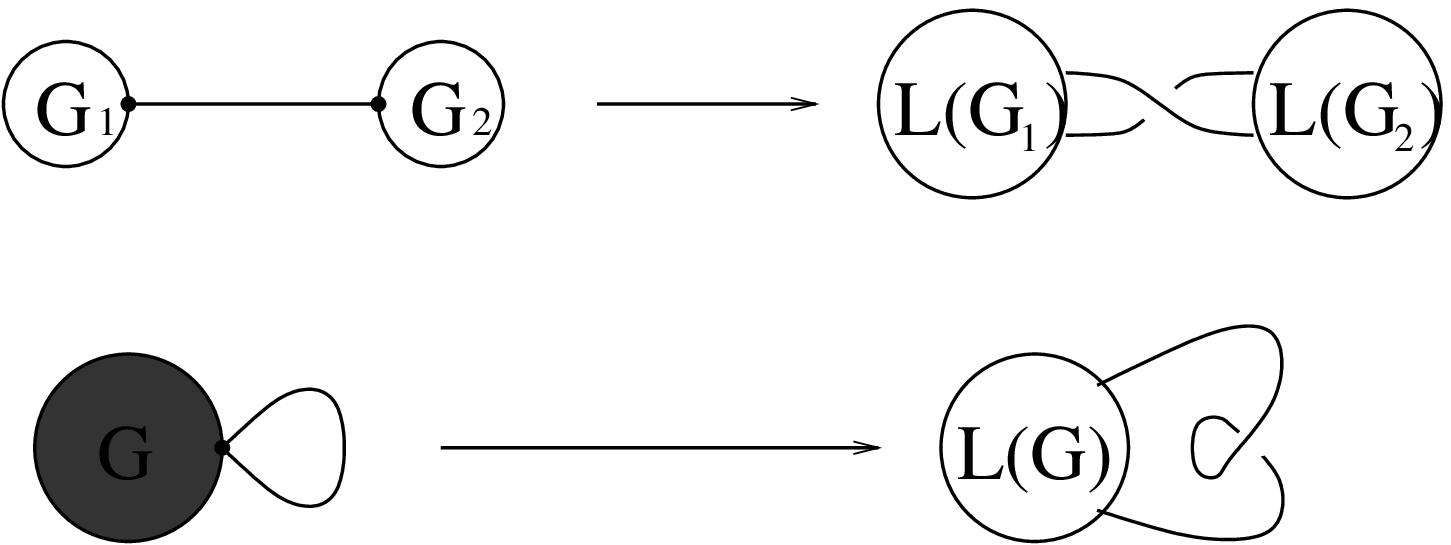,height=3.5cm}}}
\begin{center}
Fig.~5.1
\end{center}

If $e$ is an edge of a graph $G$, then $G-e$ (respectively, $G/e$ or
$G^\sigma_e$) denotes a graph obtained from $G$ by deleting (respectively,
contracting, or changing of the color of) the edge $e$. 
%(???ta definicja powinna pojawic sie w \S 1 ????)

If $G$ is a connected graph then it can be decomposed uniquely
into 2-connected components (called blocks), $G=G_1*G_2*...*G_k$. 
For a plane graph associated to a diagram of a link the above decomposition
is related to the decomposition of the associated diagram to prime diagrams
(the link is then a connected sum of them). We count here a loop as a block 
(the corresponding diagram is however that of a trivial knot with 
a nugatory crossing).

We have shown in Section 4 (Theorem 4.13)
 that the coefficient $u_{r,s}$ in the polynomial
$\Lambda_L (a,z) = \sum u_{r,s}a^r x^s$ 
is non-zero only if $r+s\leq n(L)$ and $-r+s\leq n(L)$. 
Now we will show that in some situations the above inequalities 
become strict and this will be a starting point to the proof of
the main theorem of this section, Theorem 5.1. 
Coefficients $u_{r,s}$ of the polynomial $\Lambda_L (a,x)$ 
will be called exterior if either $r+s=n$ or $-r+s=n$. 
Now let us define two auxiliary polynomials (of exterior coefficients):
$$\phi^+_L (t) = \sum_{i} u_{i,n-i}t^i,\ \ \ \ \ 
\ \ \ \phi ^-_L(t) = \sum_{i}
u_{-i,n-i} t^i.$$ 
Let us note that $\phi^+_L$ and $\phi^-_L$ are true polynomials in $t$
i.e.~they contain no negative degree monomial of $t$.

\begin{theorem}\label{4:4.1}
Let $L$ be a connected diagram of a link with at least one crossing.
Suppose that $G$ is a graph associated with $L$.
Then
$$\phi^+_L(t) = \chi_{G_b} (0,t)\chi_{\overline{G}_w}(t,0)$$
and
$$\phi^-_L (t) = \chi _{G_w} (0,t)\chi_{\overline{G}_b}(t,0).$$
\end{theorem}

The following properties of polynomials $\phi^+$ and $\phi^-$
(stated below for $\phi^+$ only) can be easily derived from properties
of Tutte polynomial.

\begin{corollary}\label{V.5.2}
Let $L$ be a connected diagram of a link
with at least one crossing. Then
\begin{enumerate}
\item[(i)] All coefficients of $\phi^+_L (t)$ are non-negative.

\item[(ii)]  $\phi^+_L (t)\neq 0$ if and only if $G_b$ has no isthmus
and $\overline{G}_w$ has no loop.

\item[(iii)] If  $\phi^+_L (t)\neq 0$, then the lowest degree term of 
$\phi^+_L (t)$ is equal to $\alpha t^b$, where $\alpha\geq 1$ 
and $b$ is the number of blocks (2-connected components)
which have at least one edge in $G_b$ and $\overline{G}_w$

\item[(iv)] If $\phi^+_L (t)\neq 0$ then the highest degree term
of $\phi^+_L (t)$ is of the form $t^{p_1(G_b)+d(\overline{G}_w)}$, where
$d(\overline{G}_w)$ is the number of edges in a tree spanning 
$\overline{G}_w$ (obviously $d(\overline{G}_w) =
|V(\overline{G}_w| -1$). 

\item[(v)] The degree of $\phi^+_L (t)$ is not greater than the number
of crossings in $L$ and it is smaller than $n(L)$ if  
$L$ has a non-nugatory crossing.
\end{enumerate}
\end{corollary}

Proof of Corollary.
\begin{enumerate}
\item[(i)]  It follows from Corollary  1.9(i).
\item[(ii)] follows from Corollary 1.9 (ii) and (iii).
\item[(iii)]  If a graph has at least one edge then $v_{0,0} = 0$. 
For a 2-connected graph (with at least two edges) we have
$v_{0,1} = v_{1,0} > 0$ and (iii) follows.
\item[(iv)] It follows from Exercise 1.11 (i) and (iii).

\item[(v)] We analyze the right side of the equality (from Theorem 5.1)
$\phi^+_L (t) = \chi_{G_b} (0,t)\chi_{\overline{G}_w} (t,0)$.
%ilo\'sci kraw\c edzi w $G_b$ a w $\chi_{\overline{G}_w} (t,0)$ ilo\'sci
The degree of $\chi _{G_b} (0,t)$  does not exceed
the number of edges in $G_b$ and, similarly, 
the degree of $\chi_{\overline{G}_w} (t,0)$ is not greater
than the number of edges in $G_{\overline{G}_{w}}$.
This implies the first part of (v).
If $L$ has a non-nugatory crossing 
then $G=G(L)$ does not consists of isthmuses and loops only.
On the other hand, the degree of $\chi_{G_b}(0,t)$ 
is equal to the number of edges in $G_b$ \ iff
$G_b$ is composed of loops alone, and the degree of 
$\chi_{\overline{G}_w}(t,0)$ can be equal to the number of edges in
$\overline{G}_w$\ if \ $\overline{G}_w$ is a tree. However, 
in this case $G$ would be a white tree with black loops.
%$G$ by\l oby drzewem (jasnym) z (ciemnymi) p\c etlami.???????mam wrazenie,
\end{enumerate}

To prove Theorem 5.1, let us begin by proving 
a lemma which will be the key to our induction argument.

%\begin{lemma}\label{V.5.3}
\begin{lemma}\label{4:4.3}
Let $L$ be a connected diagram of a link which has the following properties:
\begin{enumerate}
\item[(i)] $L$ is reduced (i.e.~it has no nugatory crossing).
%\mbox{$\ \ $}\ \ \ )

\item[(ii)] $L$ contains a bridge $B$ of length greater than 1.

\item [(iii)] We can change  crossings in $B$ so that $L$ is changed to
 an alternating diagram. 
\end{enumerate}
Then  $\phi^+_L (t) =  \phi^-_L (t) = 0$.
\end{lemma}

Let us note that the condition (i) in the above Lemma can not
be removed as it is apparent for the diagram in Fig.~5.2
for which $\Lambda_L (a,x) = a^2$, and therefore $\phi^+_L (t) = t^2$. 

%%% \vspace*{1.5in}\centerline{\Psfig{figure=Rys.4.2}}
%\centerline{{\psfig{figure=4-4-2.eps,height=3.1cm}}}
\centerline{{\psfig{figure=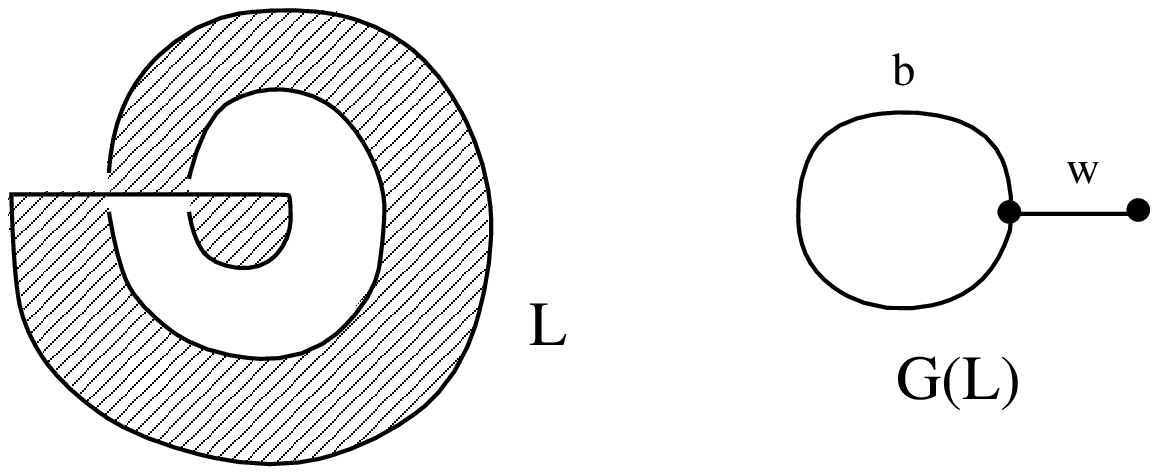,height=3.1cm}}}
\begin{center}
Fig.~5.2
\end{center}

Proof. Let us assume that $L$ satisfies the assumption of the lemma.
We define complexity of $L$ to be the ordered pair 
$(n,k)$, where $n = n(L)$ is the number of crossings in $L$ and 
$k$ is the number of crossings which have to be changed (overcrossing to
undercrossing) in order to get a diagram $L'$ 
which satisfies assumptions of the lemma
and such that its longest bridge (that is the bridge which cannot 
be further extended), say $B'$ contains the bridge $B$,
see Fig.~5.3.\\
\ \\
% \vspace*{1.5in}\centerline{\Psfig{figure=Rys.4.3}}
%\centerline{{\psfig{figure=Rys4-4-3.eps,height=3.1cm}}}
\centerline{{\psfig{figure=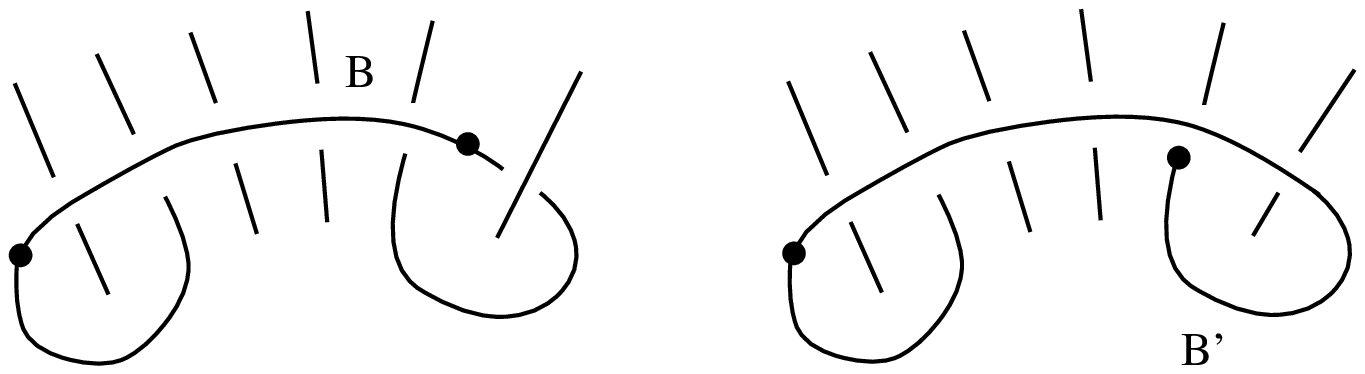,height=3.1cm}}}
\begin{center}
Fig.~5.3
%poprawic -przedluzyc luk po prawej stronie,-- extend the arc
\end{center}

Let us consider a lexicographic order of pairs $(n,k)$. 
If $n=2$ then the link $L$ is trivial with the diagram 
%{\psfig{figure=T-H.eps,height=0.4cm}},
\parbox{0.8cm}{\psfig{figure=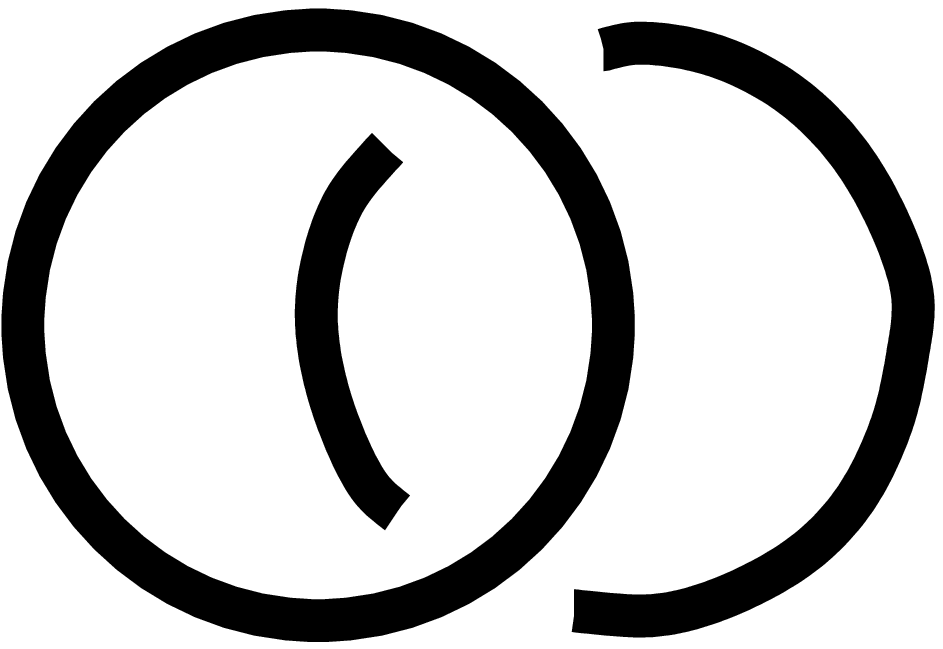,height=0.5cm}}, 
so the lemma holds.

If $k=0$ then $B$ is of maximal length already 
and therefore $L'=L$ is regularly isotopic to a diagram $L''$ 
which contains fewer crossings ($B$ is a closed component or 
it contains a loop whose shrinking can eliminate crossings; we use the fact 
that $L$ has no nugatory crossings). 
Since moreover $|r|+s\leq n(L'')< n(L)$ 
(notation as in Theorem 5.1) it follows that
 $$\phi^+_L (t) =  \phi^-_L (t) = 0.$$
Now let us consider a diagram  $L$ which satisfies the assumptions of Lemma
\ref{4:4.3} and which has complexity $(n,k)$ with $n>2$, $k>0$, 
and let us assume that the lemma is true for diagrams of smaller
complexity.
There exists a crossing $p$ in $L$ at one of the ends of the bridge
$B$ such that the change of $p$ yields a diagram
$L^p$ of smaller complexity  than this of $L$.
Let us consider diagrams $L_{{\psfig{figure=L0nmaly.eps}}}$ and 
$L_{{\psfig{figure=Linftynmaly.eps}}}$ which we obtain from
$L$ by smoothing the crossing $p$.
If we show that 
$\phi^+_{L_{{\psfig{figure=L0nmaly.eps}}}} (t) =  
\phi^-_{L_{{\psfig{figure=L0nmaly.eps}}}} (t) =  
\phi^+_{L_{{\psfig{figure=Linftynmaly.eps}}}} (t) =
\phi^-_{L_{{\psfig{figure=Linftynmaly.eps}}}} (t) = 0$ 
then applying equality 
$\Lambda_{L_{\psfig{figure=L+nmaly.eps}}} 
+ \Lambda_{L_{\psfig{figure=L-nmaly.eps}}} = 
x(\Lambda_{L_{{\psfig{figure=L0nmaly.eps}}}} + 
\Lambda_{L_{{\psfig{figure=Linftynmaly.eps}}}})$ 
we will conclude the proof.
Thus let us focus our attention on the diagram
$L_{{\psfig{figure=L0nmaly.eps}}}$ 
(the case $L_{{\psfig{figure=Linftynmaly.eps}}}$ is similar).
The diagram $L_{{\psfig{figure=L0nmaly.eps}}}$ is connected (because $L$ 
was reduced). If $L_{{\psfig{figure=L0nmaly.eps}}}$ is reduced 
then it satisfies the conditions of our lemma 
and since it has fewer crossings than $L$ thus, because of inductive 
assumption, $\phi^+_{L_{{\psfig{figure=L0nmaly.eps}}}} (t) =  
\phi^-_{L_{{\psfig{figure=Linftynmaly.eps}}}} (t) = 0$. 
If now the diagram $L_{{\psfig{figure=L0nmaly.eps}}}$ 
has a nugatory crossing $q$ then one
of the following two cases occurs. 

%\vspace*{2in}\centerline{\Psfig{figure=Rys.4.4}}
% \begin{center} Fig.~5.3 \end{center}

\begin{enumerate}
\item [(1)] 
The crossing $q$ is on the bridge $B$ but it is not next to $p$.
Then there exists a simple closed curve $C$ which meets $L$ in 4 points
(as shown in Fig. 5.4) and the diagram 
$L_{{\psfig{figure=L0nmaly.eps}}}$ may be reduced
by a regular isotopy to a diagram with fewer crossings (we use the fact that 
there was a crossing between $p$ and $q$ on the bridge).
Therefore $\phi^+_{L} (t) =  \phi^-_{L} (t) = 0$.
\\ \ \\
%\vspace*{2in}
%\centerline{{\psfig{figure=Rys4-4-4.eps,height=3.6cm}}}
\centerline{{\psfig{figure=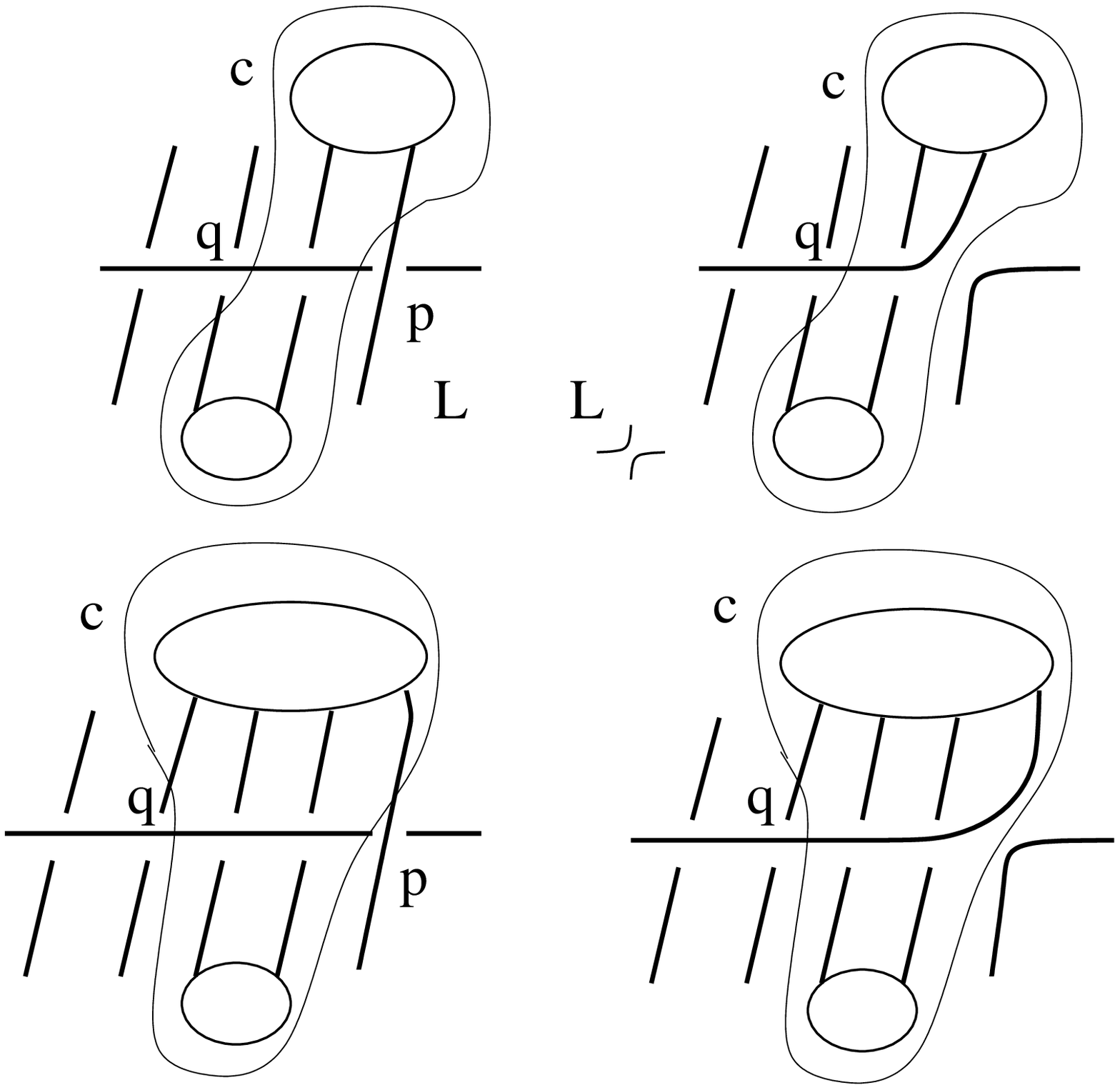,height=6.6cm}}}
\begin{center}
Fig.~5.4
\end{center}

\item [(2)] The crossing $q$ is next to $p$ on the bridge
 and no other crossing 
on the bridge $B$ is nugatory in $L_{{\psfig{figure=L0nmaly.eps}}}$. 
We will see that we can 
reduce $L_{{\psfig{figure=L0nmaly.eps}}}$, 
by removing nugatory crossings, to a diagram 
$L_{{\psfig{figure=L0nmaly.eps}}}'$ 
which either satisfies assumptions of the lemma and has
smaller complexity or $\phi^+_{L_{{\psfig{figure=L0nmaly.eps}} }'} (t) = 
\phi^-_{L_{{\psfig{figure=Linftynmaly.eps}}}'} (t) = 0$
since it is regularly isotopic to a diagram with fewer crossings.
So first let us remove all nugatory crossings 
(\parbox{2.8cm}{\psfig{figure=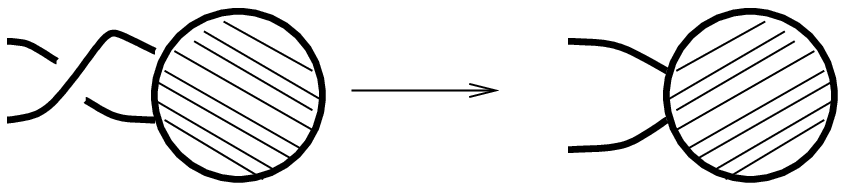,height=0.6cm}} ), except $q$,
and let us call the result 
$L_{{\psfig{figure=L0nmaly.eps}}}''$. 
Certainly this operation does not change the non-nugatory status 
of other crossings (in particular in the resulting
diagram all crossings except $q$ are not nugatory).
Let $q'$ be the first crossing after $q$ met when one travels along 
the diagram $L_{{\psfig{figure=L0nmaly.eps}}}''$ away from $p$.
%If there exists a curve in  $L_{{\psfig{figure=L0nmaly.eps}}}''$ 
%which connects the undercrossing (tunnel) in $q$ to another 
If $q'$ is an undercrossing
under the bridge $B$ (as shown in Fig.~5.5(a)) 
then $L_{{\psfig{figure=L0nmaly.eps}}}''$ 
is regularly isotopic to a diagram with fewer crossings
and we can consider $L_{{\psfig{figure=L0nmaly.eps}}}' = 
L_{{\psfig{figure=L0nmaly.eps}}}''$.  

%\vspace*{1.5in}\centerline{\Psfig{figure=Fig.V.5.5}}
\centerline{{\psfig{figure=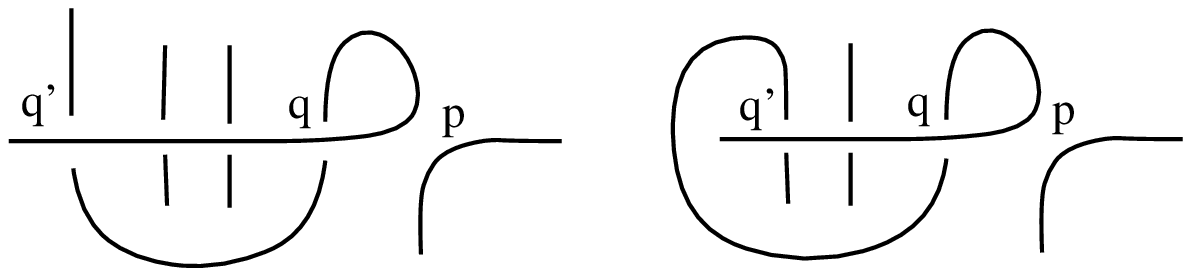,height=3.0cm}}}
\centerline{(a)}
\ \\
\centerline{{\psfig{figure=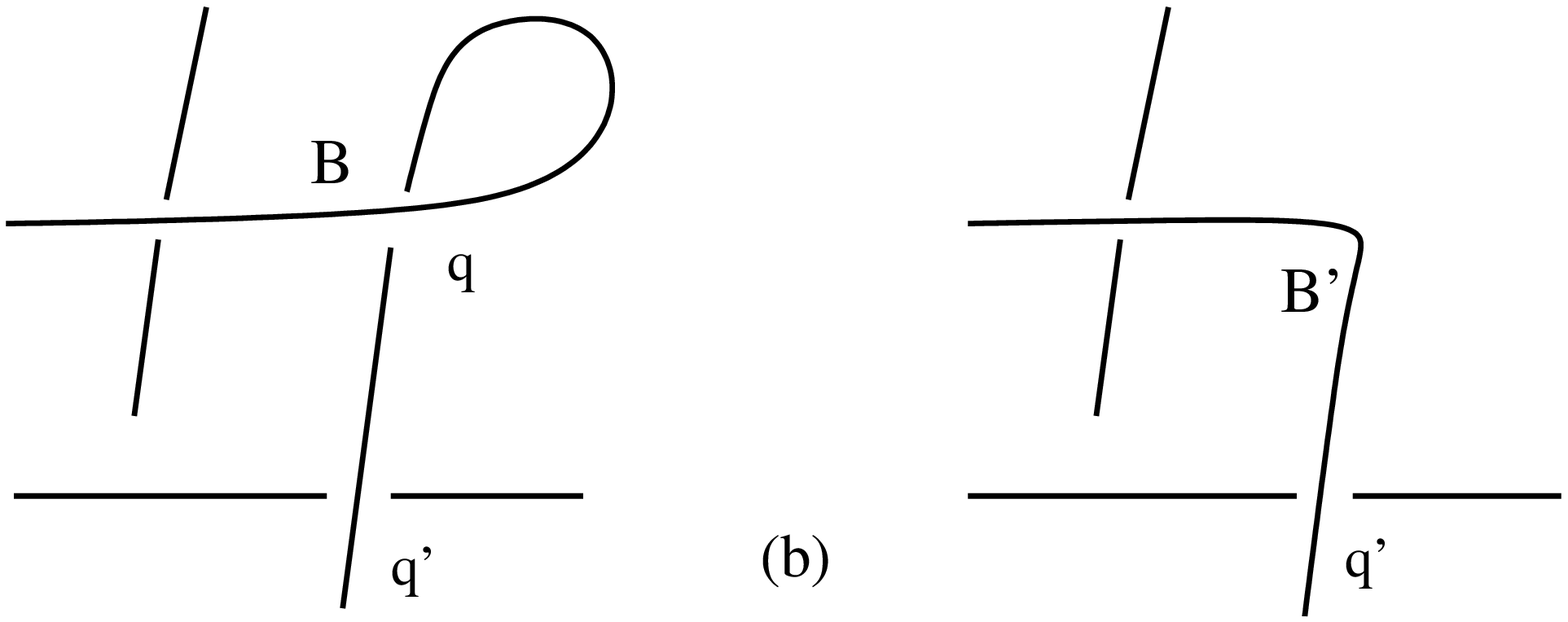,height=3.8cm}}\ \ \ \ \ \ \ \ }
\begin{center}
%Fig. 5.5  CHANGE (it was the same as 5.7) 
Fig. 5.5 
\end{center}
Otherwise the condition (iii) of the lemma implies that the 
diagram obtained from
$L_{{\psfig{figure=L0nmaly.eps}}}''$ by untwisting the crossing $q$ 
has a bridge $B'$ of the same length as $B$ and thus it 
satisfies assumptions of Lemma 5.3 with smaller complexity than 
$L$ (Fig. 5.5(b).
Let us call the diagram, which was obtained by removing all nugatory 
crossings of 
$L_{{\psfig{figure=L0nmaly.eps}}}$, 
by $L_{{\psfig{figure=L0nmaly.eps}}}'$. 
We note that  $\phi^+_{L_{{\psfig{figure=L0nmaly.eps}}}'} (t) =
\phi^-_{L_{{\psfig{figure=L0nmaly.eps}}}'} (t) = 0$. 
Since the operation of ``untwisting''  
of a nugatory crossing (as well as its inverse) is related to multiplying
$\Lambda$ by either $a$ or $a^{-1}$, then all the time we are getting
 $\phi^+ (t) = \phi^- (t) = 0$. 
This concludes the proof of Lemma 5.3.
\end{enumerate}
\begin{exercise}\label{V.5.4}
Prove Lemma 5.3 by using the notion of generalized bridge.
\end{exercise}

%\end{document}

\begin{exercise}\label{V.5.5}
Let $\hat{P}_L (a,z) = a^{Tait (L)} P_L (a,z)$ be a polynomial invariant
of regular isotopy of oriented diagrams\footnote{This version of Homflypt 
polynomial was first considered in 1985 by L.Kauffman.}.
The polynomial $\hat{P}_L(a,z)$ can be defined by the following properties:

\begin{enumerate}

\item [(1)] $\hat{P}_\bigcirc (a,z) =1$

\item [(2)] 
$\hat{P}_{\psfig{figure=R+maly.eps}} 
(a,z) =
a^{-1}\hat{P}_{{\psfig{figure=ver.eps,height=0.4cm}}} (a,z),
\hat{P}_{\psfig{figure=R-maly.eps}} = 
a\hat{P}_{\psfig{figure=ver.eps,height=0.4cm}} (a,z)$

\item[(3)] $\hat{P}_{L_+} (a,z) +\hat{P}_{L_-} (a,z) = z\hat{P}_{L_0} (a,z)$.
\end{enumerate}

Prove that Lemma 5.3 is true for the polynomial
$\hat{P}_{L} (a,z)$. That is, if $L$ satisfies the assumptions of the lemma
and we write  $\hat{P}_{L} (a,z) = \sum
e_{ij}a^i z^j$ then for $|i|+j = n(L)$ we have $e_{ij} = 0$. 
\end{exercise}

Before we conclude the proof of Theorem 5.1 we recall
the notion of adequate diagram and we prove its basic properties.

In this section we assume that for an adequate diagram $L$ any 
component of $L$ has at least one crossing. This is a technical 
condition which will be not assumed in other chapters of the book.

%\begin{definition}\label{V.5.6}
\begin{definition}\label{4:4.6}
A diagram of an unoriented link $L$
is called $+$ adequate if:
\begin{enumerate}
 \item [(1)] Any component of $L$ has at least one crossing,

 \item [(2)] If we modify $L$ to $s_+L$ by smoothing every crossing of $L$
 according to the rule
 \parbox{0.6cm}{\psfig{figure=L-nmaly.eps,height=0.6cm}}
$\rightarrow$ 
\parbox{0.6cm}{\psfig{figure=Linftynmaly.eps,height=0.6cm}} 
 (equivalently  
 \parbox{0.6cm}{\psfig{figure=L+nmaly.eps,height=0.6cm}} $\rightarrow$ 
\parbox{0.6cm}{\psfig{figure=L0nmaly.eps,height=0.6cm}})
 then the two arcs of
 \parbox{0.6cm}{\psfig{figure=Linftynmaly.eps,height=0.6cm}} 
 belong in  $s_+(L)$ to two different components.
\end{enumerate}
 
If the mirror image $\bar L$ of $L$ is $+$-adequate
then $L$ is called $-$-adequate.

A diagram is called semi-adequate if it is $+$ or $-$ adequate
and it is called adequate if it is both $+$ and $-$ adequate.
By negation, a diagram is $+$ inadequate (respectively $-$ inadequate)
if it is not $+$ (resp.~$-$) adequate. A diagram is inadequate if 
it is not semi-adequate.
\end{definition}

The following lemma provides an interpretation of adequate diagrams in
terms of associated graphs. A Kauffman state $s=sL$ of $L$, is a function 
from crossings of $L$ to $\{+,-\}$ with the convention that $s(v)=+$ 
corresponds to the smoothing 
\parbox{0.6cm}{\psfig{figure=L+nmaly.eps,height=0.6cm}} $\rightarrow$
\parbox{0.6cm}{\psfig{figure=L0nmaly.eps,height=0.6cm}} and 
$s(v)=-$ corresponds to the smoothing
\parbox{0.6cm}{\psfig{figure=L+nmaly.eps,height=0.6cm}} $\rightarrow$
\parbox{0.6cm}{\psfig{figure=Linftynmaly.eps,height=0.6cm}}.
By $|s|=|sL|$, as usually,  we denote the number of components of $L_s$,
where $L_s$ is system of circles obtained by smoothing $L$ according to 
the state $sL$.\footnote{Already J. Listing \cite{Lis} was decorating 
corners of the crossing by variables $\delta$ (deotropic) and 
$\lambda$ (leotropic) 
\parbox{0.6cm}{\psfig{figure=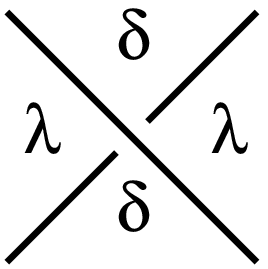,height=0.6cm}} and he observed that 
if a (connected) alternating diagram $D$ is alternating then for any 
region $R^2 -D$ corners of the region gave the same label, all $\lambda$ or 
all $\delta$. L. Kauffman was decorating corners by $A$ and $B$ \ 
\parbox{0.6cm}{\psfig{figure=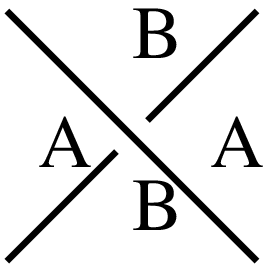,height=0.6cm}}\ and used it to 
construct (Kauffman) bracket polynomial. It is natural to 
say that a state $s$ associates to every crossing a marker $A$ 
(resp. $\lambda$)
\ \parbox{0.6cm}{\psfig{figure=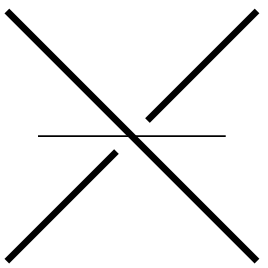,height=0.6cm}}\  
or $B$ (resp. $\delta$) \ 
\parbox{0.6cm}{\psfig{figure=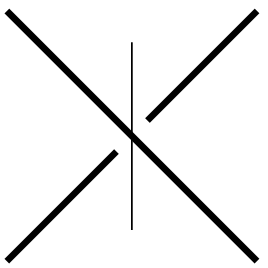,height=0.6cm}}\ . 
In Section 1 we use $b$ (black)-markers and 
$w$ (white)-markers because of black and white checkerboard coloring 
of regions of $R^2-D$. M. Thistlethwaite is using $+$ or $-$ markings so 
we follow his notation in this Section as well as in Chapter X (see Figure 
X.1.1), still, however, for corresponding graph edges we use black and white 
colors here. 
Similarly, the diagram obtained from a diagram $D$ by smoothing 
it crossings by ``applying" markers of a  Kauffman state $s$  is 
denoted by $D_s$ or $sD$ depending on the author.}. 

\begin{lemma}\label{4:4.7}
Let $G = G(L)$ be a planar graph associated to $L$ via some checkerboard
 coloring of the plane containing $L$. Then:

\begin{enumerate}
\item [(1)]
 $|s_+L| = p_0(G_b) +p_1(G_b)$\ and \ $|s_-L| = p_0(G_w) +p_1(G_w)$

\item [(2)] The diagram $L$ is $+$ adequate if and only if
$G_b$ has no isthmus and  $\overline{G}_w$ has no loop.

\end{enumerate}
\end{lemma}

Proof.
\begin{enumerate}

\item [(1)]
Notice that the modification of crossings of $L$ to $s_+L$ 
is related to collapsing of black edges (in the sense of $G//e$)
 and removing of white edges 
of the associated graph $G(L)$ (see Fig.~5.6).
\\ \ \\
% \vspace*{2in}\centerline{\Psfig{figure=Rys.4.6}}
\centerline{{\psfig{figure=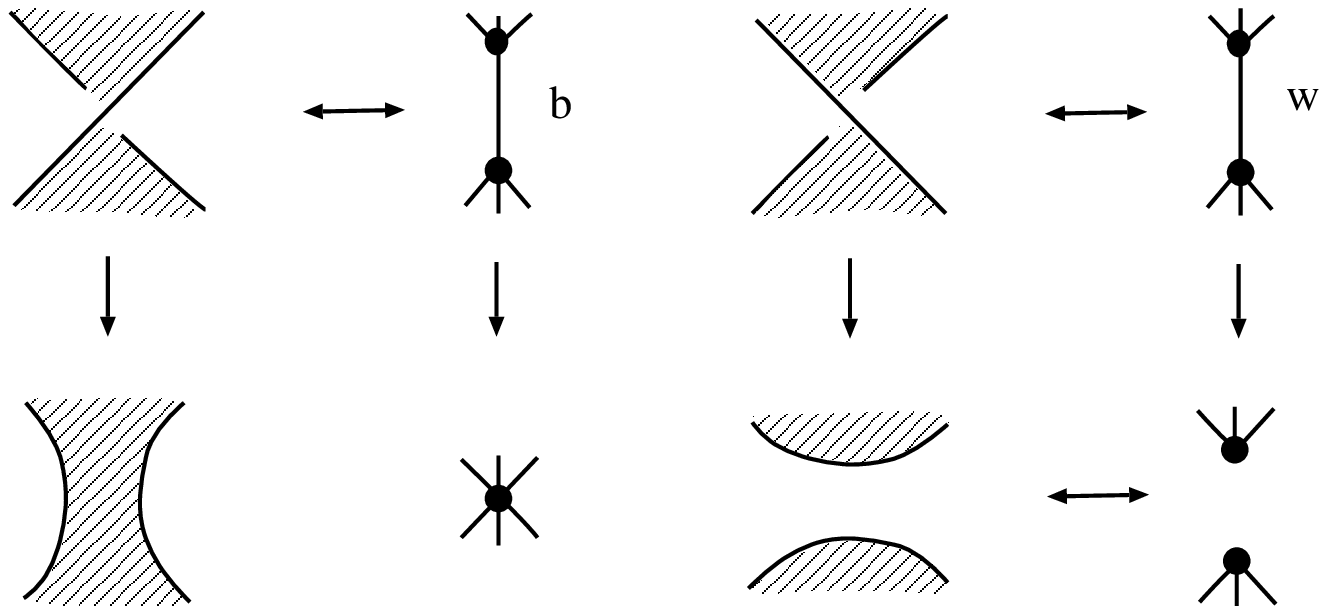,height=3.0cm}}}  
\begin{center}
Fig.~5.6
\end{center}

According to our convention the contracting of all edges in a connected
graph $\Gamma$ yields $p_1(\Gamma)+1$ vertices\footnote{From the point 
of view of Knot Theory it is convenient to consider contracting $G//e$ 
which agrees with standard contracting, $G/e$ for $e$ not being a loop.
If $e$ is a loop then $G//e$ results in the graph $G-e=G/e$ with one 
additional ``free" vertex 
\parbox{2.2cm}{\psfig{figure=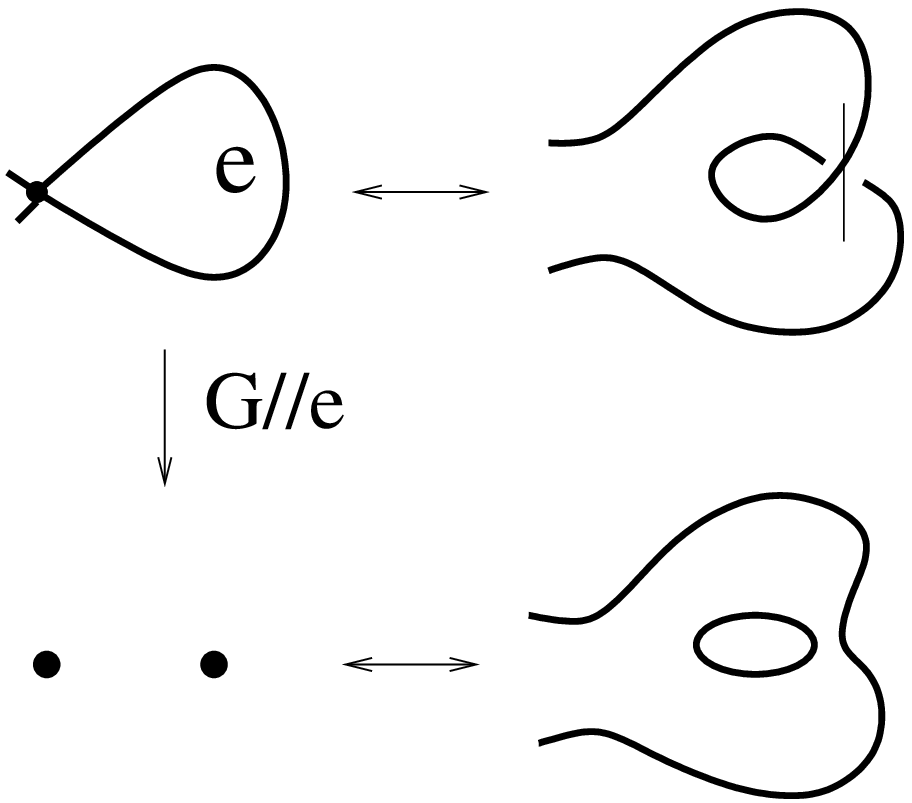,height=1.9cm}}.}.

\item [(2)]
 Notice that the condition (2) of Definition 5.6 is equivalent to the
statement that 
%the diagram $L$ is $+$ adequate if and only if 
$|s_+L^\sigma| < |s_+L|$ for any diagram $L^\sigma$ which is obtained
from $L$ by changing overcrossing to undercrossing at a crossing of $L$.
The change of $L \to L^\sigma$ corresponds to changing of the color of 
an edge in $G(L)$.
%The change of a color of an edge of $G(L)$ 
Such a change from black to white
decreases the sum $p_0(G_b) + p_1(G_b)$ unless the edge in question
is an isthmus in $G_b$ --- in such a case the number increases.
And conversely, if we change an edge from white to black 
then the number $p_0(G_b) + p_1(G_b)$ decreases if the edge
is an isthmus in $G_b$ and it increases otherwise
--- and the latter case is related to the situation when the edge 
in question becomes a loop in $\overline{G}_w$. 
Therefore the condition in (2) (that $G_b$ has no isthmus 
and  $\overline{G}_w$ has no loop) is necessary and sufficient for
the number $p_0(G_b) + p_1(G_b)$ to decrease when an edge changes its color.
>From this it follows that $L$ is $+$-adequate. 
%(???przetlumaczylem ostatnie zdanie doslownie, ale tak jak jest
\end{enumerate}

The respective criterion for $-$ adequate diagrams can be proved
similarly.

%\begin{lemma}\label{V.5.8}
\begin{lemma}\label{4:4.8}
Let $L$ be a connected $+$ inadequate (respectively, $-$ inadequate)
diagram such that $n(L)>0$.
Then $\varphi^+_L = 0$ (respectively, $\varphi^-_L = 0$).
\end{lemma}

We provide a proof for a $+$ inadequate diagram. Let $L$ be the diagram
in question and let $G=G(L)$ be a graph of $L$.
The proof is by induction with respect to the number $k$ of edges of
$G_b$ which are neither isthmuses nor loops in $G$.
If $k=0$ then, since $L$ is $+$ inadequate, 
it follows that $G$ contains either a black isthmus (isthums in $G_b$ 
is now an isthmus in $G$) or white loop (a loop in $\overline{G}_w$ is 
now a loop in $G_w$ so in $G$). In both cases 
the sign of 
the corresponding selfrossing in $L$ is negative,
see Fig.~5.7.

%\vspace*{1.8in}\centerline{\Psfig{figure=Rys.4.7}}
\centerline{{\psfig{figure=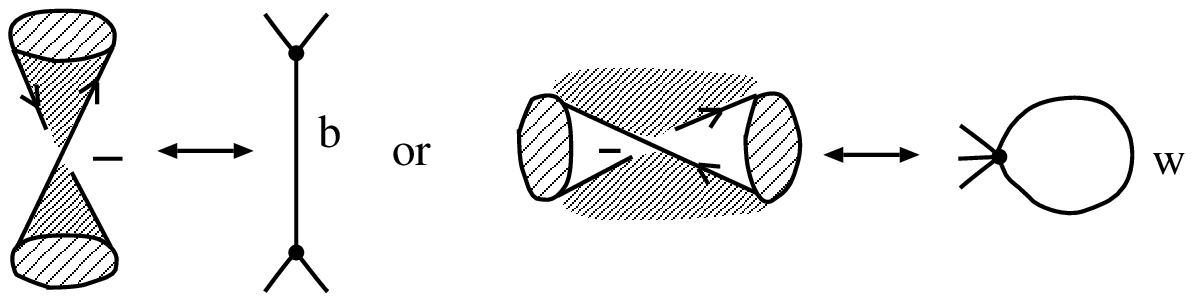,height=3.5cm}}}  
\begin{center}
Fig.~5.7
\end{center}

If we remove this edge (black isthmus or white loop) to get 
the associated diagram $L'$ with one less crossing
and $\Lambda_L'(a,x)= a \Lambda_L(a,x)$. Because in $L'$ we have 
$|r'| + s' \leq n(L') = n(L)-1$ therefore 
$|r| + s \leq n(L)-2$ and so $\phi^+_L = 0$. 

If $k=1$ then let $L_1$ denote a diagram obtained from $L$ by removing
(untwisting) all nugatory crossings. $L_1$ satisfies
the assumptions of Lemma \ref{4:4.3}. Clearly, $L_1$ has no nugatory
crossing. The graph $G_b$ has only one edge so if we change its color then 
we obtain a white graph yielding an alternating diagram and therefore the
edge in question yields the crossing being a part of a bridge 
on $L_1$ of length at least 2.  
Now by Lemma \ref{4:4.3}  we get $\varphi^+_{L_1} = 0$. 
On the other hand, the polynomial $\Lambda_L$ is obtained from
$\Lambda_{L_1}$ by multiplying by a power of $a$ of degree not exceeding
the number of crossings of $L$. 
Therefore $\varphi^+_{L} = 0$ as well.

Now suppose that $G_b$ has $k>1$ essential edges (i.e.~edges 
which are neither isthmuses nor loops) in $G$ 
and assume that Lemma \ref{4:4.8} is true for
connected $+$-inadequate diagrams (with graphs with smaller $k$).
The set of edges in $G_b$ which are essential in $G$ 
is denoted by $\varepsilon$ ($|\varepsilon|=k$).
First, let us consider the case when $\varepsilon$ contains an edge
$e$ which is essential not only in $G$ but also in $G_b$. 

It follows by Lemma \ref{4:4.7} (ii) that contraction of $e$ does not
change $+$-inadequacy of the graph 
(we do not loose neither any isthmus of $G_b$ nor any 
loop in $\overline{G}_w$), that is, $G/e$ is $+$-inadequate.
We have yet to prove that both $G-e$ and $G^\sigma_e$ 
are $+$-inadequate and then, because of our inductive assumption
and in view of the skein relation satisfied by polynomial $\Lambda (a,x)$ 
we will get $\varphi^+ = 0$. To this end, let us note that 
if $G_b$ has an isthmus then also $G_b - e$ and $G^\sigma_e$ have one 
and thus (because of Lemma 5.7(ii))
$G-e$ and $G^\sigma_e$ are $+$-inadequate. Therefore
we may assume that $G_b$ has no isthmus but $\overline{G}_w$ contains a loop.
Let us consider an edge $x$ in $G_w$ which is a loop in 
$\overline{G}_w$, which means that the ends of this edge can be joined 
by a path in $G_b -e$ ($e$ is essential in $G_b$). 
Hence $G-e$ and $G^\sigma_e$ are $+$-inadequate.
Finally, we are left with the case when $\varepsilon$ 
consists only of isthmuses of $G_b$. But, since  $k>1$, then for any edge
$e\in\varepsilon$, black subgraphs of $G/e$, $G-e$ and $G^\sigma_e$
have isthmuses as well. Therefore, any of these graphs is $+$-inadequate
and $\varphi^+_{L} = 0$. 

This concludes the proof of Lemma 5.8.

%\begin{corollary}\label{V.5.9}
\begin{corollary}\label{4:4.9}
Let $G$ be a $+$-adequate (respectively, $-$-adequate) graph. If $e$
is an edge of  $G$, then the graph $G^\sigma_e$ is $+$ inadequate 
(respectively, $-$ inadequate).
\end{corollary}

We present the proof for $+$-adequate diagrams.
The edge  $e$ is not an isthmus in $G_b$,
thus if $e\in G_b$ then $e$ is a loop in $\overline{(G^\sigma_e)_w}$
hence $G^\sigma_e$ is $+$-inadequate.
If $e\in G_w$ then, because $G$ is $+$-adequate, the $e$ is not 
a loop in $\overline{G}_w$ and thus $e$ is an isthmus in 
$(G^\sigma_e)_b$. Therefore $G^\sigma_e$ is $+$-inadequate.
This completes the proof of Corollary 5.9. We should remark here that 
Corollary 5.9 has natural explanation if we consider associated 
$+$-diagrams and recall that $+$-adequate can can be interpreted as 
saying that in $s_+D$ no circle touches itself, thus in $s_+L^{\sigma}$ 
one has a self-touching circle. In fact this property of reduced alternating 
diagrams was the main reason for defining $+$-adequate diagrams.

Now we can proceed with the proof of Theorem 5.1.
%\ref{4:4.1}. 
We will consider the case of $\varphi^+_L (t)$, the case
of $\varphi^-_L (t)$ is similar. 
The proof will be done by induction with respect
to the number of edges of the graph $G$. First, let us assume that
the connected diagram  $L$ is $+$-inadequate. By Lemma
\ref{4:4.8} we have $\varphi^+_L = 0$ and because of Lemma \ref{4:4.7} (ii) 
either $G_b$ contains an isthmus and $\chi_{G_b} (0,t) = 0$ or
$\overline{G}_w$ contains a loop and then $\chi_{\overline{G}_w} (t,0) = 0$. 
Therefore $\varphi^+_L = \chi _{G_b} (0,t)\cdot\chi_{\overline{G}_w} (t,0)$.

The next case to consider is when $L$ is $+$ adequate but neither 
 $G_b$ nor $\overline{G}_w$ contains a essential edge. This means that
$G$ consists of black loops and white isthmuses. 
If the number of loops
and isthmuses is denoted by $p$ and $q$, respectively,
then $\chi _{G_b} (0,t) = t^p$ and $\chi_{\overline{G}_w} (t,0) = t^q$. 
On the other hand, the diagram $L$, which is associated to such a graph $G$,
represents a trivial knot with $p+q$ positive twists
(c.f.~Fig.~5.8), hence $\Lambda_L (a,x) = a^{p+q}$ and the theorem is true 
in this case as well. 

%\vspace*{1.5in}\centerline{\Psfig{figure=Rys.4.8}}
\centerline{{\psfig{figure=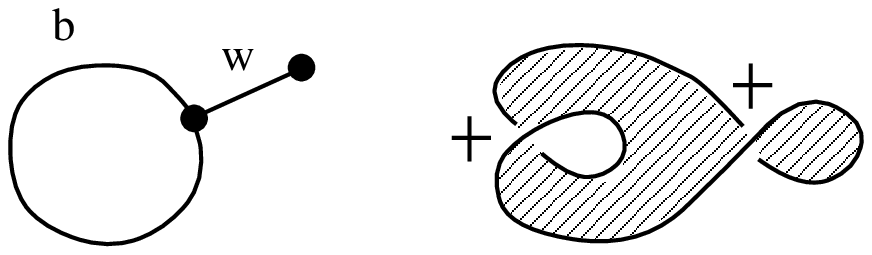,height=3.0cm}}}  
\begin{center}
Fig.~5.8
\end{center}

In particular, the theorem is true for any graph with one edge.
Now let us consider a $+$-adequate diagram $L$ and the associated graphs.
Let us assume that either $G_b$ or $\overline{G}_w$ 
has a essential edge $e$. By inductive assumption Theorem
\ref{4:4.1} is true for graphs with fewer edges.
Let us deal with the case $e\in G_b$ since the case of
$e\in \overline{G}_w$ can be done similarly.

By Corollary 5.9 the graph $G_e^{\sigma}$ is $+$-inadequate, so by
 Lemma 5.9 it follows that
$\varphi^+_{G^\sigma_e} = 0$. Therefore, by skein relation 
 for $\Lambda_L (a,x)$, we get 
$\varphi^+_G (t) = \varphi^+_{G-e} (t) + \varphi^+_{G/e} (t)$. 

Moreover $(G-e)_b = G_b - e$, $(G/e)_b = (G_b)/e$, and graphs
$\overline{(G/e)_w}$ and $\overline{G}_w$ are isomorphic (because $e$ is
essential in $G_b$). Now we get
\begin{eqnarray*}
&\varphi^+_G (t) &=\\
&  \varphi^+_{G-e} (t) + \varphi^+_{G/e} (t)& =\\
&\mbox{(because of inductive assumption)} &\\
&\chi_{(G-e)_b} (0,t)\cdot\chi_{\overline{(G-e)_w}} (t,0) +
\chi_{(G/e)_b} (0,t)\cdot\chi_{\overline{(G/e)_w}} (t,0)&=\\
&\chi_{(G_b-e)_b} (0,t)\cdot\chi_{\overline{G}_w} (t,0) +
\chi_{(G_b)/e} (0,t)\cdot\chi_{\overline{G}_w} (t,0)& =\\
&(\chi_{G_b-e} (0,t)+\chi_{(G_b)/e)_b} (0,t))\cdot\chi_{\overline{G}_w}
(t,0)& =\\
&\chi_{G_b} (0,t)\cdot\chi_{\overline{G}_w} (t,0).&\\
\end{eqnarray*}

This concludes the proof of Theorem \ref{4:4.1}.

If we apply chromatic polynomial $C(G,\lambda)$ (see Section 1)
then we can reformulate Theorem 5.1 in the following
way.

\begin{corollary}\label{4:4.10}
Assume that $L$ is a connected diagram of a link with $n\geq 1$
crossings. 
Let $G$ be a planar graph associated to $L$ via some 2-color (checkerboard)
coloring of the plane containing $L$
and let $G^*$ be the graph associated to $L$ via 
the opposite coloring of the plane. (The graph $G^*$ is dual to 
$G$ and colors of the edges are reversed.)
Then
$$\varphi^+_L (t) =
(-1)^n(1-t)^{-2}C(\overline{G}_w,1-t)C((G^\star)_w;1-t)$$

$$\varphi^-_L (t) = 
(-1)^n(1-t)^{-2}C(\overline{G}_b,1-t)C(\overline{(C^\star)_b},1-t)$$
\end{corollary}

Proof. We apply a formula relating chromatic polynomial to 
Tutte polynomial (see Section 1) and we use relations 
$G_b = \overline{(G^\star)_w}$ , $G_w = \overline{(G^\star)_b}$
and  $\chi(G;x,y) = \chi(G^\star;y,x)$.

Now we will describe some applications of adequate and semi-adequate
diagrams of links.

%\begin{corollary}\label{V.5.11
\begin{corollary}\label{4:4.11}
Let $L$ be a connected diagram of a link
with at least one crossing. Then:

\begin{enumerate}

\item $\varphi^+_L\neq 0$ if and only if $L$ is $+$ adequate,\\

 $\varphi^-_L\neq 0$ if and only if $L$ is $-$ adequate.

\item If $L$ has a non-nugatory crossing 
and it is semi-adequate then the degree of $x$ 
in  $\Lambda (a,x)$ is positive and thus $L$ represents a non-trivial link.
\end{enumerate}
\end{corollary}

Proof.
\begin{enumerate}

\item Follows immediately because of Lemma \ref{4:4.7} and Corollary 
5.2.

\item Follows from (1) and Corollary 5.2(v).
\end{enumerate}

M.~Thistlethwaite checked that all diagrams of minimal crossing-number 
up to 11 are semi-adequate and among knots which have diagrams with at most
12 crossings only few possibly do not admit semi-adequate diagrams. 
For example a 12 crossing knot of Fig. 5.9 does not admit 
a semi-adequate diagram of 12 crossings but it is an open problem 
whether it admits a semi-adequate diagram with more 
than 12 crossings \cite{This-5}. 
\\ \ \\
\centerline{{\psfig{figure=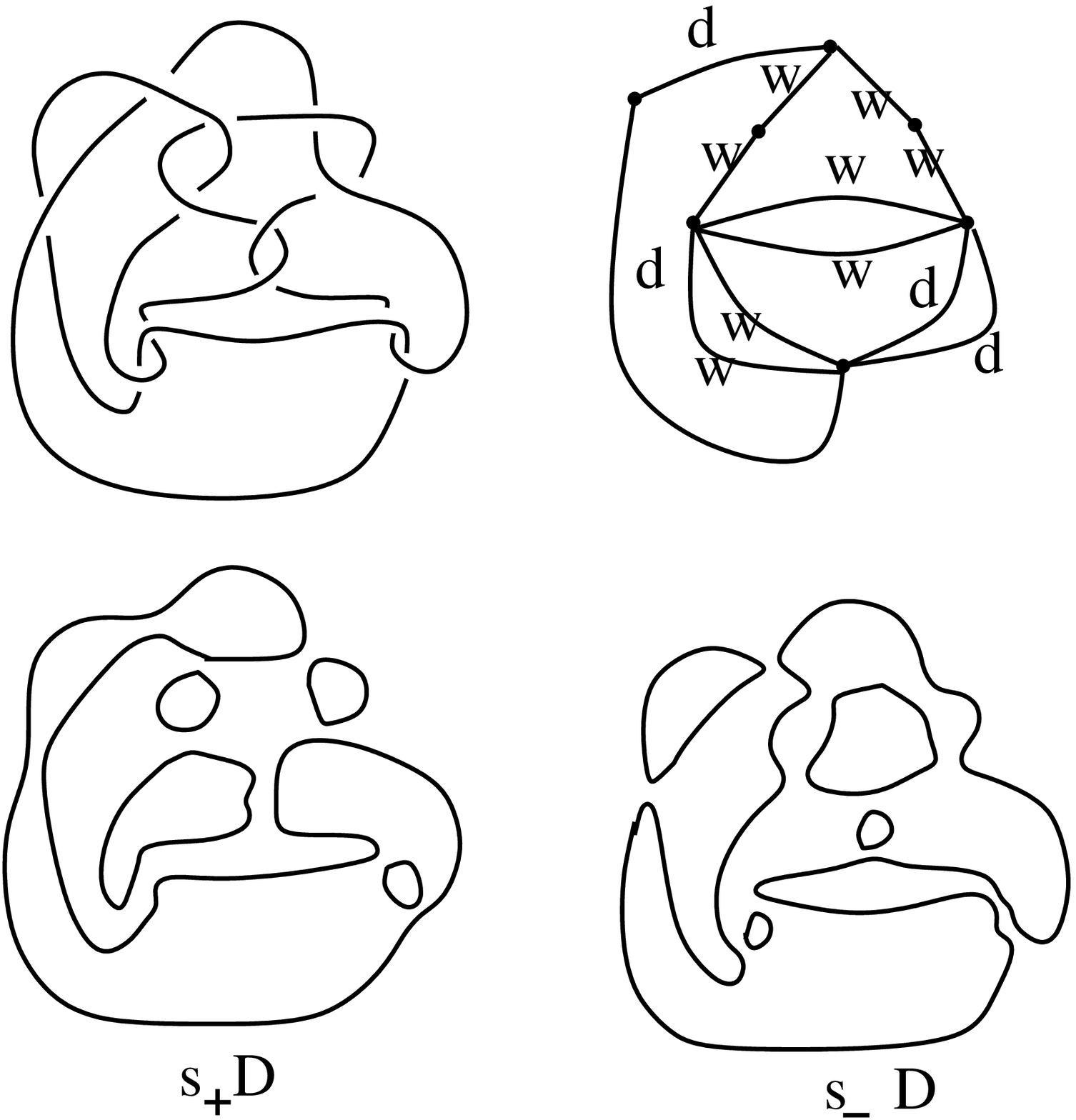,height=7.0cm}}}
\begin{center}
Fig.~5.9
\end{center}

\begin{corollary}\label{4:4.12}
A semi-adequate diagram is not regularly isotopic to a diagram with 
a smaller number of crossings. Furthermore it is not isotopic (equivalent) 
to a diagram with a smaller number of crossings and the same Tait number. 

\end{corollary}

Proof. If $L$ is a connected semi-adequate diagram with $n>1$
crossings then because of Corollary \ref{4:4.11} (i) the polynomial
$\Lambda_L$ has a non-zero exterior coefficient 
%($\varphi^+$ or $\varphi^-$ $\neq 0$).
that is, it contains a non-zero term $u_{r,s}\cdot a^r z^s$, 
where either $r+s=n$ or $-r+s = n$. Thus $L$ is not regularly isotopic
to a diagram with fewer crossing (we recall that $\Lambda_L$
is a regular isotopy invariant). If the diagram $L$ has components
$\row{L}{c}$ then
${\Lambda}_L = 
(\frac{a+a^{-1}}{z} -1)^{c-1}\cdot 
\Lambda_{L_1}\cdot\dots\cdot\Lambda_{L_c}$
and the corollary follows. The second part of the corollary follows from 
the fact that if two diagrams are isotopic and have the same Tait number 
then they are related by second, third and balanced Reidemeister moves (in 
which pair of kinks of opposite signs is created or deleted). All these 
moves preserve ${\Lambda}_L$.

The conclusion remains true if we merely assume that any component 
of the diagram is semi-adequate.

\begin{corollary}\label{4:4.13}
A connected semi-adequate diagram cannot describe a split link,
that is, a link that can be 
separated in $S^3$ by a sphere $S^2$.
\end{corollary}

Proof. Let us assume that a link $L$ can be decomposed into two
sublinks $L_1$, $L_2$ and suppose that $D$ is a connected
diagram of $L$. Certainly, the sum of signs of crossings between
$D_1$ and $D_2$ (associated to $L_1$ and $L_2$, respectively)
is equal to 0. Therefore the link $L$ can be presented by a diagram
$D'$ consisting of $D_1$ and $D_2$, which lie on different
levels (one above another). The diagram $D$ is isotopic to 
$D'$ and has the same Tait number ($(Tait (D) = Tait(D'))$) hence
$\Lambda (D) = \Lambda (D')$.  Now diagrams $D_1$ and $D_2$ can be 
moved apart in $D'$, by regular isotopy, to decrease the number of crossings.
But then the diagram $D'$ (and also $D$) is not semi-adequate.

\begin{corollary}\label{4:4.14}
Let $L_1$ and $L_2$ be two diagrams with $n$ crossings which represent
the same link. 
If $L_1$ is $+$-adequate then $Tait(L_1)\geq Tait (L_2)$,
and if $L_1$ is $-$-adequate then $Tait (L_1)\leq Tait (L_2)$. 
If $L_1$ is adequate then also $L_2$ is adequate and  
$Tait (L_1) = Tait (L_2)$.
\end{corollary}

Proof. Since $L_1$ and $L_2$ are isotopic then 
$$a^{-Tait (L_1)}\Lambda_{L_1} = a^{-Tait (L_2)}\Lambda_{L_2}$$
and thus
$$\Lambda_{L_1} = a^{Tait L_1 -Tait (L_2)}\Lambda_{L_2}.$$

If $L_1$ is $+$ adequate then $\deg_{a,z}\Lambda_{L_1} =n$,
and since
$\deg_{a,z}\Lambda_{L_2}\leq n$ it follows that $Tait (L_1) - Tait 
(L_2)\geq 0$.                                  
%(???pojawilo sie oznaczenie na stopien, czy bylo zdefiniowane?
%jesli tak to mozna bylo skrocic pewne rozwazania juz wczesniej???)

Similarly, if $L_1$ is $-$ adequate then $Tait (L_1) \leq Tait
(L_2)$. Therefore, if $L_1$ is adequate then $Tait (L_1) = Tait 
(L_2)$ and consequently $\Lambda_{L_1} = \Lambda_{L_2}$, and because 
of Corollary \ref{4:4.11} $L_2$ is adequate as well.

\begin{corollary}\label{4:4.15}
If a link $L$ has an adequate diagram with $n$ crossings then
the link $L$ does not admit a diagram with fewer crossings.
\end{corollary}

We cannot claim however that semi-adequate diagram without a nugatory crossing 
cannot be isotoped to a diagram with smaller number of crossings.
Figure 5.10 presents a $+$-adequate diagram of the right handed trefoil 
knot with 4 crossings but without a nugatory crossing.
\\ \ \\
\centerline{{\psfig{figure=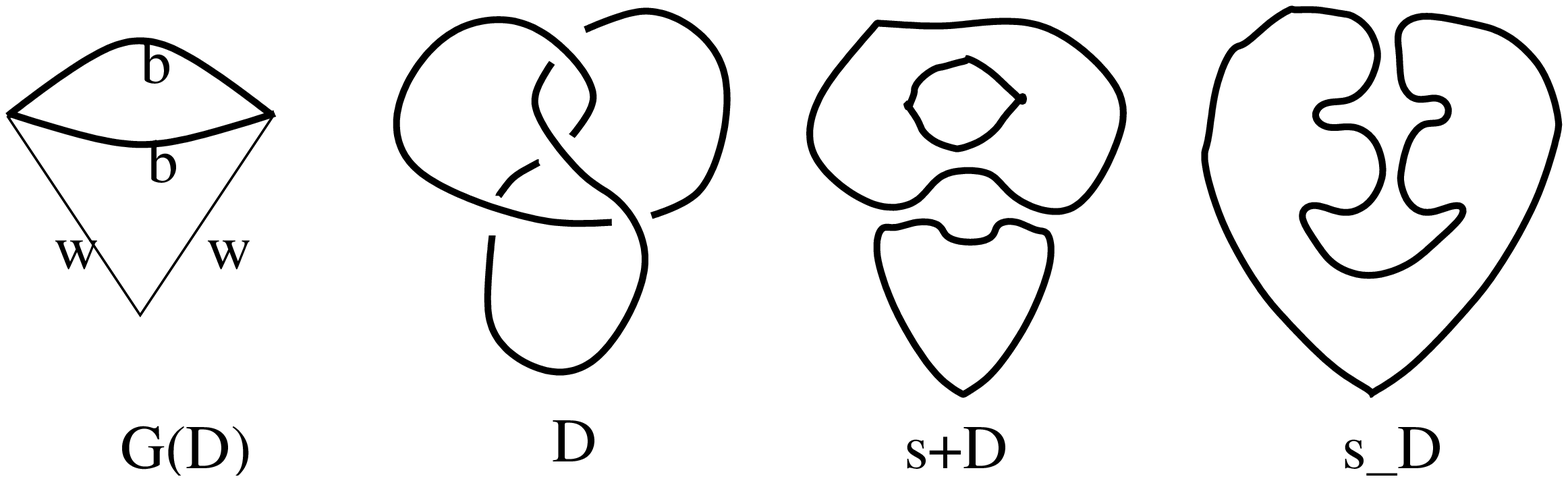,height=3.9cm}}}
\begin{center}
Fig.~5.10
\end{center}

\begin{corollary}\label{4:4.16}
For any positive integer $n$ there exists a prime diagram of a 
link $L$ 
which has the minimal number of crossing among diagrams representing
$L$ and which has a bridge longer than $n$.
\end{corollary}

\begin{proof}
 Consider the infinite family  of 2-color graphs
$G_1,G_2,\ldots$; of which the first three members are illustrated 
in Fig. 5.11 ($G_n$ consists of two nested families of circles, $n+1$ 
on the left and $n$ on the right). 
The black subgraph 
$(G_i)_b$ is to consists of the left-hand family of nested circles. 
Because of Lemma \ref{4:4.7} the graph $G_i$
is associated to an adequate diagram of a link
and therefore the diagram has the minimal number of crossings.

%\vspace*{2in}\centerline{\Psfig{figure=Rys.4.9}}
\centerline{{\psfig{figure=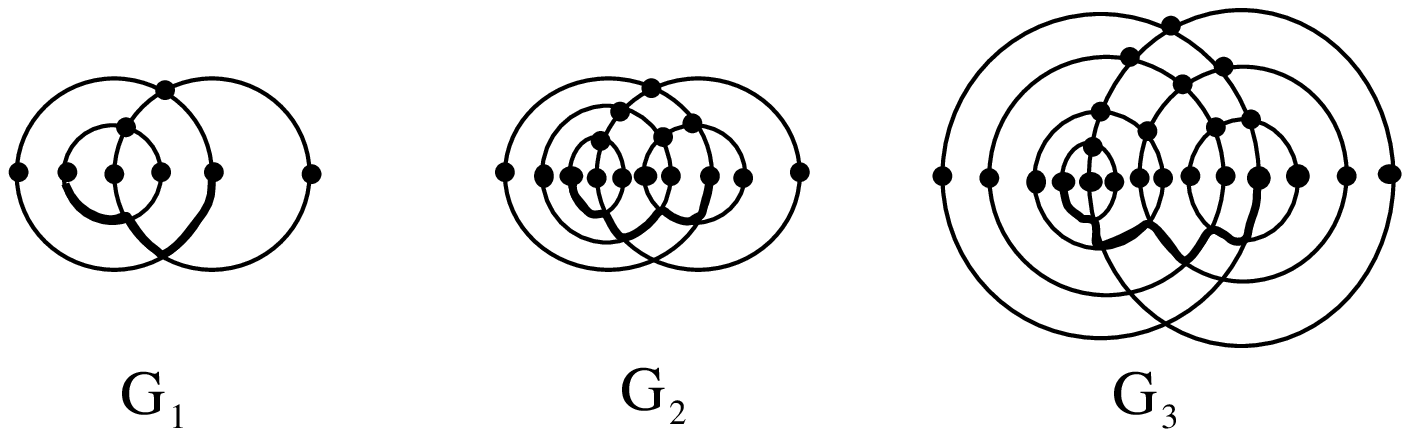,height=3.9cm}}}  
\begin{center}
Fig.~5.11
\end{center}

The thick edges  of graphs $G_1$, $G_2$ and $G_3$ 
represent crossings of a bridge of length $2i+1$ in the associated
diagram of a link. Figure 5.12 presents links related to $G_1$ and $G_2$.
\\ \ \\
%\vspace*{3in}\centerline{\Psfig{figure=Rys.4.10}}
\centerline{{\psfig{figure=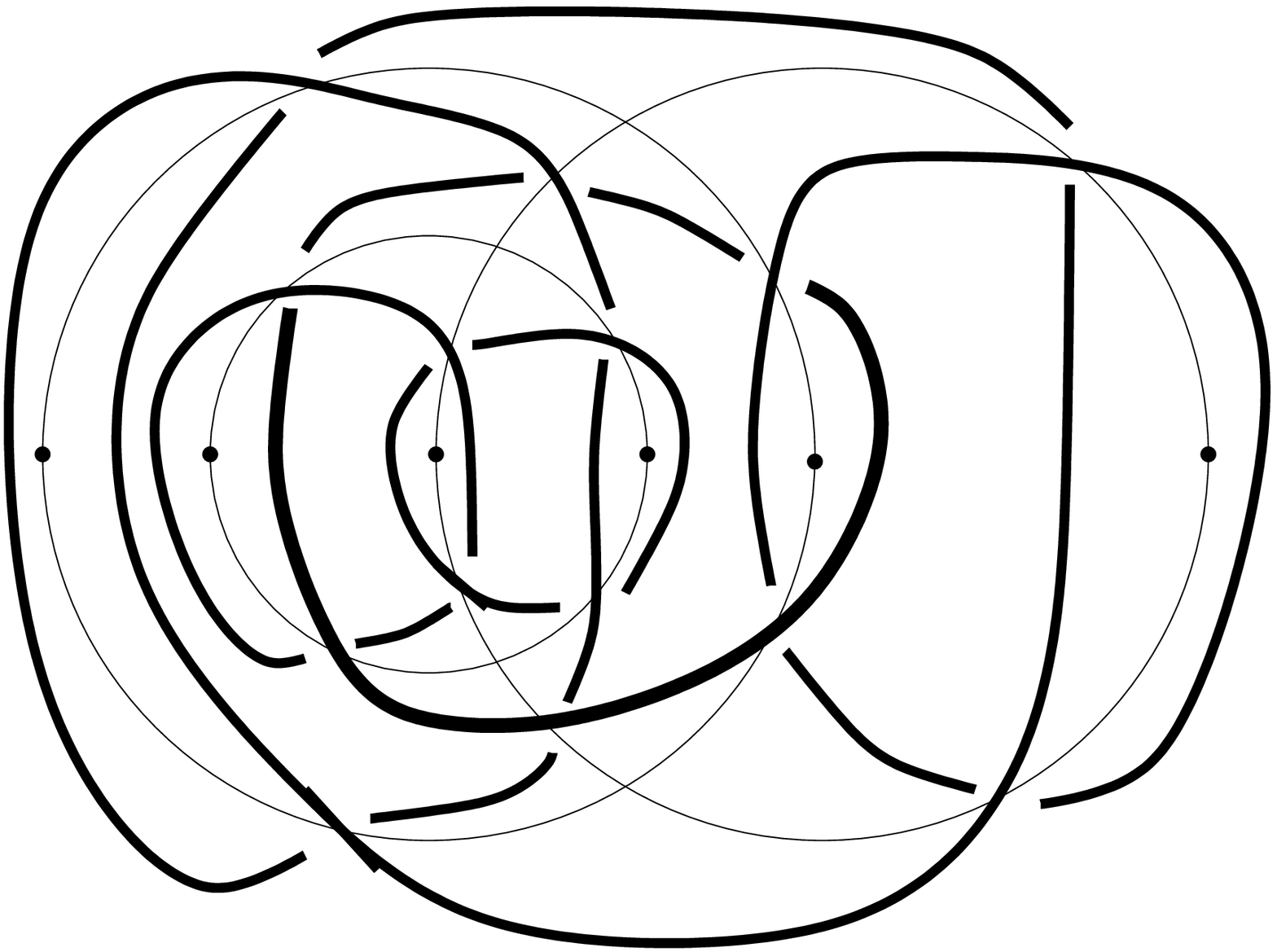,height=5.0cm}}\ \ \ 
{\psfig{figure=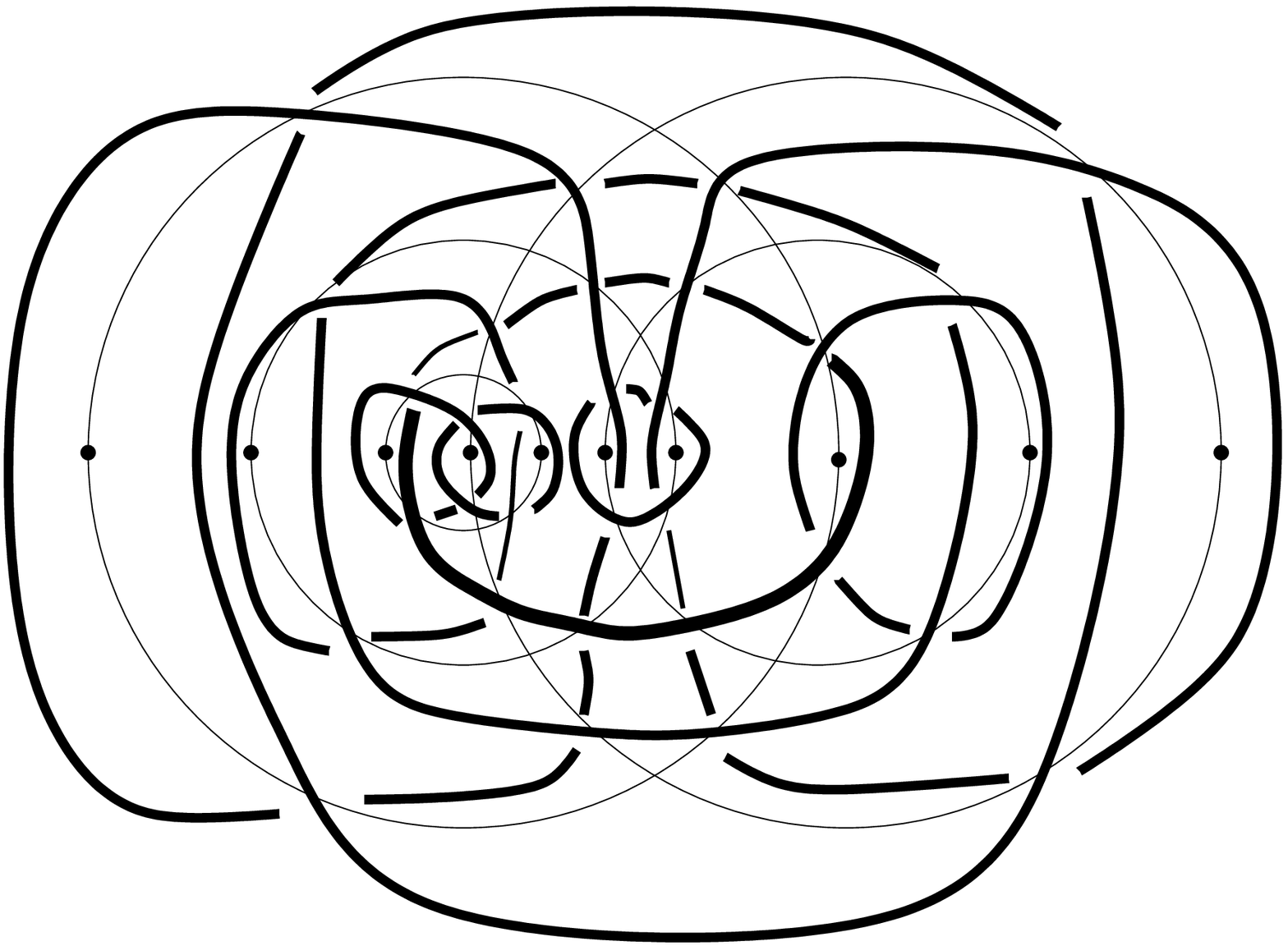,height=5.0cm}}}
%\centerline{{\psfig{figure=Rys4-4-10.eps,height=5.0cm}}}  
\begin{center}
Fig.~5.12 
\end{center}

If we want to get knots instead of links in the lemma then it can be done
by adding to graphs $G_i$ some appropriately chosen midpoints of edges.  
\end{proof}

\begin{corollary}\label{4:4.17}
Computation of a the Kauffman polynomial of a link is
$NP$-hard.
If the well known conjecture that  $NP\neq P$ holds then
the Kauffman polynomial cannot be computed in polynomial time
with respect to the number of crossings.
\end{corollary}

Proof. 
Computation of a chromatic polynomial of a planar
graph is $NP$-hard, 
see (\cite{G-J}). Because of Theorem \ref{4:4.1} and Corollary \ref{4:4.11}
it follows that the computation of exterior coefficients of 
the Kauffman polynomial of an alternating diagram can be reduced to 
computation of chromatic polynomial of the 
associated (one-color) planar graph.

Therefore, the computation of the Kauffman polynomial of an alternating 
link is $NP$-hard.

\begin{example}[Thistlethwaite]\label{V.5.18}
Consider two 2-colored graphs $G_1,G_2$ of Fig. 5.13. 
Their associated diagrams $D_1,D_2$, are
the famous Perko pair which for many years was though to represent 
different knots. Both diagrams  have 10 crossings and they have  
different Tait numbers.
In 1974 K.A.Perko (who did master degree with Fox at Princeton and later 
became a lawyer in New York), 
 noticed that they represent the same knot. 
%$D_1$ is $-$-adequate and $D_2$
The outermost polynomials  $\varphi^+_{D_i}(t)$ and 
$\varphi^-_{D_i}(t)$ are as follows:
$$\varphi^+_{D_1}(t) = 0 = \varphi^-_{D_2}(t),\ \ 
\varphi^-_{D_1}(t) = t^2(t^2+t^3),\ \varphi^+_{D_2}(t)= t^4.$$ 
The full Kauffman polynomials of diagrams are
$$\Lambda_{D_1}(a,z)= a^{-2}\Lambda_{D_2}(a,z)= (a^{-4}+a^2)z^6 +
(a^{-5} + a^{-3})z^5 + (-4a^{-4}+a^{-2}-1 - 6a^2)z^4 +$$ 
$$(-4a^{-5} - 3a^{-3} - a)z^3 + (3a^{-4}-a^{-2}+ 3+9a^2)z^2 + 
(3a^{-5} + a^{-3} +2a)z + a^{-2}-1 - 3a^2.$$
By Corollary 5.14 the Perko knot has no adequate diagram, but 
it has a $+$-adequate diagram, $D_2$ and $-$-adequate diagram $D_1$.
\end{example} 
\ \\ 
\centerline{{\psfig{figure=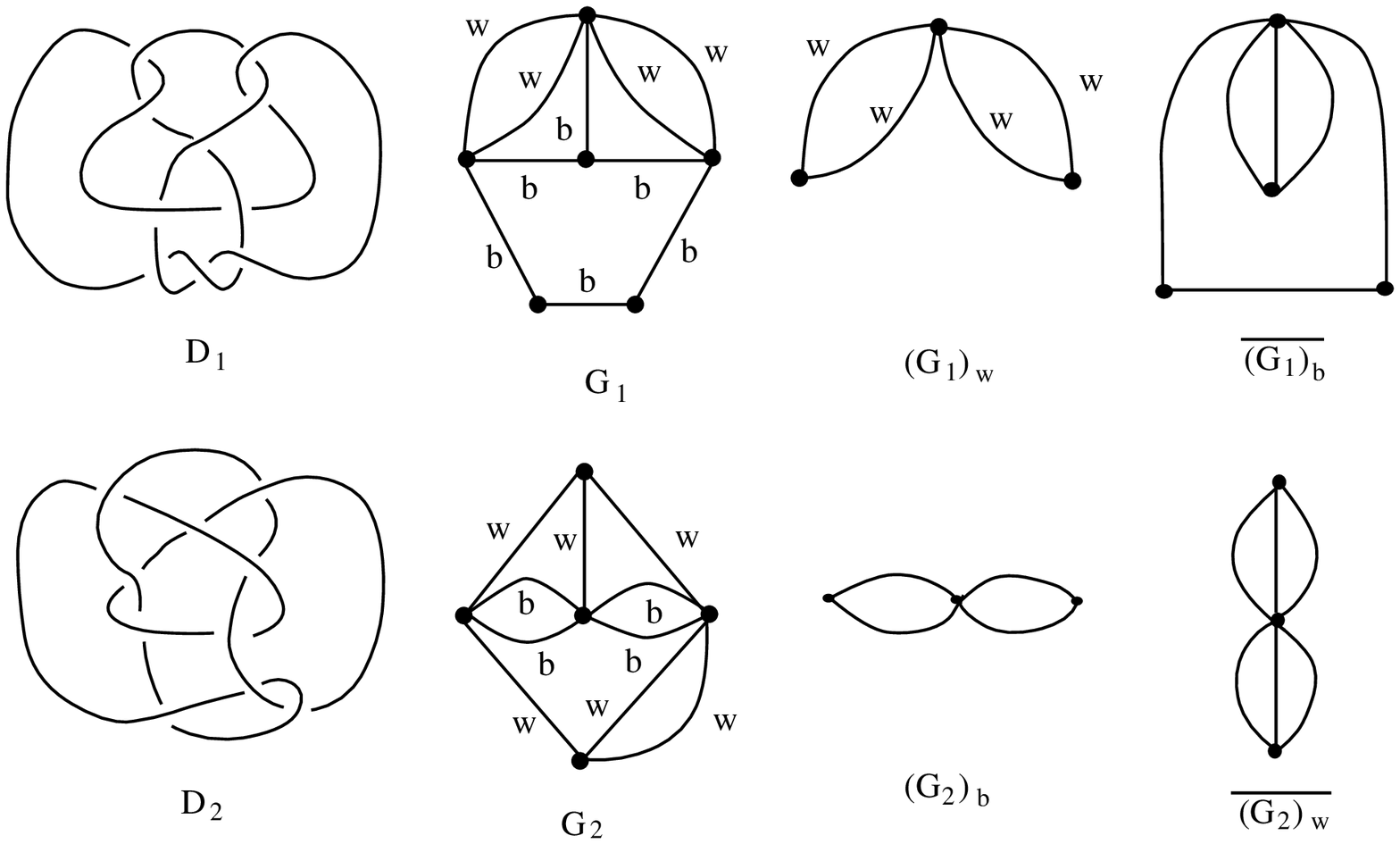,height=6.8cm}}}
\begin{center}
Fig. 5.13
\end{center}

%\begin{exercise}\label{V.5.19}
\begin{exercise}\label{4:4.18}
Let us consider a substitution $x= a+a^{-1}$ in Kaufman polynomial
$\Lambda_L (a,x)$ of an non-oriented link. Let us note that if 
$T_n$ ($n>1$) is a diagram of a trivial link with  $n$ components then 
$\Lambda_{T_n} (a,a+a^{-1}) = 0$. This way $\Lambda_L
(a,a+a^{-1})$ is similar to Alexander-Conway polynomial. 
We will denote this polynomial by $\Lambda_L (a)$, and 
the related polynomial for oriented links will be denoted by 
$F_L (a)$, that is $F_L (a) = a^{-Tait (L)}\Lambda_L (a)$.

Prove that
\begin{enumerate}
\item[(1)] If $K$is a knot then $F_K (a) = 1 + (a+a^{-1})(G_K (a))$
for some Laurent polynomial $G_K(a)$. If $L$
is a link consisting of $\mu (L)>1$ 
components then $F_L(a) =(a+a^{-1})^{\mu (L)-1}(G_L (a))$ 
where $G_L(a)$ is a Laurent polynomial.

\item[(2)] If $L$ is a connected $+$-adequate (respectively, $-$-adequate)
diagram of a link then $max\ deg \Lambda_L (a) = n(L)$ (respectively,
$min\ deg \Lambda_L (a) = -n(L)$).

\item[(3)] For any connected diagram $L$
$$\Lambda_L (a) = \sum^{n(k)}_{-n(L)} u_i a^i$$
and moreover
$u_{n(L)} \neq 0$ if and only if $L$ is $+$ adequate.

\item[(4)] We have $u_{n(L)} = \varphi^+_L (1) =
(-1)^{n(L)}\cdot C'(\overline{G}_w,0)\cdot C'(\overline{G}^\star_w,0)$
where $C'$ denotes the derivative of chromatic polynomial.
\end{enumerate}
\end{exercise}

The polynomials $F_L (a)$ and $\Lambda_L (a)$ should be easier to analyze
than the general the Kauffman polynomial. They seem to be particularly useful
to examine periodicity of links (see Chapter VII).
This may be an interesting research problem.
%research problem related to the fact that some Vassiliev invariants
%related to Kauffman polynomial are products of smaller Vassiliev 
%invariants as outermost polynomials are products!!

%%%%%%%%%%%%%%%%%%%%%%%%%%%%%%%%%%%%%%%%%%%%%%%%%%%%%%%%%%%%%%%%

%\section{Coefficients of Homflypt polynomial}
\section{Coefficients of Jones-Conway polynomial}
\markboth{\hfil{\sc Graphs and links}\hfil}
{\hfil{\sc Coefficients of Jones-Conway polynomial}\hfil}

We start this section with a general theorem characterizing coefficients
of Jones-Conway polynomial $P_L (a,z)$. This result was essentially proved
by Morton \cite{Mo-2,Mo-3} and independently by Franks and Williams 
\cite{F-W}. It is used to give a good approximation of the 
braid index of a link. In our exposition we rely on \cite{Mo-3}. 

Let $L$ be a diagram of an oriented link. By $n^+(L)$ and $n^-(L)$ 
we denote the number of positive and, respectively, negative
crossings of $L$.  Moreover, in this section $n(L)$ 
denotes the number of crossings and $\tilde{n}(L)$ 
denotes the algebraic number of crossings (we often use the notation 
$Tait (L)$ for this number).
Thus we have:
$$n(L) = n^+(L) + n^-(L)$$
$$\tilde{n}(L) = n^+(L) - n^-(L).$$

If we smooth all crossings of $L$ (respecting the orientation)
then we obtain a family of simple closed curves called Seifert circles 
of the diagram $L$ (c.f.~Fig. 6.1; compare Chapter IV). 
The number of these circles is denoted by $s(L)$. Furthermore, 
define, after P.Cromwell, the {\it Seifert graph}, $G_{\vec{s}}(L)$ of an
oriented diagram $L$ as a signed graph whose vertices are in bijection with
Seifert circles of $L$ and signed edges correspond to crossings of $L$.
If a positive (resp. negative) crossing connects Seifert circles then
the positive (resp. negative) edge connects related
vertices of $G_{\vec{s}}(L)$.
\\
\centerline{\psfig{figure=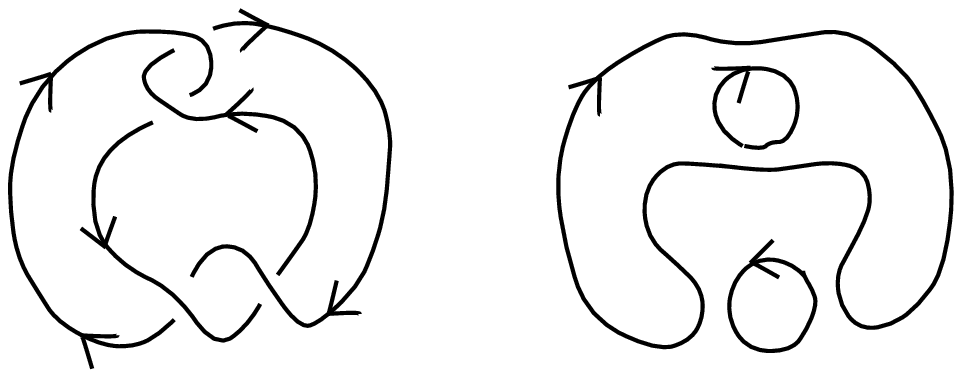,height=3.0cm}\ \ \ \ 
\ \ \ \ \ \  
\psfig{figure=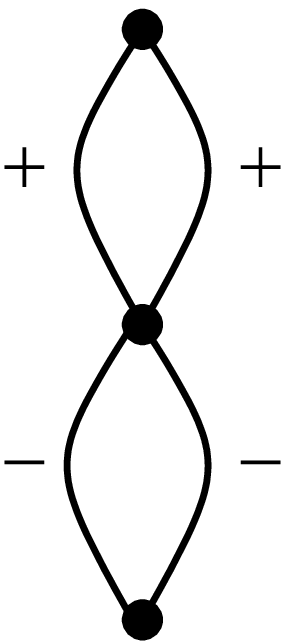,height=3.0cm}}  
\begin{center}
Fig. 6.1
\end{center}

Recall that the version of the Jones polynomial we use 
satisfies the skein relation $aP_{L_+} + a^{-1}P_{L_-}=zP_{L_0}$.
We have 
$$P_L(a,z) = \sum^E_{k=e} a_k(z)a^k,$$ 
where $a_e,\ a_E\neq0$, when 
written as a polynomial in the variable $a$. 
We define the $a$-span of $P_L$ as $\beta_a(L) = E-e$.
Similarly, we can write 
$$P_L(a,z) = \sum^M_{k=m} b_k(a)z^k,$$ where $b_m,\ b_M\neq 0$ 
and we define the $z$-span of $P_L$ as $\beta_z (L) = M-m$.

\begin{lemma}[\cite{L-M-1}]\label{4:5.1}
In the above notation $m=1-\mu (L)$ where $\mu (L)$ is the number
of components of the link  $L$. 
If $L$ is a knot then $a_0(i) = 1$, where $i=\sqrt{-1}$. 
\end{lemma}

We begin by proving the second part of the lemma.
First, recall (see  Chapter III) 
that the substitution $a=i$ and $z = i(\sqrt{t} - \frac{1}{\sqrt{t}})$
 in the Jones-Conway polynomial $P_L(a,z)$
leads to the Alexander polynomial $\bigtriangleup_L (t)$ as 
normalized by Conway. Therefore, we are in the situation of proving
 a known result on Alexander polynomial.
Let $P_L(z) = P_L (i,z) = \sum^T_{k=t} c_kz^k$ where 
$c_t,\ c_T\neq 0$. We have, therefore, $P_{L_+} - P_{L_-}=-izP_{L_0}$.
With this notation we first prove the following.
\begin{formulla}\label{V.6.2}
$t\geq\mu (L) -1$ and if
$L$ is a diagram of a knot then $t=0$ and $c_0 =1$.
\end{formulla}
\begin{proof}
We consider diagrams with chosen (ordered) base points and 
use induction with respect to
the  lexicographically ordered pair $X(L)=$ (number of crossings, number of 
``bad'' crossings) for the diagram $L$. 
Recall, that for a given choice of base points, a crossing is 
called ``bad'' if it has to be changed in order to make the diagram 
descending, see Chapter III. 

If the diagram $L$ represents a trivial link (e.g.~if $L$ is a descending
diagram) then
$$P_L(z) = \left\{
\begin{array}{lll}
1&\mbox{if}&\mu (L) =1\\
0&\mbox{if}&\mu (L) > 1\\
\end{array}
\right.
$$
And therefore 6.2 is true. In particular, it holds for diagrams 
with $n(L) \leq 1$ crossings.
Suppose now that 6.2 is true for diagrams with fewer than
$n(L)$ crossings (where $n(L)>1$) and for diagrams 
%(??po polsku bylo ciagle sploty, ale chodzi-jak rozumiem o ich diagramy ????) 
with $n(L)$ crossings but with fewer than $b(L)$ bad crossings
(where $b(L)\geq 1$).

Let $p$ be the first bad crossing of $L$. Suppose that  $L'$ is a diagram 
obtained from
$L$ by changing the crossing $p$ and  $L_0$ is a diagram  obtained 
by smoothing of $p$; both $L'$ and $L_0$ satisfy the inductive assumption
so 6.2 holds for them.
Moreover  $P_L(z) = P_{L'} (z) - sgn(p) iz P_{L_{0}} (z)$. 
Now 6.2 follows immediately if we note that 
$\mu (L) = \mu (L') = \mu (L_{0})\mp 1$ and 
$\mu(L_{0}) = 2$ when $L$ is a knot.
\end{proof}
The first part of Lemma 6.1 follows from 6.2 as $a_0(i)=c_0=1$.
The next step is to show that $m\geq 1-\mu (L)$. We use again an 
induction on complexity $X(L)$ noting that it holds for a trivial 
link $T_{\mu}$ as $P_{T_{\mu}}=(\frac{a+a^{-1}}{z})^{\mu -1}$, and 
applying skein relation. To show that $m=1-\mu (L)$, we establish 
the formula from \cite{L-M-1} using 6.2 and inequality $m\geq 1-\mu (L)$.
\begin{formulla}\label{V.6.3}
$$b_{1-\mu (L)}(a) = (-a^2)^{-\lk (L)} (a+a^{-1})^{\mu (L)
-1}\cdot\prod^{\mu(L)}_{i=1} b_0^{L_i} (a)$$
\end{formulla}
where $\row{L}{\mu(L)}$ are components of $L$ and $\lk (L)$
is the global linking number of $L$.
Here $b_0^{L_i} (A)$ denotes free coefficient of $P_{L_i} (a,z)$
which, because of 6.2, is non-zero ($b_0^{L_j}(i)=1$).

\begin{proof}
We use an induction with respect to the
number of crossings which have to be changed in the diagram of $L$ in order
to make 
$L_1$ lying over $L_2\cup...\cup L_{\mu(L)}$, $L_2$  over 
$L_3\cup...\cup L_{\mu(L)}$,..., $L_{\mu(L)-1}$ over 
$L_{\mu(L)}$ in $R^3=R^2 \times R$. 
%$L_{\mu(L)}$ (???po pierwsze: chyba sa tu misprinty,
If we do not have to change any crossing then 
$\lk (L) =0$ and $L$ is a split sum of 
$\row{L}{\mu(L)}$, therefore 
$$P_L(a,z) =
(\frac{a+a^{-1}}{z})^{\mu(L) -1}\prod^{\mu(L)}_{i=1} P_{L_i} (a,z)$$
and formula 6.3 holds for $L$. 

Finally, if $p$ is a crossing between two different components of $L$
then $b_{1-\mu (L_+)}(a) = -a^{-2}b_{1-\mu(L_-)}(a)$ and then the
inductive step follows (note that in this case $P_{L_0}(a,z)$ do not 
contribute to the formula as $1-\mu(L_0)= 1-\mu(L)+1$).
\end{proof}

\begin{theorem}\label{V.6.4} 
For any diagram $L$ we have
$$-\tilde{n}(L) - (s(L) -1)\leq e\leq E\leq -\tilde{n}(L)+(s(L)-1).$$
\end{theorem}

\begin{theorem}\label{4:5.5}
For any diagram $L$ the highest degree in $z$ in Jones-Conway polynomial 
satisfies:
$$M\leq n(L) - (s(L)-1).$$
\end{theorem}

\begin{corollary}\label{4:5.6}
For any diagram $L$ we have
$$s(L)\geq\frac{E-e}{2}+1 = \frac{1}{2}\beta_a (L) +1.$$
\end{corollary}

\begin{corollary}\label{4:5.7}
For any diagram $L$ of an amphicheiral link we have  
$$|\tilde{n}(L)|<s(L).$$
\end{corollary}

Proof of Corollary 6.7.\ 
If a link is amphicheiral then $e = -E$, hence $e\leq
0\leq E$, and by Theorem 6.4
$$\tilde{n}(L) \leq s(L) -1\ \ \mbox{ and }\ \ \tilde{n}(L) \geq -(s(L) -1),$$
and therefore 
$$|\tilde{n}(L)|< s(L).$$

Proof of Theorem 6.4

It is enough to prove the inequality $-\tilde{n}(L) - (s(L) -1)\leq e$.
Indeed, for the mirror image $\overline{L}$ we have
$$\tilde{n}(\overline{L}) = -\tilde{n}(L),\ \ \ 
s(\overline{L}) = s(L), \ \ \ 
e(\overline{L}) = -E(L)$$ and therefore the inequality
$-\tilde{n}(\overline{L}) - (s(\overline{L}) - 1) \leq e(\overline{L})$
implies $-\tilde{n}(L) +s(L) -1\geq E(L)$.

Let us consider the function $\varphi (L) = \tilde{n}(L)  +(s(L) -1)$.

To prove Theorem 6.4 we have to show that 
$a^{\varphi(L)}P_L(a,z)$ is a polynomial in the variable $a$ 
(i.e.~$a$ does not occur with negative exponent).
The Seifert circles of  $L_+$, $L_-$ and $L_0$ are the same
and therefore
$$\varphi(L_+) = \varphi(L_0)+1 = \varphi(L_-)+2,$$ 
which implies
$$a^{\varphi(L_+)}P_{L_+} (a,z) + a^{\varphi(L_-)}P_{L_-} (a,z) =
za^{\varphi(L_0)}\cdot P_{L_0} (a,z).$$
So if $a^{\varphi(L)}\cdot P_{L} (a,z)$ is a polynomial in $a$
for two of the three diagrams 
$L_+$, $L_-$ and $L_0$ then the same property holds for the third one. 
Now the standard induction with respect to the number of all crossings
and the number of bad crossings of the diagram $L$ reduces our problem to the 
case of descending diagrams (c.f.~proof of Theorem III.??1.2). 
%\ref{2:1.2}). 
Since for a descending diagram with $\mu  (L)$ components
we have $P_L(a,z) = (\frac{a+a^{-1}}{z})^{\mu (L)-1}$ it follows that
we have to prove that $\varphi (L)\geq\mu (L)-1$.

The proof is by induction with respect to the number
of crossings in $L$ (we assume that the theorem holds for all 
diagrams with smaller number of crossings).

If $L$ has no crossing at all, then $\tilde{n}(L) = 0$ and $s(L)=\mu
(L)$ hence $\varphi (L) = \mu(L) -1$. Suppose that $n(L)>0$ 
and let us assume that the claim is true for all diagrams 
with smaller number of crossings.
Let  $b = (\row{b}{\mu(L)})$ be a sequence of base points 
on $L = L_1\cup\cdots\cup L_{\mu(L)}$ such that $L$ is descending 
with respect to $b$.
We have to consider the following two cases.

\begin{enumerate}
\item[(1)] There exists a self-intersecting component $L_i$.
Let $p$ be the first self-intersection point of $L_i$ (that is 
we start drawing our diagram from $b_i$ and we do not have intersection 
till we reach $p$).
Let $L_0$, as usual, denote the diagram obtained from $L$ by smoothing at $p$.
Then  $s(L_0) = s(L)$ and $\varphi(L_0) = \varphi(L)\mp 1$, 
depending on the sign of the crossing $p$.
On the other hand $L_0$ represents a trivial link 
(even a descending diagram) 
with fewer number of crossings than $L$.
Therefore, by inductive assumption,
$\varphi(L_0)\geq\mu(L_0)-1=\mu(L)$ hence $\varphi(L)\geq\mu(L)-1$.

\item[(2)] All crossings of $L$ are
crossings between different components of $L$ (each $L_i$ is a simple closed 
curve).

A change of height of components preserves 
$\varphi (L)$ and $\mu(L)$. Therefore we may assume
that crossings occur between components $L_1$ and $L_2$.
Let us consider a positive crossing $p$ of these two components
(there exists a positive
crossing because $\lk(L_1,L_2) = 0$).
Now, after smoothing $p$, we obtain a descending diagram  $L_0$ with
$\mu(L_0) = \mu(L) -1$ components. To be sure that $L_0$ is descending 
we can choose the base point $b'$ of the new component of $L_0$ to be
on the ``old" $L_1$ just after the crossing $p$ (smoothing and the 
choice of $b$ is illustrated here  
\parbox{2.3cm}{\psfig{figure=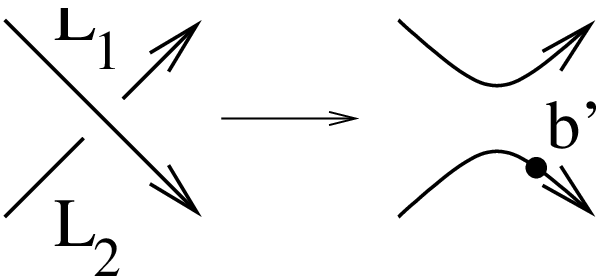,height=0.9cm}}).
 
Clearly $\tilde{n}(L_0) = \tilde{n}(L) -1$.
Because of the inductive assumption $\varphi (L_0)\geq\mu(L_0)-1$
and therefore $\varphi(L) =\tilde{n}(L)+s(L)-1=\tilde{n}(L_0)+1+s(L_0)-1=\varphi(L_0)+1\geq\mu(L_0)
= \mu(L)-1$. This concludes the proof of Theorem 6.4.
\end{enumerate}

The braid index, $b(L)$,  of an oriented link $L$ is 
the minimal number of strings needed so that $L$ is represented by a closure 
of a braid. It is always well defined, as Alexander proved that every 
link can be deformed (by ambient isotopy) to a braid 
form (in an implicit form it was 
already observed by Brunn in 1897). We discuss Alexander theorem in 
APPENDIX.

\begin{corollary}\label{4:5.8}
Let $L$ be an oriented link. 
Then
$$b(L)\geq\frac{1}{2}(E-e)+1$$
\end{corollary}

Proof. The diagram of a closed braid with $k$ strings has exactly $k$
Seifert circles. Therefore \ref{4:5.8} follows
by Corollary \ref{4:5.6}.

Corollary 6.8 turns out to be an exceptionally
efficient tool to determine $b(L)$. 
Jones checked that it was sufficient to determine the braid 
index of 265 knots among 270 prime knots which have diagrams 
with at most 10 crossings \cite{Jo-2}.
The exceptional knots are $9_{42}, 9_{49}, 10_{132},
10_{150}, 10_{156}$ (according to Rolfsen's notation \cite{Ro}). 
	Let us discuss briefly the case of $9_{42}$ (the knot turned out 
			to be exceptional already having the same Jones-Conway and Kauffman 
			polynomials as its mirror image $\bar 9_{42}$ but not being 
			amphicheiral). 
	It can be represented as a braid with four strands
	$\sigma^3_2 \sigma_3\sigma^{-1}_1 \sigma_2 \sigma^{-2}_{3} 
\sigma^{-1}_1 \sigma_2 \sigma^{-1}_3$
	(Fig.~6.2) and its Jones-Conway polynomial is equal to
	$$P_{9_{42}} (a,z) = a^{-2}(-2+z^2)+(-3+4z^2-z^4)+a^2(-2+z^2)$$
	and therefore by \ref{4:5.8} we have $b(9_{42})\geq\frac{1}{2}(2-(-2))+1 = 3$.
	\\
	\ \\
%\centerline{{\psfig{figure=4-5-2.eps,height=3.1cm}}}
\centerline{{\psfig{figure=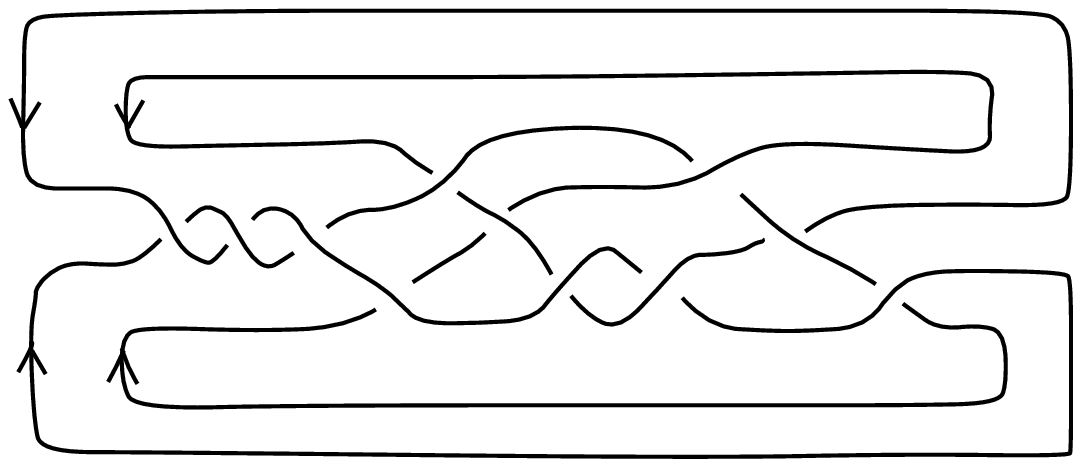,height=4.1cm}}}
\begin{center}
Fig.~6.2
	\end{center}
%Rotate Figure 6.2)
Therefore Corollary \ref{4:5.8} cannot exclude the possibility 
that the braid index of $9_{42}$ is equal 3. 
In order to prove that actually $b(9_{42}) = 4$,  Morton and 
Short \cite{M-S-1,M-S-2} 
applied the following argument
which used Jones-Conway polynomial together with \ref{4:5.8}.
If $9_{42}$ had a braid representation with three strands
then the cable satellite of type $(2,0)$ around the braid 
(compare Chapter VI)
would have a braid representation with 6 strands.
On the other hand,
Jones-Conway polynomial of this satellite is equal to
$$a^{-5}(4z^{-1}+30z+76z^3+85z^5+45z^7+11z^9+z^{11})+
\cdots+a^7(-7z-14z^3-7z^5-z^7)$$
and therefore by Corollary 6.8 
the braid index of the satellite is at least equal to
 $\frac{1}{2}(7-(-5))+1 = 7$, hence the satellite 
does not have a braid representation  with 6 strands. 
Thus we conclude that $b(9_{42}) = 4$. 

It was conjectured that the inequality from Corollary 6.8 becomes 
an equality for alternating links. 
However, K. Murasugi with the author 
found counterexamples \cite{M-P-2}. The simplest link we found has 15 crossings,
Fig. 6.3(a) and the simplest knot 18 crossings Fig.6.3(b). We challenge 
the reader to prove (or disprove) that these are the smallest 
counterexamples.\\
\ \\
\centerline{{\psfig{figure=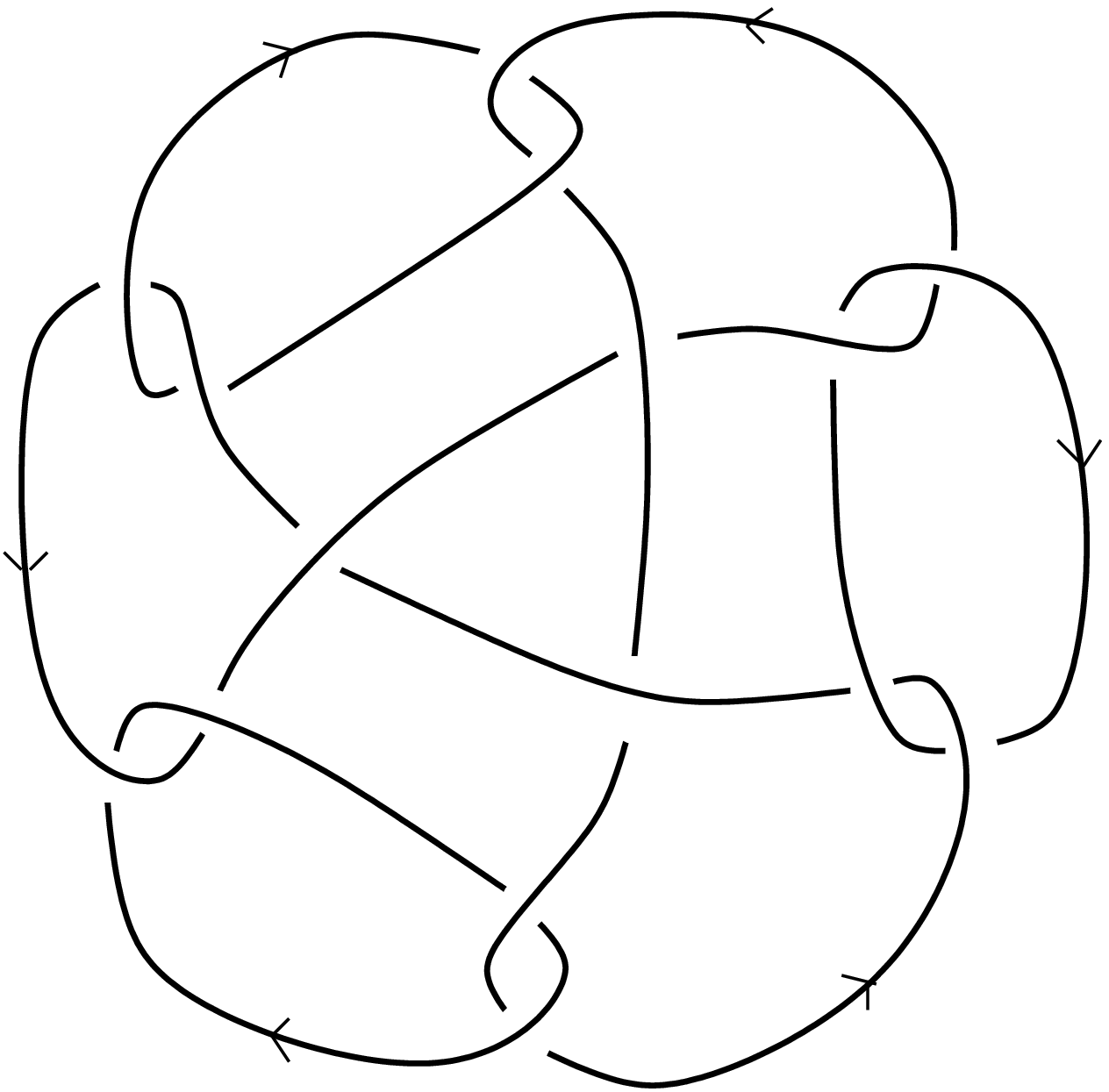,height=4.6cm}}\ \ \ \ \ \ \ \ 
{\psfig{figure=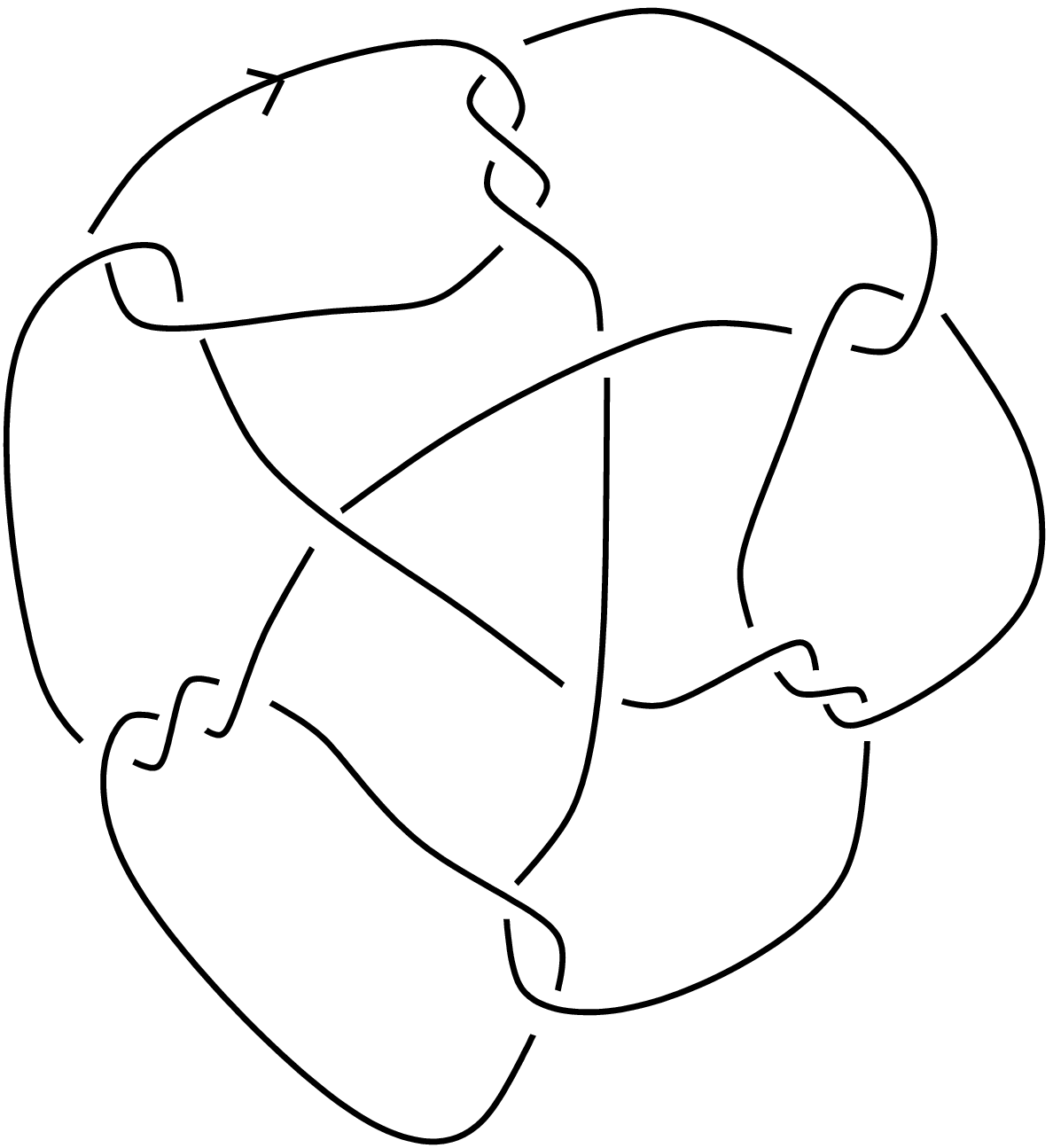,height=4.6cm}} }
\centerline{ \ (a)\ \ \ \ \ \ \ \ \ \ \ \ \ \ \ \ \ \ \ \ \ \ \ \ \ \ \ \ 
\ \ \ \ \ \ \ (b)   \ }
\begin{center}
Fig.~6.3
\end{center}

One could think that Corollary 6.8 is weaker than Corollary 6.6,
since it might seem natural to expect the existence of a link 
whose braid index is greater than the number of Seifert circles in 
some diagram of the link. This expectation, however, was proven to be wrong 
by S.Yamada.
Let $s_{\min} (L)$
%Let $s_{\min (L)}$ (????chyba lepiej $s_{\min} (L)$????)
denote the minimal number of Seifert circles
of possible diagrams of a given link $L$.

\begin{theorem}[Yamada.]\label{4:5.9}
$$s_{\min (L)} = b(L).$$
\end{theorem}

The above theorem was proved for $s_{\min}\leq 7$ by Murakami and
Nakanishi; the general version was proved by Yamada in \cite{Ya}.

Furthermore, Yamada construction applied to a diagram $L$ which realizes 
$s_{min}(L)$ produces a braid with the same Tait number as $L$. 
We present Vogel-Traczyk prove of Yamada theorem in APPENDIX.

\begin{exercise}\label{4:5.10}
Prove that for any oriented diagram $L$ the following inequalities hold 
$$n(L)\geq s(L)-\mu '(L)\geq s(L)-\mu(L),$$
where $\mu '(L)$ is the number of connected components of diagram $L$
(i.e. the number of components of the projection of the diagram, that 
is the graph obtained from $L$ by identifying at every crossing the 
overcrossing and undercrossing).
\end{exercise}

The number $n(L) - (s(L)-1)$ which bounds $M$ from above
 in Theorem \ref{4:5.5} 
is equal to $1-\chi(F_L)$, where $\chi(F_L)$ is the Euler characteristic
of a Seifert surface $F_L$ which is built using
the Seifert circles of $L$.
If $g(F_L)$ is the genus of $F_L$ then $1-\chi(F_L) = 2g(F_L)+\mu(L)-1$.
Let $g_p(L)$ be the minimal genus of a Seifert surface of $L$ 
built using the Seifert circles of some diagram of $L$.
The number $g_p(L)$ is called planar genus of $L$.
The Theorem \ref{4:5.5} provides a bound $M\leq 2g_p(L) +\mu(L)-1$. 

We note, however, 
that the planar genus of $L$ can not be replaced by the genus $g(L)$ of
$L$ (we recall that $g(L)$ is the minimal genus of a Seifert surface of $L$).
For example, the untwisted Whitehead double of a trefoil knot
(Fig.~6.3)

%\vspace*{2in}\centerline{\Psfig{figure=Rys.5.3}}
\centerline{{\psfig{figure=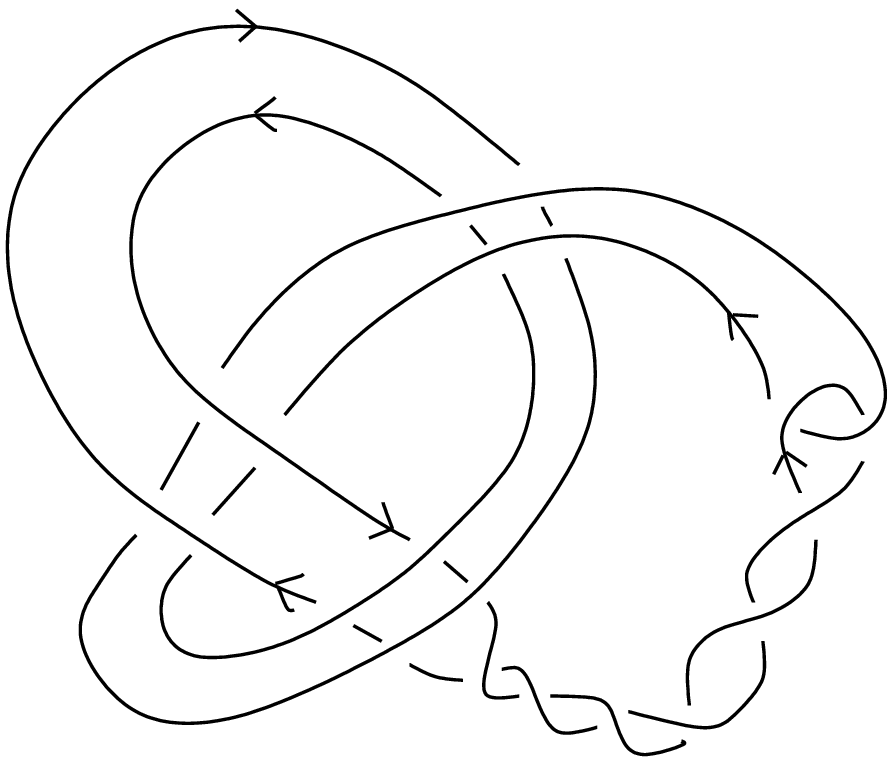,height=3.5cm}}}
\begin{center}
Fig.~6.3
\end{center}
has $M=6$ while $g=1$ (which is true for any Whitehead double
of a non-trivial knot).

On the other hand it is known that the degree of Alexander polynomial 
(and thus also of its version equal to $P(i,z)$)
of the Whitehead double  does not exceed $2g+\mu-1$. Therefore $M$
can be substantially bigger than the degree of $P(i,z)$. 

Theorems 6.4 and  6.5 imply that
$$E+M\leq n(L)-\tilde{n}(L)$$
and 
$$-e+M\leq n(L)+\tilde{n}(L).$$
Now it is natural to ask when the above inequalities become equalities.

%\begin{theorem}\label{V.6.11}
\begin{theorem}\label{4:5.11}
\begin{enumerate}
\item[(i)]
The polynomial $P_L(a,z)$ contains the term $a^Ez^M$, with 
$$E = -\tilde{n}(L) + (s(L)-1),\ \ M=n(L)-(s(L)-1)$$
 if and only if  $n(L)=\tilde{n}(L)$, that is $L$ is a positive diagram 
(all crossings are positive).
\item[(ii)] The polynomial $P_L(a,z)$ contains the term $a^ez^M$, with
$$e= -\tilde{n}(L) - (s(L)-1),\ \ M=n(L)-(s(L)-1)$$
 if and only if  $n(L)=-\tilde{n}(L)$, that is $L$ is a negative 
diagram (all crossings are negative).
\item[(iii)] A positive diagram is always $+$-adequate while
a negative diagram is $-$-adequate.
\end{enumerate}
\end{theorem}

\begin{corollary}\label{4:5.12}
\begin{enumerate}
\item[(i)]
If $L$ is a positive diagram then 
$$M = -E =n(L)-(s(L)-1)$$ 
\item[(ii)] A positive diagram represents 
a trivial link 
if and only if it can be reduced to a trivial diagram by a finite number of
positive first Reidemeister 
moves (\ \parbox{0.8cm}{\psfig{figure=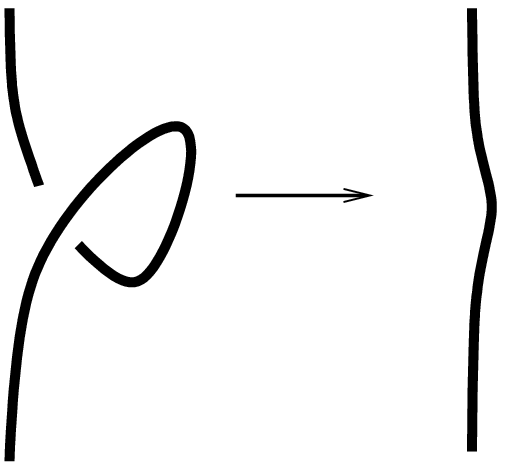,height=0.6cm}}).
\end{enumerate}
\end{corollary}

Proof of Corollary 6.12. Part (i) follows from Theorems 6.11, 6.4 
and 6.5. For the second part we assume that $L$ is a positive 
diagram of a trivial link. Then $1-\mu = M= n(L)-(s(L)-1)$ and the 
conclusion of (ii) follows. To clarify this point it is useful to 
notice that in our situation the Seifert graph $G_{\vec{s}}(L)$ is 
a positive forest so $L$ can be reduced to a trivial link by 
positive first Reidemeister moves.  

\begin{corollary}[Murasugi, Traczyk]\label{4:5.13}
\begin{enumerate}
\item[(i)]
A non-trivial link cannot admit both negative and positive diagram.
\item[(ii)]
If a non-trivial link allows a positive diagram then it is not 
amphicheiral.
\end{enumerate}
\end{corollary}

\begin{proof} Corollary 6.12 implies that, if a nontrivial link allows
a positive connected diagram then $E<0$, and if it admits a connected negative
diagram then $e>0$, and of course $e\leq E$. Let $L$ be a non-connected 
positive (resp. negative) diagram of $k$ components. Then the 
corresponding link is a split link of $k$ components (Corollary 5.13)
 so if a link 
allows positive and a negative diagram then each non-split component of the 
link allows a positive and a negative diagram and (ii) follows in full 
generality. (i) follows from (ii) by the definition of amphicheirality.
\end{proof}
Corollary 6.13 has been also proven (for knots) by T. Cochran and E. Gompf
 who 
noticed that a positive diagram $D$ of a nontrivial knot dominates the 
right-handed trefoil knot (one changes only some positive crossings), 
so it has a negative  signature ($\sigma(D) \leq \sigma(\bar 3_1)=-2$).
Therefore a notrivial positive knot is not amphicheiral, \cite{Co-Go}; 
compare also \cite{Ru,T-9,P-33}. We discuss generalization of a Cochran method 
by K.Taniyama \cite{Tan-1,Tan-2}
 and generalization of Corollary 6.12(ii) and Corollary 6.13 in 
Theorem 6.27.

\begin{corollary}\label{4:5.14}
A positive diagram of a link has a minimal number of crossings
if and only if it has a minimal number of Seifert circles.
In particular, if $L$ is a positive diagram of a link and 
$$s(L) = \frac{E-e}{2} +1$$ then $L$ has the minimal number
of crossings and $$n(L) = M + \frac{E-e}{2}$$. 
\end{corollary}
\begin{proof} 
The first part of the corollary follows from Corollary 6.12
 (i.e.~$n(L) = M + (s(L)-1)$ ); in the second part we apply
also Corollary 5.6 (i.e.~$s(L)-1\geq\frac{E-e}{2}$).
\end{proof}

{\bf Proof of Theorem 6.11}. 
We model our proof on Thistlethwaite approach to Kauffman polynomial 
(section 5).  Let us consider the polynomial 
$\hat{P}_L(a,z) = a^{Tait (L)}P_L(a,z)$, which is
a regular isotopy invariant. 
The polynomial $\hat{P}_L(a,z)$ can be defined by the following conditions:

\begin{enumerate}

\item $\hat{P}_{T_k}(a,z) = a^{Tait (T_k)}(\frac{a+a^{-1}}{z})^{k-1}$, 
where $T_k$ is
any diagram of the trivial link of $k$ components.

\item $\hat{P}_{L_+}(a,z) + \hat{P}_{L_-}(a,z) = z\hat{P}_{L_0}(a,z)$. 
\end{enumerate}

Let us write $\hat{P}_L(a,z) = \sum e_{ij}a^i z^j$. 

\begin{lemma}\label{4:5.15}
\ \\
\begin{enumerate}
%Let $L$ be a connected diagram with $n(L)>0$
\item[(1)] If $e_{ij}\neq 0$ 
then $|i|+j\leq n(L)$.

\item[(2)] If $e_{ij}\neq 0$ and $L$ is $+$-inadequate diagram then
$i+j<n(L)$. 

\item[(3)] If $e_{ij}\neq 0$ and $L$ is $-$-inadequate diagram then
$-i+j < n(L)$.
\end{enumerate}
\end{lemma}

\begin{proof}
For the part (1) we apply standard induction as in the proof of Theorem
\ref{4:3.12}. The parts (2) and (3) are proven in the same manner as
 Lemma 5.8 (using Exercise 5.5 in place of Lemma 5.3).
\end{proof}
To proceed with the proof of Theorem \ref{4:5.11} we introduce
exterior polynomials 
$\psi^+ =\sum_{i+j=n}e_{ij}a^i$ and $\psi^-=\sum_{-i+j=n}
e_{ij}a^i$. 
%(????wyglada na misprint w granicach sumowania???)

If $L$ is $+$ adequate then $L^\sigma_e$ is inadequate
(Corollary 5.9), consequently for a $+$-adequate diagram $L$
it follows that 
$\psi^+_L (a,z) =\psi^+_{L_0} (a,z)$. 
Thus $\psi^+_L$ is non-zero if we can change (reduce) $L$ to a $+$-adequate
diagram $L'$ which has only nugatory crossings, and the change was 
achieved by consecutive smoothings of crossings of $L$ 
in such a way that all intermediate diagrams are adequate.

We see immediately that
\begin{formulla}\label{4:5.16}
$$\psi^+_L = \left\{
\begin{array}{ll}
0&\mbox{if the above reduction is not possible}\\
a^{s(L)-1}&\mbox{if the above reduction is possible}\\
\end{array}
\right.$$
\end{formulla}
%(???moze czytelniej bedzie tak:
%We see immediately that if $\psi^+_L$ is non-zero then it is equal to 
%$a^{s(L)-1}$. bo 
%po co robic tabelke i jeszcze pisac "possible" itd ---
%co jest mniej czytelne ????)

We proved already that if $\psi^+_L neq 0$ then $L$ is a 
$+$-adequate diagram. We will show now that $L$ is also a 
 positive diagram by analyzing Formula 6.16. 
First we associate to the diagram $L$ a 2-color signed
planar graph $G(L)$. That is, the edges of $G(L)$ are either black ($b$)
or white ($w$) as defined in Section 1 and the sign of the edge is
equal to the sign of the respective crossing.
(let us note that not all 2 color signed graphs
are associated to links, e.g.~the graph
$\stackrel{w^-}{\bullet\!\!\!\!\!\longrightarrow\!\!\!\!\!\bullet}$ 
is not related to any diagram; Exercise 6.19); compare Fig.~6.4.

%\vspace*{1.5in}\centerline{\Psfig{figure=Rys.5.4}}
\centerline{{\psfig{figure=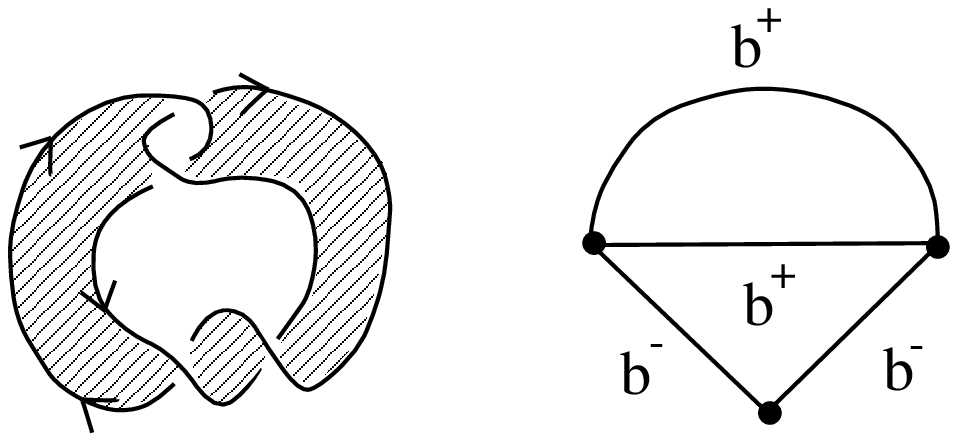,height=2.5cm}}}  
\begin{center}
Fig.~6.4
\end{center}

%\begin{center}
%Fig.~6.4
%\end{center}

Now we use the notation introduced in Section 5. That is, $G_b$ denotes
a subgraph consisted of black edges and $\overline{G_b}$
is the graph obtained from $G$ removing all white edges and identifying 
their endpoints.
%by removing white loops and collapsing of the remaining white edges.
Furthermore, $G_{b^+}$ (respectively $G_{b^-}$) 
is a subgraph of $G$ consisting of positive (resp.~negative) black edges.
Recall (Lemma \ref{4:4.7} (ii)), that $L$ is a $+$ adequate diagram
if and only if  $G_b$ has no isthmus and  $\overline{G}_w$ has no loop.
The smoothing of a crossing is related to either collapsing of $b^+$
or $w^-$, or removing of $b^-$ or $w^+$; see Fig.~6.5. 
%(???chyba $u$ zastapilo $w$, a moze jest ``uhite'' zamiast ``white''?? \ :--)   

%\centerline{{\psfig{figure=Rys4-5-5.eps,height=7.5cm}}}  
\centerline{{\psfig{figure=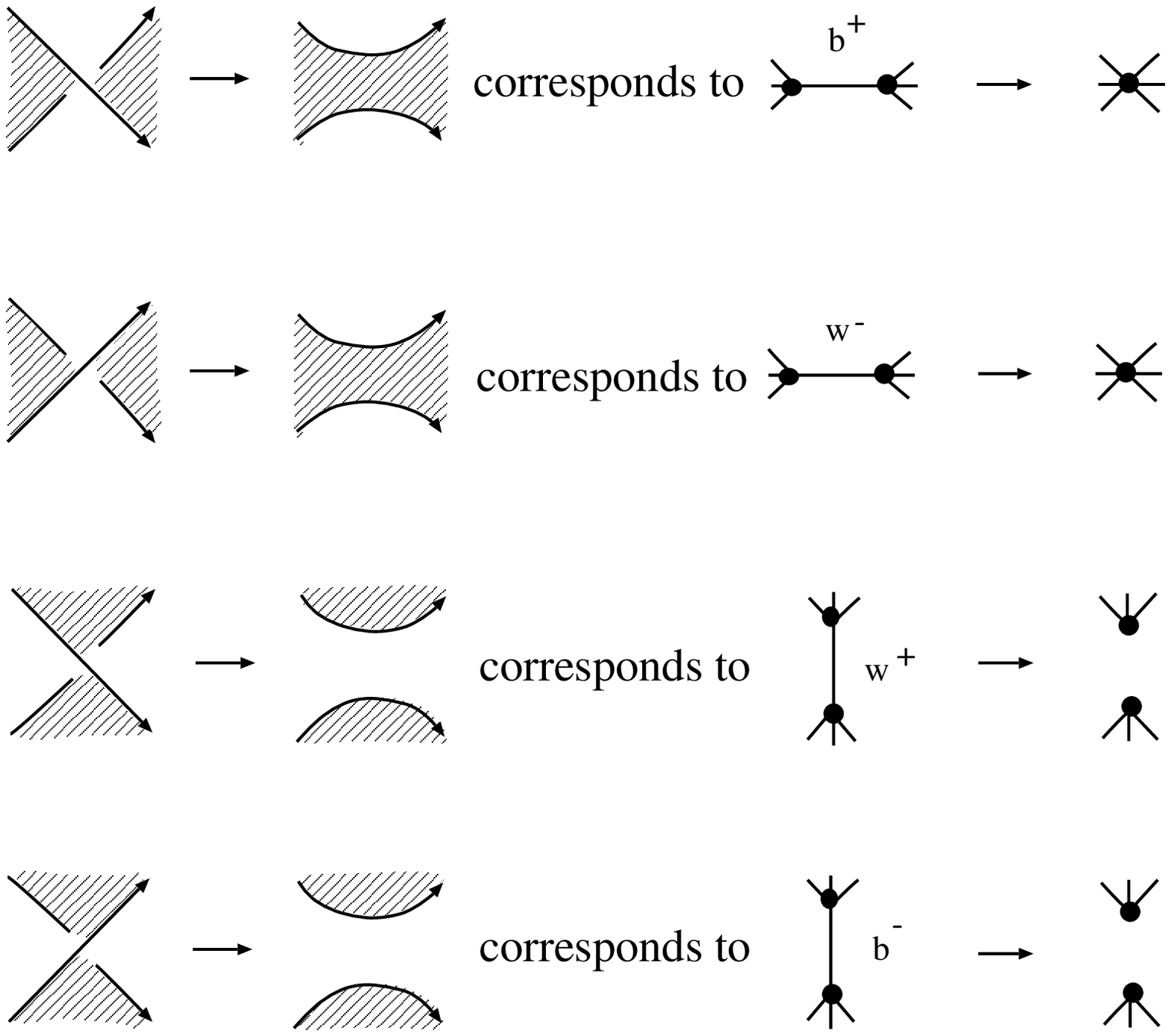,height=8.5cm}}}
\begin{center}
Fig.~6.5
%Fig.~6.5 (odpowiada=corresponds to Old 4.5.5)
\end{center}

%\begin{center}
%Fig.~6.5
%\end{center}

The formula 6.16 is true in a general situation (for any $+$-adequate diagram)
 but for the sake of convenience let us assume that $L$ is a connected 
diagram and on the way from $L$ to a diagram with only nugatory crossing 
we meet only only connected diagrams
(i.e.~we do not remove isthmuses and we do not collapse loops).
Then 6.16 implies that if $L$ can be modified to $L'$
(with all nugatory crossings) and the intermediate diagrams are $+$ 
adequate (so $\psi^+(L)\neq 0$,
then this property does not depend on the order we perform smoothings
(i.e.~if we change the order then the intermediate diagrams are $+$ adequate).
That is the case because if we meet an $+$-inadequate diagram on the way 
the result would be that $\psi^+(L)=0$. 

Now let us analyze the process of smoothing of crossings of $G$.
We begin by collapsing all positive black edges ($b^+$) except
$b^+$. Notice that the resulting graph contains no negative black edges 
($b^-$) because in the process of removing it we would obtain 
a black graph with an isthmus and thus not $+$-adequate graph.
Similarly,  $\overline{G}_w$ has no negative edges. 
Indeed, let $G'$ be a graph obtained by ``eliminating" black edges 
(collapsing $b_+$) except
$b^+$ loops, then
$\overline{G}_w= G_w'$.

Now from $G_w '$ remove $w^+$ edges,
then the resulting diagram is still $+$-adequate
and it consists of $w^-$ edges, which we can contract still keeping 
adequate diagram till we reach white loop which cannot happen in 
$+$-adequate diagram. We should stress that no graph coming from 
an oriented diagram can have $w^-$ isthmus (Exercise 5.19).

Therefore we have proved that if $\psi^+_L\neq 0$ then $L$ 
is a positive diagram.

\begin{lemma}\label{4:5.17}
If $L$ is a positive diagram then $\psi^+_L\neq 0$ and in particular
$L$ is $+$-adequate.
\end{lemma}

\begin{proof} The claim follows essentially from the previous argument,
however, we can prove it in another way. 
Namely, let us note that $P_L (a,a+a^{-1}) = 1$ (Lemma III.3.38(ii)), 
and thus for a positive diagram $L$ we have
$\hat{P}_L(a,a+a^{-1}) = a^{n(L)}$ hence $\psi^+_L\neq 0$.
\end{proof}
This completes the proof of Theorem \ref{4:5.11}.

\begin{exercise}\label{4:5.18}
Prove Theorem \ref{4:5.11} directly, that is, without using
(or, at least, partially eliminating) 
Thistlethwaite's method. In particular, do not use adequate diagrams.
\end{exercise}
To illustrate methods which I have in mind in Exercise 6.18 let us 
prove part of Lemma 6.17 that a positive diagram is $+$-adequate:\ \ 
For an oriented diagram $L$ we define the Seifert state $\vec s$ as 
a Kauffman state in which all markers agree with orientation of $L$.
For a Kauffman state $s$ we say that $s$ is adequate if $sD$ has no 
circles touching itself (compare Chapter X were we construct a graph 
$G_s(L)$ for any state $s$ and $s$ is adequate if $G_s(L)$ has no loop).
Thus $L$ is $+$-adequate if $s_+$ is adequate and and $L$ is 
$-$-adequate if $s_-$ is adequate. We have 
\begin{proposition}\label{V.6.19}
\begin{enumerate}
\item[(i)] $\vec s$ is always an adequate state.
\item[(ii)] $\vec s= s_+$ iff $L$ is a positive diagram.
\end{enumerate}
\end{proposition}
\begin{proof} (i) A Seifert circle cannot touch itself, otherwise 
Seifert surface would be unorientable.\\
(ii) It is illustrated in Fig. 6.6.
\end{proof}
\centerline{{\psfig{figure=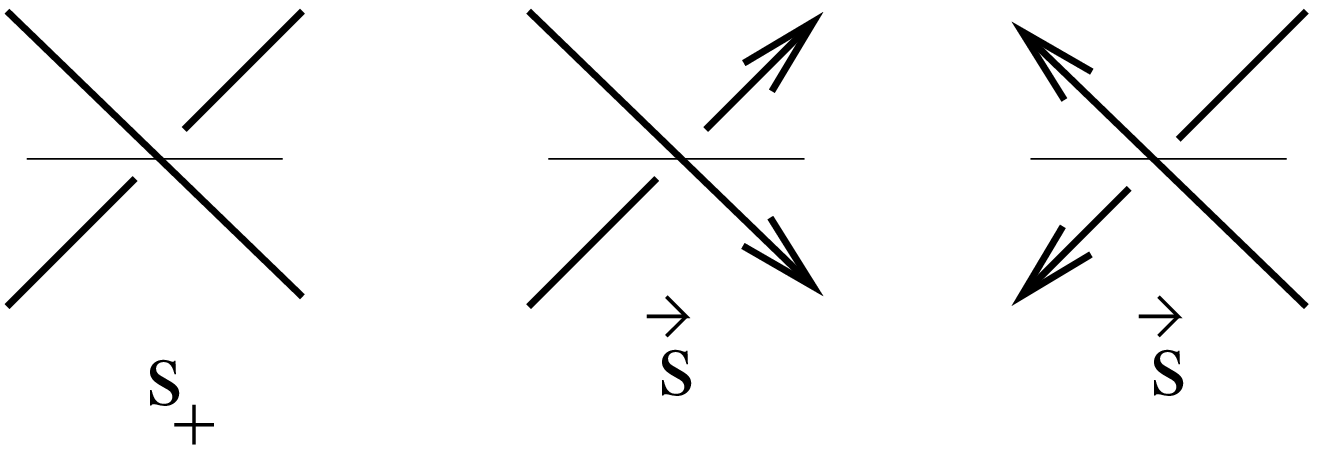,height=3.5cm}}}
\begin{center}
Fig.~6.6
\end{center}
 
%\begin{exercise}[\cite{Ko}]\label{V.6.20} 
\begin{exercise}[\cite{Ko}]\label{4:5.19}
Find necessary conditions which have to be satisfied by a 
2-color signed
planar graphs associated to an oriented links.
In particular, show that\\
(i) No isthmus can be of type $w^-$ or $b^+$
while no loop can be of type $w^+$ or $b^-$.\\
(ii) Show that any vertex of the graph has an even number of
$b^+$ and $w^-$ edges. Formulate a ``dual" statement for 
$b^-$ and $w^+$ edges. 
\end{exercise}

An interesting problem concerns finding conditions on diagrams
of oriented links for which the inequalities from Theorem
\ref{4:5.5} become equalities, that is, when $M= n(L)-(s(L)-1)$. 
We dealt with this problem (with additional assumptions) in Theorem 
\ref{4:5.11}. Some other partial results were obtained by  Kobayashi
\cite{Ko} and Traczyk \cite{T-3}. Traczyk proved that if $L$ is 
an alternating diagram
of a fibered link, that is a link whose complement in $S^3$ is
 a (Seifert) surface-bundle over a circle then
$M=n(L)-(s(L)-1)$ and the coefficient $b_M(a)$ with the highest degree
monomial $z^M$ is equal to $a^{\sigma(L)}$, where  $\sigma(L)$  
is the Trotter-Murasugi signature of $L$ (see Remark 6.23). 

Below, we present a slight generalization of Traczyk's result.

\begin{definition}\label{V.6.21} \  
\begin{enumerate}
\item
An oriented diagram $L$ of a link is called simplified tree-like diagram 
if it satisfies the following conditions:
\begin{enumerate}
\item For any pair of Seifert circles of the diagram $L$
all crossings ``connecting" this pair have the same sign (equivalently 
multiple-edges of the Seifert graph have the same sign).
\item Let $\Gamma (L)$ be a graph associated to $L$ 
with vertices representing Seifert circles of $L$ and
the vertices joined by a (single) edge if the respective
Seifert circles touches in the diagram  
$L$ (c.f.~Fig.~6.7)\footnote{In   
other words $\Gamma (L)$ is obtained from the Seifert graph by 
replacing every multiple edge by a singular edge.}.
Then $\Gamma (L)$ is a tree.

%\vspace*{1.8in}\centerline{\Psfig{figure=Rys.5.6}}
\centerline{{\psfig{figure=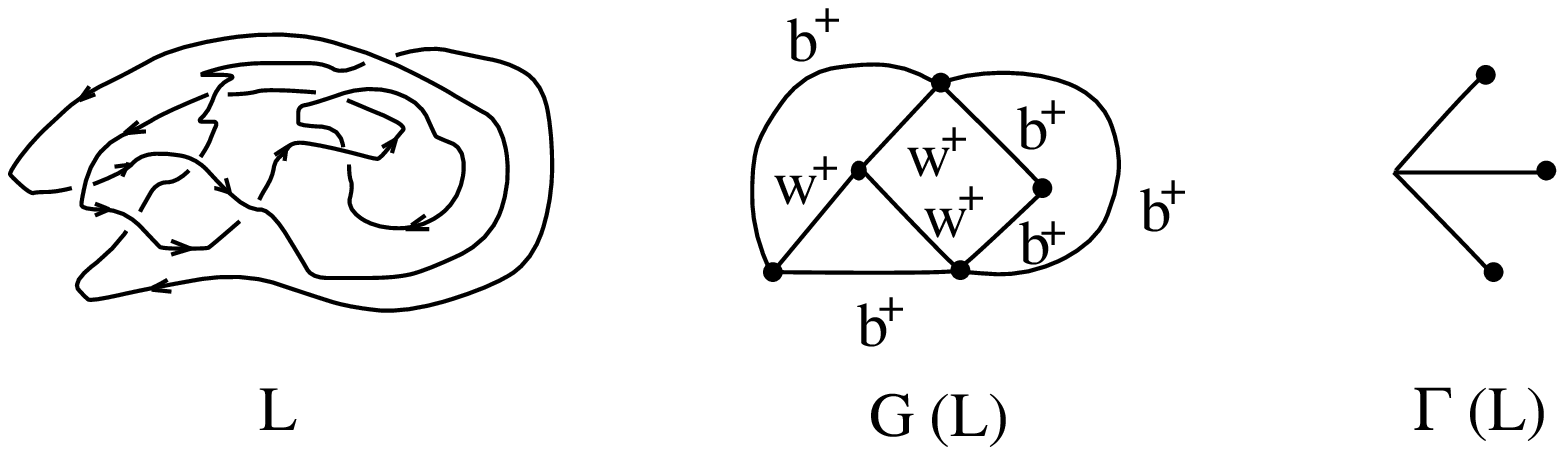,height=3.5cm}}}  
\begin{center}
Fig.~6.7
\end{center}
\end{enumerate}

\item An oriented diagram of a link is called a tree-like diagram
if it can be obtained from a simplified tree-like diagram by 
replacing any crossing with an odd number
of crossings (half-twists), as it is illustrated in Fig. 6.8
(note that the sign of the crossing is preserved).\\
\ \\
\centerline{\psfig{figure=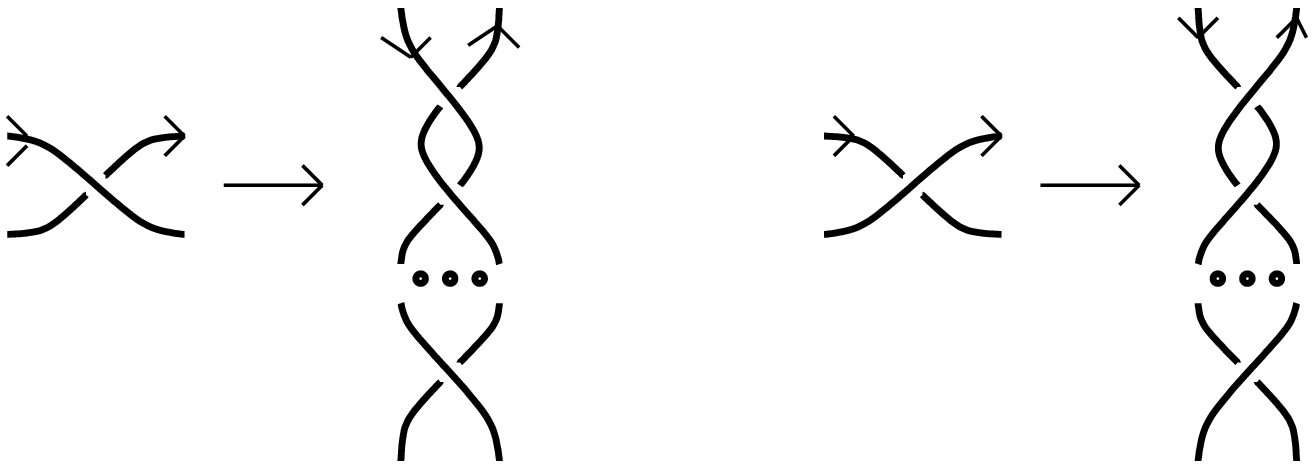,height=3.1cm}}
\ \\ \ \\
%\centerline{\psfig{figure=Rys4-5-7.eps,height=3.5cm}}  
\centerline{\psfig{figure=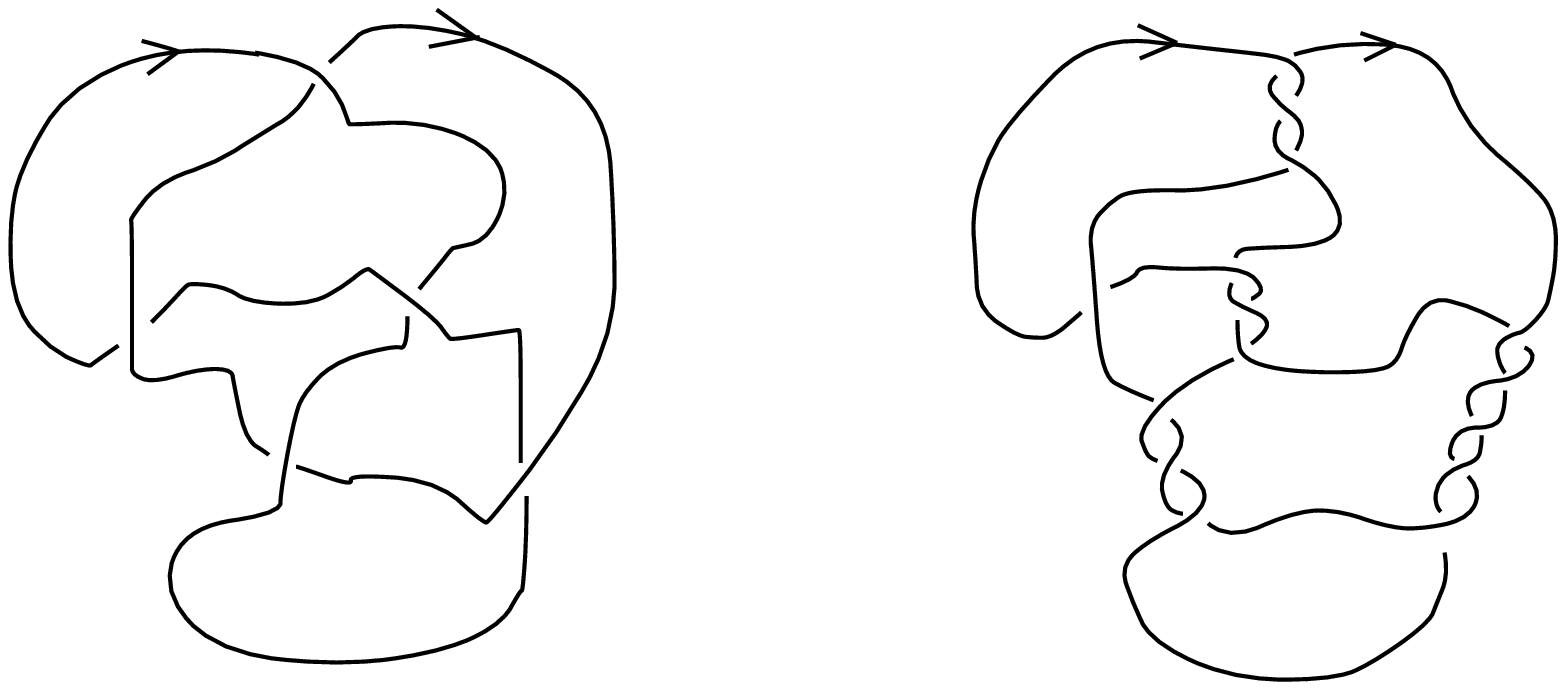,height=3.5cm}}
\begin{center}
Fig.~6.8
\end{center}

\end{enumerate}
\end{definition}

%\begin{theorem}\label{V.6.22}
\begin{theorem}\label{4:5.21}
If $L$ is a tree-like diagram of an oriented link then
\begin{enumerate}
\item[(1)] $M= n(L)-(s(L)-1)$

\item[(2)] The coefficient $b_M(a)$ of the monomial $z^M$ in $P_L(a,z)$ 
is an alternating (Laurent) polynomial of the variable $a$
and moreover $\deg_{\min} b_M(a)$ and $\deg_{\max} b_M(a)$ 
can be computed in the following way:
let any two Seifert circles $C_i, C_j$ of a simplified
tree-like diagram (related to $L$) meet in $k(i,j)$ crossings
which in the tree-like diagram $L$ are replaced 
by $d_1(i,j),d_2(i,j),\ldots,d_k(i,j)$ crossings, respectively, then 

\begin{eqnarray*}
\deg_{\max} b_M(a)&=&-\tilde{n}(L)+\sum_{i,j;d_r(i,j)>0} 
((\sum_r (d_r(i,j)-1))+1) - \sum_{i,j;d_r(i,j)<0} m(i,j)\\
\deg_{\min} b_M(a)&=&-\tilde{n}(L)+\sum_{i,j;d_r(i,j)<0} 
((\sum_r (d_r(i,j)+1))-1) + \sum_{i,j; d_r(i,j)>0} m(i,j)\\
\end{eqnarray*}
where $m(i,j) = \min_r |d_r(i,j)|$.

\item[(3)] 
If moreover $L$ is a simplified tree-like diagram then
$$b_M(a) = a^{-\tilde{n}(L) + d^+(L) -d^-(L)},$$ 
where $d^+(L)$ (respectively $d^-(L)$) is the number of pairs of Seifert 
circles joined by positive (respectively, negative) crossings.
If moreover $L$ is alternating then $b_M(a)= a^{\sigma(L)}$. 
\end{enumerate}
\end{theorem}

\begin{remark}\label{4:5.22}
In \cite{M-3} Murasugi proved that simplified tree-like alternating
diagrams correspond exactly to alternating fiber links.
% complement in $S^3$ of which
%has a structure of a bundle over $S^1$ (they are called fiber links) 
He showed also that all simplified tree-like diagrams 
represent fiber links (the inverse remains open).
\end{remark}

Proof of Theorem \ref{4:5.21} (1) and (3). Consider a diagram $L$
which differs from a simplified tree-like diagram at some bridge $B$.
That is after changing some crossings at $B$ we get 
a simplified tree-like diagram. Furthermore, one of the crossings at $B$
 which have to be changed is non-nugatory.
Then:
%\begin{formulla}\label{V.6.24}
\begin{formulla}\label{4:5.23}
$$M<n(L) - s(L)+1.$$
\end{formulla}

If the diagram $L$ satisfies above conditions then we define
the complexity of $L$ to be equal to an ordered pair 
$(n(L), k(L))$, where $k(L)$ is the number of crossings which 
have to be changed (i.e.~undercrosing to overcrossing) in order to
obtain a diagram for which the claim is true and which has the longest 
bridge $B'$ containing  $B$ (similarly as in the proof of Lemma 5.3)
Let us consider the lexicographic order of pairs $(n,k)$. 
We will use induction with respect to the complexity.
If $n=2$ then the link $L$ is trivial with the diagram
\parbox{0.8cm}{\psfig{figure=Tr-H.eps,height=0.5cm}},
so \ref{4:5.23} holds.

If $k=0$ then $B$ is of maximal length already
and therefore it is reduced by Reidemeister moves to a diagram
which contains $n(B)$ fewer crossings ($B$ is a closed component or
it contains a loop whose shrinking can eliminate crossings).
Furthermore, the number of Seifert circle can decrease at most 
by $n(B)-1$. We use here the fact that $B$ contains a non-nugatory crossing
 hence there is a 
Seifert circle which meets $B$ at  two crossings  at least
(Seifert circles which do not meet $B$ are unchanged).
Since the resulting link is ambient isotopic to the old one, then by 
Morton's inequality (Theorem \ref{4:5.5}) we get $M(L)<n(L)-(s(L)-1)$. 

Next we consider $L$ which satisfies our assumptions and which has
complexity $(n,k)$, where $n>2$, $k>0$, 
and we assume that the inequality \ref{4:5.23} holds for diagrams
with smaller complexity.
Let $p$ be a crossing at one of the ends of $B$, the change of which will
decrease the complexity of $L$. The diagram obtained by changing $p$ 
will be denoted $L'$. Because of the inductive assumption
$M(L')<n(L)-(s(L)-1)$ and since 
$a^{\mp 1}P_L +a^{\pm 1}P_{L'} = zP_{L_{0}}$
we will be done if we prove that
$$M(L_{0})< n(L_{0}) - (s(L_{0})-1).$$

If $L_{0}$ satisfies assumptions of \ref{4:5.23} then we 
conclude by inductive assumption. Otherwise 
%when after smoothing of $p$ 
all crossings of $L'$ on $B$ which have to be changed in order to get  
a simplified tree-like diagram are nugatory.
This, however, means that the two Seifert circles of $L$ meeting at $p$ 
meet at only one crossing different from $p$, let us denote
 that crossing by $x$.  The crossing $x$ is on $B$ and it is 
the only crossing on $B$ which have to be changed in order to make $L$  
a simplified tree-like diagram.
Therefore $x$ and $p$ have opposite signs (Fig.~6.9) which means that 
$p$ was not an obstruction to extend the bridge $B$. Hence $p$ did not have
to be changed. This completes the proof of 
 \ref{4:5.23}.\\
 \ \\

%\vspace*{1.5in}\centerline{\Psfig{figure=Rys.5.8}} 
\centerline{{\psfig{figure=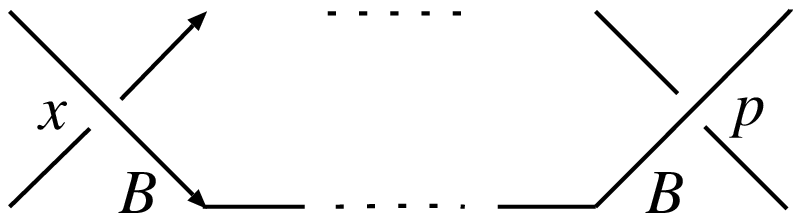,height=2.5cm}}}  
\begin{center}
Fig.~6.9
\end{center}

Now we can conclude the proof of Theorem \ref{4:5.21} (3).

One could use induction with respect to the number of crossings 
however, we prefer to consider global properties 
of the diagram. 
Smoothing of a  non-nugatory crossing leads to 
a simplified tree-like diagram while changing the crossing gives
a diagram with Jones-Conway polynomial containing only 
``small'' degree monomials of $z$ (by \ref{4:5.23}. 
Therefore, to get the maximal exponent of $z$ in Jones-Conway polynomial
we have to smooth the maximal number of crossings
so that the diagram remains connected.
If Seifert circles $C_i$ and $C_j$ were joined by 
$k(i,j)$ crossings then we have to smooth $k(i,j)-1$ of them.
The smoothing of a positive crossing brings the factor $a^{-1}z$, 
while the smoothing of a negative crossing contributes 
the factor $az$ to the Jones-Conway polynomial
 Therefore, in the resulting
polynomial we get a monomial
$$a^{-\tilde{n}(L)+d^+(L)-d^-(L)}z^{n(L)-(s(L)-1)},$$ as required.
Traczyk proved (see Chapter IV and\cite{T-2}) 
that if $L$ is an alternating link then 
$$\sigma(L) = -\tilde{n}(L)+d^+(L)-d^-(L),$$  which explains the last 
statement of Theorem 6.22(3).

Now we prove Theorem 6.22(2).

We start with the following simple formula (compare 
\cite{P-2})

\begin{formulla}\label{4:5.24}\
\begin{enumerate}
\item[(i)] 
$P_{{\parbox{2.9cm}{\psfig{figure=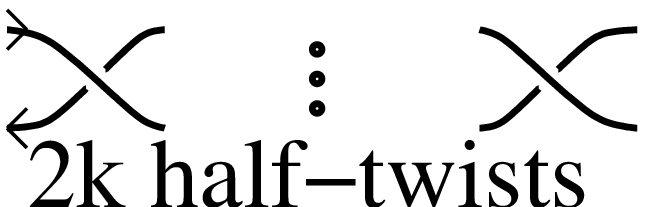,height=0.8cm}}}}
(a,z)= z(a-a^3+a^5 + ... +a(-a^2)^{k-1})
P_{{\parbox{1.2cm}{\psfig{figure=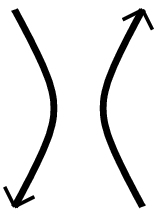,height=0.6cm}}}}+$

\ \ $(-a^2)^{k} 
P_{{\parbox{1.2cm}{\psfig{figure=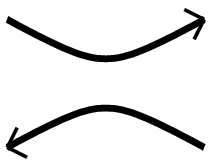,height=0.6cm}}}}
= za^k\frac{a^{-k} +(-1)^{k+1}a^k}{a+a^{-1}}
P_{{\parbox{1.2cm}{\psfig{figure=vert-anti.eps,height=0.6cm}}}}+ 
(-a^2)^{k}
P_{{\parbox{1.2cm}{\psfig{figure=horiz-anti.eps,height=0.6cm}}}}$
\item[(ii)]
$P_{{\parbox{2.9cm}{\psfig{figure=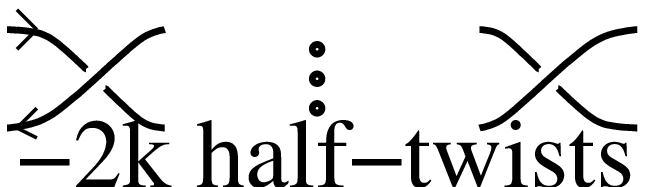,height=0.8cm}}}}
(a,z)= z(a^{-1} -a^{-3}+ a^{-5}+ ... +a^{-1}(-a^{-2})^{k-1})
P_{{\parbox{1.2cm}{\psfig{figure=vert-anti.eps,height=0.6cm}}}}+$

\ \ $(-a^{-2})^{k}
P_{{\parbox{1.2cm}{\psfig{figure=horiz-anti.eps,height=0.6cm}}}}
= za^{-k}\frac{a^{k} +(-1)^{k+1}a^{-k}}{a+a^{-1}}
P_{{\parbox{1.2cm}{\psfig{figure=vert-anti.eps,height=0.6cm}}}}+
(-a^{-2})^{k}
P_{{\parbox{1.2cm}{\psfig{figure=horiz-anti.eps,height=0.6cm}}}}.$

\end{enumerate}
\end{formulla}

We consider Seifert circles  $C_i, C_j$ joined by  $k(i,j)$
bands, each one of them twisted $d_r(i,j)$ times ($1\leq r\leq
k(i,j)$) creating $d_r(i,j)$ positive crossings, $d_r(i,j)$ an odd number. 
Let us consider the case when $d_r(i,j)>0$. 

%Let $\alpha_r = a(a^{-1}-a^{-3}+a^{-5}+\cdots
%+a^{-1}(-a^{-2})^{\frac{1}{2}(d_r(i,j)-1)}$ and 
% $\beta_r := a(-a^{-2})^{\frac{1}{2}(d_r(i,j)-1)}$. 
Let $\alpha_r = 1 + (-a^2)^{-1} + (-a^2)^{-2}+ \cdots 
(-a^2)^{-(\frac{1}{2}(d_r(i,j)-1)-1)}$ and \\
$\beta_r = (-a^{2})^{-(\frac{1}{2}(d_r(i,j)-1))}$. In this notation 
we get:
$$P_{{\parbox{3.5cm}{\psfig{figure=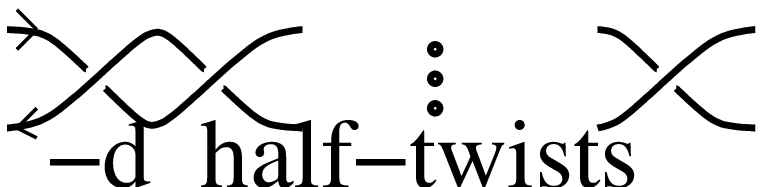,height=0.8cm}}}}
(a,z)= za^{-1}\alpha_r
P_{{\parbox{0.9cm}{\psfig{figure=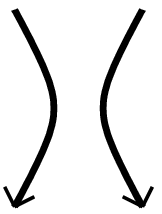,height=0.6cm}}}}+
\beta_r
P_{{\parbox{1.2cm}{\psfig{figure=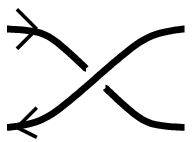,height=0.6cm}}}}
$$

Consider links $L_i$ and $L_j$ being split components of the link 
obtained from $L$ by removing
 all bands which join $C_i$ and $C_j$ (Fig. 6.10). \\
 \ \\
\centerline{{\psfig{figure=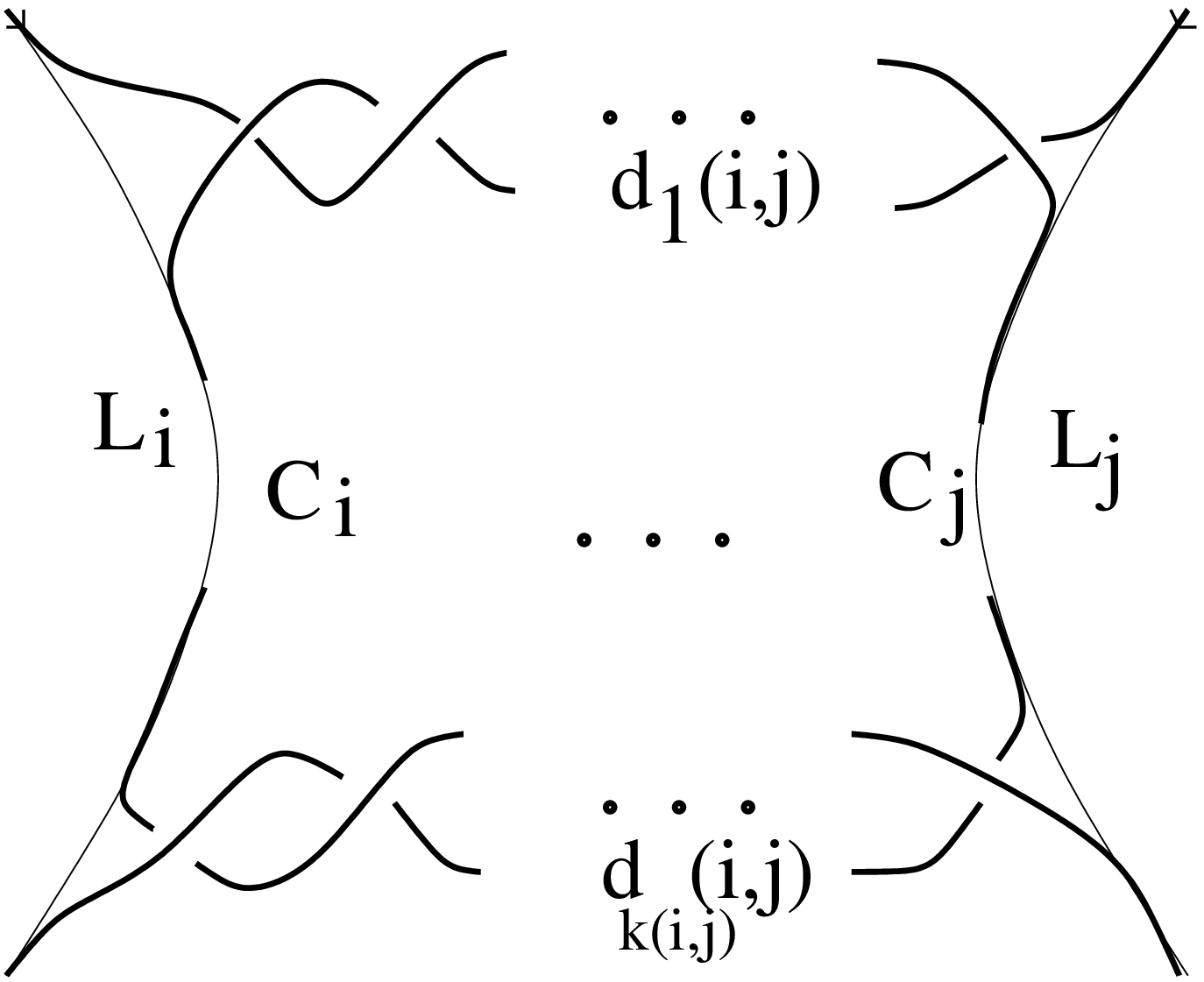,height=5.9cm}}}
\begin{center}
Fig.~6.10
\end{center}

Then
\begin{eqnarray*}
%&P_L(a,z) &=\\
%&z^{k(i,j)-1)} P_{L_i}(a,z)\cdot P_{L_j}(a,z)&\cdot\\
&b_{M(L)}(a)z^{M(L)}&=\\
&z^{k(i,j)-1)}b_{M(L_i)}(a)z^{M(L_i)}\cdot b_{M(L_j)}(a)z^{M(L_j)}&\cdot\\
&a^{1-k(i,j)}
((\alpha_1 +\beta_1)(\alpha_2+\beta_2)\cdots(\alpha_{k(i,j)}+\beta_{k(i,j)})+
(a+a^{-1}-1)
\alpha_1\alpha_2\cdots\alpha_{k(i,j)})&=\\
&z^{k(i,j)-1)}b_{M(L_i)}(a)z^{M(L_i)}\cdot b_{M(L_j)}(a)z^{M(L_j)}&\cdot\\
&a^{1-k(i,j)}
((\alpha_1 +\beta_1)(\alpha_2+\beta_2)\cdots(\alpha_{k(i,j)}+\beta_{k(i,j)})+
(1-\alpha_{k(i,j)}-\beta_{k(i,j)})
\alpha_1\alpha_2\cdots\alpha_{k(i,j)-1}).&
\end{eqnarray*}

The above formula, which is a simple consequence
of \ref{4:5.24}  and Theorem \ref{4:5.21}(3), 
allows us to compute the coefficient $b_M(a)$ of $z^M$
in Jones-Conway polynomial. In particular, we can find out
that $b_M(a)$ is alternating as substituting $b=-a^2$ gives, up to 
an invertible monomial, a polynomial with positive coefficients. We  
we can also compute 
$\deg_{\max} b_M(a)$ and $\deg_{\min} b_M(a)$. 
This completes the proof of Theorem \ref{4:5.21}.

\begin{exercise}\label{V.6.26}
Let us consider the following generalization of a tree-like diagram.
A diagram $L$ of a link is called generalized tree-like diagram if
it can be obtained from a simplified tree-like diagram by replacing
of any half-twist by an odd number of half-twists (now we allow the 
sign of twisting to be changed). Find conditions when 
$M(L) = n(L) - (s(L)-1)$ for a generalized tree-like diagram. 
For example if $k(i,j)=2$ and $d_1(i,j)=-d_2(i,j)$ then 
$M(L) <n(L) - (s(L)-1)$.
%in particular, if $k(i,j)=2$, then $b_1(i,j)\neq  -b_2(i,j)$. YES
%In fact one can use isotopy for reduction.
\end{exercise}

\subsection{Almost positive links}\label{V.6.1}
It was asked by Birman and Williams \cite{Bi-Wi} and L.Rudolph whether 
nontrivial Lorenz knots have always positive signature. Lorenz knots 
are examples of positive braids\footnote{In older conventions positive 
braids had all crossing negative, in this book positive braids are 
defined to have positive crossings.}. It was shown by Rudolph \cite{Ru} 
that positive braids have positive signature (if they represent 
nontrivial links). Murasugi has shown that nontrivial, alternating, 
positive links have negative signature. Cochran and  Gompf proved 
that a nontrivial positive knot has negative signature 
\cite{Co-Go,T-9,P-33}. From this it followed that a notrivial positive 
knot is not amphicheiral. Another proof, by Murasugi and Traczyk, 
was given in  Corollary 6.13. 
%we proved that a notrivial positive link is not amphicheiral. 
We conjectured (\cite{P-33}, Conjecture 5) that 
 if we allow one negative crossing in $D$ (i.e. D is almost positive) 
then the link is not amphicheiral as well. The conclusion of the conjecture 
followed easily from the master thesis of K. Taniyama that almost positive 
nontrivial link dominates the right handed trefoil knot or the positive 
Hopf link \cite{Tan-1,Tan-2}. Therefore it has negative signature and 
cannot be amphicheiral. We made several generalizations of this result 
in \cite{P-Ta} as described below.

A link is {\it $m$-almost positive} if it has a
diagram with all but $m$ of its crossings being positive.\\
The {\it unknotting number} (Gordian number) of a positive link
is equal to $\frac{1}{2}(c(D)-s(D)+com(D))$, where $D$ is a positive
diagram of the link, $c(D)$ is the number of crossings, $s(D)$ is the
number of {\it Seifert circles} of $D$, and $com(D)$ is the number
of components of the link (this generalizes the {\it Milnor's
unknotting conjecture}\footnote{The unknotting number of a $(p,q)$ torus 
knot is equal to $\frac{(p-1))(q-1)}{2}$.}, 
1969 and the {\it Bennequin conjecture}, 1981.
Furthermore for a positive knot the
unknotting number is equal to the {\it 4-ball genus} of the knot\footnote{If 
$S^3=\partial D^4$ and $K$ is a knot in $S^3$ then the 4-ball genus is 
the minimal genus of a surface in $D^4$ bounding $K$},
to the  genus of the knot,
to the {\it planar genus} of the knot (from Seifert construction),
to the the minimal degree of the Jones polynomial and
to the half of the degree of the Alexander polynomial \cite{Kr-Mr}.
An elementary proof of the formula for an unknotting number of positive 
knots, using Khovanov homology, was given in 2004 by 
J. Rasmussen \cite{Ras} (compare Chapter X).

One can define a relation $\geq $ on links by $L_1 \geq L_2$ iff $L_2$
can be obtained from $L_1$ by changing some positive crossings of $L_1$.
This relation allows us to express several fundamental properties of
positive (and $m$-almost positive) links.
\begin{theorem}
\begin{enumerate}
\item[(1)] 
If $K$ is a positive knot then $K \geq (5,2)$ positive torus knot
unless $K$ is a {\it connected sum} of {\it pretzel knots} $L(p_1,p_2,p_3)$,
where $p_1$, $p_2$ and $p_3$ are positive odd numbers
\begin{enumerate}
\item[(a)]
If $K$ is a nontrivial positive knot
then either the signature $\sigma(K)\leq-4$ or $K$ is a pretzel
knot $L(p_1,p_2,p_3)$ (and then $\sigma(K)=-2$).
\item[(b)] If a positive knot has unknotting number one then it is a
positive twist knot.
\end{enumerate}
\item[(2)]
Let $L$ be a
nontrivial 1-almost positive link. Then $L\ge$ right-handed trefoil
knot (plus trivial components), or $L\ge$ right-handed Hopf link
(plus trivial components). In particular $L$ has a negative signature.
\item[(3)] If $K$ is a $2$-almost positive knot then either
\begin{enumerate}
\item[(i)] $K \geq $ right handed trefoil, or
\item[(ii)] $K \geq \bar 6_2$ (mirror image of $6_2$ knot)
($\sigma_1^3\sigma_2^{-1}\sigma_1\sigma_2^{-1}$ in the {\it braid} notation)
 or
\item[(iii)] $K$ is a twist knot with a negative clasp.
\item[(a)] If $K$ is a $2$-almost positive knot different from a twist
  knot with a negative clasp then $K$
has negative signature and $K(1/n)$ (i.e. $1/n$ surgery on
  $K$, $n>0$) is a homology 3-sphere that does not bound a compact, smooth
  homology $4$-ball, \cite{Co-Go,P-Ta}.
\item[(b)] If $K$ is a non-trivial  2-almost positive knot different from
  the {\it Stevedore's knot} then $K$ is not a slice knot.
\item[(c)] If $K$ is a non-trivial  2-almost positive knot different from
  the {\it figure eight knot} then $K$ is not amphicheiral.
\end{enumerate}
\item[(4)]
Let $K$ be a 3-almost positive
knot. Then either $K\geq$ trivial knot or $K$ is the left-handed
trefoil knot (plus positive knots as connected summands). In
particular, either $K$ has a non-positive signature or
$K$ is the left-handed trefoil knot.
\end{enumerate}
\end{theorem}

\ \\ \ \\ \ \\
\noindent \textsc{Dept. of Mathematics, Old Main Bldg., 1922 F St. NW \\
The George Washington University, Washington, DC 20052}\\
e-mail: {\tt przytyck@gwu.edu}

\end{document}